\documentclass[conf]{new-aiaa} 
\usepackage[utf8]{inputenc}

\usepackage[version=4]{mhchem}
\usepackage{longtable,tabularx}
\usepackage{graphicx}
\usepackage{bbold}
\usepackage{algorithm}
\usepackage{algpseudocode}
\usepackage{diagbox}
\usepackage{booktabs}

\usepackage{todonotes}

\makeatletter
\newcommand*\bigcdot{\mathpalette\bigcdot@{.4}}
\newcommand*\bigcdot@[2]{\mathbin{\vcenter{\hbox{\scalebox{#2}{$\m@th#1\bullet$}}}}}
\makeatother

\newcommand{\autoscvx}{\textsc{\tt Auto-SCvx}}
\newcommand{\ptr}{\textsc{\tt PTR}}

\newcommand{\boxing}[6]{
	\begin{figure*}[#1]
		\noindent\makebox[\textwidth][c] {
			\fbox{
				\begin{minipage}{#4}
					\vspace{0.2cm}
					\begin{#2} {\bf \textit{#5}} \label{#3} \end{#2}
					#6
				\end{minipage}
			}
		}
	\end{figure*}
}

\newcommand{\lin}[2]{{\ell_{#1}(\Delta #2, \bar{#2})}}
\newcommand{\lintwo}[3]{{\ell_{#1}([\Delta #2^{\top}, \Delta #3^{\top}]^{\top}, [\bar{#2}^{\top}, \bar{#3}^{\top}]^{\top})}}
\newcommand{\linfull}[2]{{\ell_{#1}(#2 - \bar{#2}, \bar{#2})}}
\newcommand{\argmin}[1]{ \operatorname{arg} \underset{#1}{\operatorname{min.}}}
\newcommand{\prox}[2]{{\frac{1}{2 #1}\|#2-\bar{#2}\|_{2}^{2}}}
\newcommand{\proxdev}[2]{{\frac{1}{2 #1}\|\Delta #2\|_{2}^{2}}}

\newcommand{\hn}{{h_{\text{ncvx}}}}
\newcommand{\gn}{{g_{\text{ncvx}}}}
\newcommand{\hc}{{h_{\text{cvx}}}}
\newcommand{\gc}{{g_{\text{cvx}}}}
\newcommand{\nequ}{n_{\text{eq}}}
\newcommand{\nineq}{n_{\text{ineq}}}
\newcommand{\nH}{n_{\text{H}}}
\newcommand{\nG}{n_{\text{G}}}
\newcommand{\Mh}{M_{\text{h}}}
\newcommand{\Mg}{M_{\text{g}}}
\newcommand{\MH}{M_{\text{H}}}
\newcommand{\MG}{M_{\text{G}}}
\newcommand{\one}[1]{\mathbb{1}_{#1}}
\newcommand{\zero}[1]{\mathbb{0}_{#1}}
\newcommand{\diag}[1]{{\operatorname{diag}\left(#1\right)}}
\newcommand{\Wh}{{W_{h}}}
\newcommand{\Wg}{{W_{g}}}

\newcommand{\wh}[1]{{\omega_{h_{#1}}}}
\newcommand{\wg}[1]{{\omega_{g_{#1}}}}
\newcommand{\epsh}[1]{{\varepsilon_{h_{#1}}}}
\newcommand{\epsg}[1]{{\varepsilon_{g_{#1}}}}
\newcommand{\epsopt}{{\varepsilon_{\text{opt}}}}
\newcommand{\epsoptz}{\epsopt_{z}}
\newcommand{\epsoptJ}{\epsopt_{J}}
\newcommand{\epsoptin}[1]{\epsopt_{#1}}
\newcommand{\epsfeas}{{\varepsilon_{\text{feas}}}}
\newcommand{\epsfeasin}[1]{\epsfeas_{#1}}

\newcommand{\nz}{n_z}
\newcommand{\nx}{n_x}
\newcommand{\nnu}{n_u}
\newcommand{\dz}{\Delta z}
\newcommand{\dx}{\Delta x}
\newcommand{\dT}{\Delta T}
\newcommand{\dxhat}{\Delta \hat{x}}
\newcommand{\du}{\Delta u}
\newcommand{\dpp}{\Delta p}
\newcommand{\dq}{\Delta q}
\newcommand{\dlam}{\Delta \lambda}
\newcommand{\dmu}{\Delta \mu}
\newcommand{\dnu}{\Delta \nu}
\newcommand{\zstar}{z^{*}}

\newcommand{\pstar}{p^{*}}
\newcommand{\qstar}{q^{*}}

\newcommand{\dzstar}{\dz^{*}}

\newcommand{\dlamstar}{\dlam^{*}}
\newcommand{\dmustar}{\dmu^{*}}

\newcommand{\xprop}{x^{\text{prop}}}
\newcommand{\zbar}{{\bar{z}}}
\newcommand{\xbar}{{\bar{x}}}
\newcommand{\ubar}{{\bar{u}}}
\newcommand{\Tbar}{\bar{T}}
\newcommand{\pbar}{{\bar{p}}}
\newcommand{\qbar}{{\bar{q}}}
\newcommand{\lambar}{{\bar{\lambda}}}
\newcommand{\mubar}{{\bar{\mu}}}
\newcommand{\nubar}{{\bar{\nu}}}
\newcommand{\lamhat}{{\hat{\lambda}}}
\newcommand{\muhat}{{\hat{\mu}}}
\newcommand{\forallnodes}{\forall k \in \{1,\dots,N\}}
\newcommand{\forallintervals}{\forall k \in \{1,\dots,N-1\}}
\newcommand{\Ak}{A_{k}}
\newcommand{\Bmk}{B^{-}_{k}}
\newcommand{\Bpk}{B^{+}_{k}}
\newcommand{\Sk}{S_{k}}
\newcommand{\PhiA}{\Phi_{A}}
\newcommand{\PhiBm}{\Phi_{B}^{-}}
\newcommand{\PhiBp}{\Phi_{B}^{+}}
\newcommand{\PhiS}{\Phi_{S}}
\DeclareMathOperator{\Lagr}{{\mathcal{L}}}
\DeclareMathOperator{\Lncvx}{\Lagr_{\text{ncvx}}}

\renewcommand{\And}{\textbf{and}}
\newcommand{\Or}{\textbf{or}}
\newenvironment{method}[1][H]
{
\floatname{algorithm}{Algorithm}
\begin{algorithm}[#1]%
\small
}{\end{algorithm}}

\newcommand{\Xk}{\bar{\tilde{X}}_{k}}
\newcommand{\eval}[2]{\left. #1 \right\rvert_{#2}}
\newcommand{\evalfull}[2]{\left. \left[ #1 \right] \right\rvert_{#2}}

\newcommand{\dotc}[1]{\overset{\circ}{#1}}

\newcommand{\real}{\mathbb{R}}

\newcommand{\infnorm}[1]{\left\|#1\right\|_{\infty}}

\title{\Large 
Auto-tuned  Primal-dual Successive Convexification for \\
Hypersonic Reentry Guidance
\footnote{Based on worked presented at the AIAA SciTech Conference 2025 in Orlando, FL.} 
}

\author{
Skye Mceowen\footnote{Ph.D.\ Candidate, UW Dept. of Aero. and Astro., \texttt{skye95@uw.edu}, AIAA Member.}, 
Daniel J. Calderone\footnote{Postdoctoral researcher, UW Dept. of Elec. and Comp. Eng., \texttt{djc@uw.edu}.},
Aman Tiwary\footnote{M.S.\ Student UW Dept. of Mech. Eng., \texttt{amant10@uw.edu}, AIAA Member.}, 
Jason S. K. Zhou\footnote{M.S.\ Student, UW Dept. of Aero. and Astro., \texttt{jskzhou@uw.edu}, AIAA Member.}, 
Taewan Kim\footnote{Ph.D.\ Candidate, UW Dept. of Aero. and Astro., \texttt{twankim@uw.edu}, AIAA Member.}, 
Purnanand Elango\footnote{Ph.D.\ Candidate, UW Dept. of Aero. and Astro., \texttt{pelango@uw.edu}, AIAA Member.}
and
Beh\c{c}et A\c{c}\i{}kme\c{s}e\footnote{Professor, UW Dept. of Aero. and Astro., \texttt{behcet@uw.edu}, AIAA  Fellow.}}
\affil{\vspace{0.75em} University of Washington, Seattle, WA 98105, USA}

\begin{document}

\renewcommand{\thetable}{\arabic{table}}

\maketitle

\begin{abstract}
This paper presents \textit{auto-tuned primal-dual successive convexification} (\autoscvx), an algorithm designed to reliably achieve dynamically-feasible trajectory solutions for constrained hypersonic reentry optimal control problems across a large mission parameter space. In {\autoscvx}, we solve a sequence of convex subproblems until convergence to a solution of the original nonconvex problem. This method iteratively optimizes dual variables in closed-form in order to update the penalty hyperparameters used in the primal variable updates. A benefit of this method is that it is auto-tuning, and requires no hand-tuning by the user with respect to the constraint penalty weights. Several example hypersonic reentry problems are posed and solved using this method, and comparative studies are conducted against current methods. In these numerical studies, our algorithm demonstrates equal and often improved performance while not requiring hand-tuning of penalty hyperparameters. 
\end{abstract}

\section{Introduction}

Hypersonic reentry is an increasingly relevant application for trajectory optimization. Flight vehicle reusability is gradually becoming the new standard for rocket and spacecraft design, especially in the context of human-rated missions. Methodologies that can reliably design trajectories satisfying a restrictive set of mission constraints while maximizing vehicle performance are a growing need in the spaceflight and defense industries. In addition, an increasing number of more maneuverable, high lift-to-drag vehicles are in development. The combination of these factors makes the development of trajectory design algorithms for hypersonic reentry vehicles an increasingly critical field. Early methods for vehicle reentry trajectory design were offline algorithms that tracked the reference drag-acceleration profiles with classical feedback control \cite{harpold1979shuttle}. These techniques were extended to include energy-dependent drag profiles \cite{roenneke1994re}, nonlinear tracking control schemes \cite{mease1994shuttle}, and drag profile optimization \cite{lu1997entry}. Determining feasible trajectories from drag profiles was a focus \cite{leavitt2007}. However, these offline profiles are not robust in the face of mission changes or anomalies. In response, predictor-corrector methods emerged in an effort to recompute desired drag profiles online \cite{rea2007,lu2008predictor}, where a bank angle trajectory is recomputed numerically by integrating the dynamics subject to desired boundary conditions. However, these techniques are limited by their inability to handle inequality constraints and modeling errors. To handle inequality constraints in the context of predictor-corrector methods, quasi-equilibrium glide conditions are derived to create a unified framework \cite{shen2003onboard,xue2010constrained,lu2014entry}. Stochastic predictor-corrector guidance algorithms have also been proposed \cite{mcmahon2021}. 

Online entry guidance algorithms have been flown for highly challenging or constrained missions, such as an extended Apollo entry guidance algorithm used for Mars Science Laboratory (MSL), to minimize range error while ensuring parachute deployment at sufficiently safe altitudes for a variety of mission configurations \cite{msl-entry2011}. Other proposed approaches include an adaptive disturbance-based sliding mode controller to improve tracking in the presence of parameter dispersions and modeling errors \cite{sagliano2017}. Pseudospectral methods have been widely applied to compute feasible solutions for hypersonic reentry flight. One proposed approach determines trajectories onboard using adaptive multivariate pseudospectral interpolation \cite{sagliano2016}. Other approaches find solutions based on using pseudospectral techniques to solve optimal control problems for hypersonic flight \cite{ross2012,malyuta2021advances,sagliano2018,jorris2012,fahroo2008}. Numerical methods for optimal control are attractive because they allow for explicit specification of both the objective to be optimized and analytical constraints that must remain satisfied during flight \cite{rao2009,patterson2014gpops}.

Trajectory planning through convex optimization has demonstrated real-time performance capabilities in a wide variety of applications in recent years including rocket landing \cite{LiuSurvey2017,kamath2023customized}. In this framework, an optimal control problem is formulated that specifies a mission performance objective, as well as a constraints restricting the vehicle dynamics, state, and control variables to a feasible region.  The speed and guarantees of convex optimization to perform real-time, online trajectory optimization are been central to modern-day precision landing applications for both suborbital and orbital rockets \cite{behcet2007jgcd,Blackmore2010,blackmore2016autonomous}. The first real-time rocket landing algorithm used to compute an online trajectory harnessed the lossless convexification algorithm, which guarantees global optimality for the original problem despite relaxing nonconvex control constraints into a convex form. More general nonconvex optimal control problems for rocket landing, with arbitrary nonconvexities in the dynamics, state and control constraints have been successfully solved with sequential convex programming (SCP) techniques by iteratively solving a series of convex subproblems until convergence to an approximate local optimal solution of the original problem \cite{mao2016cdc,mao2017aut}. In these techniques, convex approximations of the original nonconvex problem are modeled about a reference trajectory; this is either the initial guess or the solution to a previous iteration. As each approximate subproblem is solved, the algorithm walks towards an optimal solution in the feasible region of all constraints.

Successive convexification is a variant of SCP that is widely used for trajectory optimization \cite{szmuk2016scitech}. Two flavors of the successive convexification algorithm have been presented; one method employs hard trust region constraints to preserve the validity of the local convex approximations, and another method employs soft constraints to penalize trust regions in the cost function (\ptr{}) \cite{malyuta2021convex}. The {\ptr} algorithm has been most widely used in practice. This method augments the nonconvex constraints with virtual buffer variables, relaxing these nonconvex constraints to allow violation in the early iterations of the algorithm. A penalty function is added to the subproblem cost which heavily penalizes these virtual buffer variables having nonzero magnitude. If implemented properly, this makes infeasibility more expensive than the true cost, driving the virtual buffers to zero and incentivizing convergence to the feasible region of the nonconvex constraints \cite{szmuk2018successive}. In addition to being successfully demonstrated on a wide range of 6-degree-of-freedom (DoF) rocket landing problems \cite{szmuk2018jgcdarxiv,Reynolds2019}, and real-time aerial drone flight experiments  \cite{szmuk2019iros,mceowen2022aero}, these methods have also been increasingly relevant for hypersonic reentry. All of these approaches have historically used hand-tuned, constant hyperparameter virtual buffer penalty weights of large positive magnitude. 

SCP techniques applied to reentry guidance applications \cite{wang_paper,liu2016, 2021scvxentry} have often modeled subproblems as second-order cone programming problems \cite{LiuSurvey2017} or applied pseudo-spectral transcription methods to model the vehicle dynamics \cite{wang2019rapid,yu2019,sagliano-6dof-entry-scp}. Such algorithms aim to achieve fast solution speeds and accurate trajectory guidance and control that adhere to relevant constraints \cite{wang2018autonomous,wang2019rapid,wang2020improved,han2019rapid,zhao2017reentry}.  To mitigate high-frequency jitter in the control profiles produced by SCP techniques in the reentry problem, work has been done both to convert the nonlinear reentry dynamics into a control affine system, and to prove that these relaxed control parameterizations converge to a solution of the original problem even in the presence of virtual buffers \cite{bae2022convex}. However, convergence to a feasible solution is often highly sensitive to scaling, a good approximation of the trajectory used to initialize the algorithm, and well-tuned hyperparameters for penalizing the constraint violation. These feasibility hyperparameters have historically been painstakingly hand-tuned via a trial-and-error approach. This is not guaranteed to work, and when infeasibility occurs it may be unclear whether this is due to an ill-posed optimal control problem, scaling and conditioning, or simply a poor guess for the hyperparameters. Although frameworks such as the augmented Lagrangian update linear penalty hyperparameters (modeled as dual variables), the quadratic penalties are still hand-tuned with a scalar weight \cite{bertsekas2014constrained,oguri2023}. In addition, accuracy of the nonconvex constraints, especially the vehicles dynamics, often degrades during discretization creating susceptibility to intersample constraint violation. Dense time grids on the order of multiple hundreds of nodes are often necessary for retaining the dynamic feasibility and integrity of the solutions. 


\textbf{Contributions}. 
In this work, we present \textit{auto-tuned primal-dual successive convexification} (\autoscvx) as a framework for reliably solving nonconvex optimal control problems, such as hypersonic reentry guidance, over a wide problem parameter space with sparse time grids and high accuracy. This method optimizes dual variables in closed-form within the SCP framework in order to update the penalty hyperparameters used in the primal variable update. The main contribution of this method is that it is auto-tuning, and does not require either hand-tuning nor an initial guess with respect to the constraint penalty weights by the user. This algorithm is motivated to enhance reliability and accuracy for hypersonic reentry trajectory optimization solutions across changes in mission parameterization with the following features: (1) closed-form penalty weight updates for the primal subproblem using dual variable information; (2) a deviation variable model for improved scaling of each convex subproblem; (3) an inverse-free exact discretization technique that permits arbitrary satisfaction (up to machine precision) of the dynamics over large time horizons (using multiple shooting); and (4) a virtual buffer penalty-term approach to eliminate artificial infeasibility while constructing each subproblem as a quadratic program (QP) amenable to real-time applications. This approach also allows for a variable time grid that leaves the timesteps between nodes (and thus the time-horizon) as free variables for the optimizer. 

This work extends the control parameterization of a reusable launch vehicle (RLV) model presented in \cite{wang_paper}. In prior work, this model was used to pose and solve a hypersonic reentry trajectory optimization with {\ptr} using bank angle as the control input \cite{mceowen2023scitech3DoF}. In this paper, we present an extension to include angle-of-attack as an additional control input, allowing modulation of the lift-to-drag ratio for maneuvering the vehicle. To remove jitter from the resulting control solution, a control rate limit constraint is applied across each timestep due to the continuous first-order-hold parameterization of the control. 

To demonstrate the performance of the proposed method, a study is performed for trajectory optimization on a variety of hypersonic reentry models formulated as nonconvex optimal control problems. These problems are subject to various 3-DoF dynamic models for hypersonic reentry, multiple control parameterizations, and nonconvex path constraints. Comparative studies against the existing {\ptr} algorithm are conducted. In these studies, {\autoscvx} demonstrates an ability to reliably solve a wide array of complex problems with equal and often improved performance against existing methods. The remainder of the paper is structured as follows. First,  the various models used for 3-DoF hypersonic reentry optimal control problems are discussed. Next, the {\autoscvx} methodology extensions are formulated for a generic optimal control problem, before being tailored to the hypersonic reentry problem. Finally, numerical comparative results of the algorithm for various reentry models are shown, before discussion in the context of other methods and final remarks.

\section{Problem Formulation: 3-DoF Hypersonic Reentry Model}

\subsection{Nondimensionalization} \label{sec:nondim}

The nonlinear system dynamics for modeling a 3-DoF hypersonic reentry vehicle are expressed in general as:
\begin{align}
    \dot{x} = f(t, x(t), u(t)), \label{eq:gen-ncvx-dyn} 
\end{align}
where time $t \in \real$, state vector $x \in \real^{\nx}$, and control input vector $u\in \real^{\nnu}$ have all been nondimensionalized as described in Table \ref{tab:nondim}. Here $g_\oplus=9.81 \text{m}/\text{s}^2$ is the gravitational acceleration at Earth's surface, and $R_\oplus=6378 \cdot 10^3\text{m}$ is the Earth's radius. All corresponding dimensional quantities and constraints in the formulation containing $t$, $x$ and $u$ have been nondimensionalized appropriately following the same procedure \cite{lu2014entry}. This step is critical to ensure proper numerical conditioning of the problem, as the individual elements of the dimensional state and control vectors tend to have large relative magnitude differences in both their values and rates of change. 
Two specific models are presented below using this nondimensionalization technique, one with bank angle as the only control input, and the other where both bank angle and angle-of-attack are consider as controls.

\begingroup
\centering
\begin{table}
\centering
\caption {\label{tab:nondim} \text{Dimensional quantities are scaled to become unitless quantities.}} 
\begin{tabular}{lrrr}
 \hline \hline 
 Dimensional Quantity               & Units                        & Divide By & \\ 
\hline 
     Time                                & $[\text{s}]$                  & $(R_\oplus / g_\oplus)^{\frac{1}{2}}$ \\ 
    Distance                            & $[\text{m}]$                  & $R_\oplus$ \\
    Velocity                            & $[\text{m} / \text{s}]$       & $(R_\oplus g_\oplus)^{\frac{1}{2}}$ \\
    Acceleration                        & $[\text{m} / \text{s}^2]$     & $g_\oplus$ \\
    Angles                        & $[\text{rad}]$     & $1$ \\ \hline \hline
\end{tabular}
\end{table}
\endgroup

\subsection{Reentry Vehicle Dynamics}

The dynamics for an unpowered flight vehicle reentering over a spherical, rotating Earth are modeled as:
\begin{align}
    \dot{x}(t) = \left\{\begin{array}{l} 
                \dot{r}= v \sin \gamma \\
                \dot{\theta}= \frac{v \cos \gamma \sin \psi}{r \cos \phi} \\ 
                \dot{\phi}= \frac{v \cos \gamma \cos \psi}{r}\\
                \dot{v}= - D - \left(\frac{\sin \gamma}{r^{2}}\right)+\Omega^{2} r \cos \phi (\sin \gamma \cos \phi - \cos \gamma \sin \phi \cos \psi ) \\
                \dot{\gamma}= \frac{1}{v}\left(L \cos \sigma + \left(v^{2}-\frac{1}{r}\right)\left(\frac{\cos \gamma}{r}\right)+2 \Omega v \cos \phi \sin \psi
                +\Omega^{2} r \cos \phi(\cos \gamma \cos \phi+\sin \gamma \cos \psi \sin \phi)\right) \\ 
                \dot{\psi}= \frac{1}{v}\left(\frac{L \sin \sigma}{\cos \gamma}+\frac{v^{2}}{r} \cos \gamma \sin \psi \tan \phi-2 \Omega v(\tan \gamma \cos \psi \cos \phi-\sin \phi)
                +\frac{\Omega^{2} r}{\cos \gamma} (\sin \psi \sin \phi \cos \phi)   \right) 
                \end{array}\right. ,
                \label{eq:nl-dyn-spherical}
\end{align}
with states $r$ representing orbital radius of the vehicle, $\theta$ and $\phi$ representing vehicle longitude and latitude (respectively), $v$ representing vehicle velocity, $\gamma$ representing flight path angle, and $\psi$ representing vehicle heading measured clockwise from the north in the local horizontal plane \cite{vinh1980hypersonic,lu2014entry}. These states comprise vector $x = [r,\theta,\psi,v,\gamma,\psi]^\top$. The control input is $u$, which can be parameterized as the bank angle $\sigma$, angle-of-attack $\alpha$, or both. The earth's rotation rate is given as $\Omega \approx 7.292 \cdot 10^{-5} \ \mathrm{rad} /\mathrm{~s}$. All quantities have been nondimensionalized as described in Section \ref{sec:nondim}. Lift and drag accelerations, also nondimensionalized, are given as:
\begin{align}
    L =& \frac{R_\oplus \rho v^2 S_{\text{ref}} C_L }{2 m}, \\
    D =& \frac{R_\oplus \rho v^2 S_{\text{ref}} C_D }{2 m},
\end{align}
which are dependent on atmospheric density: 
\begin{align}
    \rho = \rho_\oplus \exp(-\beta R_\oplus (r - 1 ) ),
\end{align}
where Earth's radius $R_\oplus$ is defined as in Section \ref{sec:nondim} and $\beta=1/H$ is computed assuming atmospheric scale height $H \approx 7000$ m. For this work, a vertical-takeoff, vertical-landing reusable launch vehicle (RLV) model is adapted from \cite{wang2018autonomous}, with assumed reference area $S_{\text{ref}} \approx 391.2 \text{m}^2$ and mass $m \approx 104.3 \cdot 10^3 \text{kg}$.

\subsection{Control Input and Aerodynamic Coefficient Modeling} \label{sec:ctrl-model}

Two control parameterizations are considered in this work. The first model assumes bank angle as the only control input, which modulates the direction of lift about the velocity vector. In this case the angle-of-attack is assigned a pre-designed profile, modeled as a nonlinear function of velocity. The second model is extended to include both bank angle and angle-of-attack as control inputs. The original velocity-dependent angle-of-attack profile is relaxed into an inequality constraint, such that angle-of-attack control is allowed to modulate within $\pm 5^\circ$ of the original design. 
In both cases, the aerodynamic lift coefficient is then determined as a quadratic function of angle-of-attack, and the drag coefficient in turn is determined as a quadratic function of the lift coefficient. For the case where bank angle is the only control input, and the angle-of-attack is a direct function of velocity, the lift and drag coefficients can be reformulated as a function of velocity. This is discussed in detail below.

\subsubsection{Bank Angle Control} \label{subsec:bank-ctrl}

\begin{figure}
    \centering
    \includegraphics[width=1\linewidth,trim={5cm 0.5cm 5cm 0.5cm},clip]{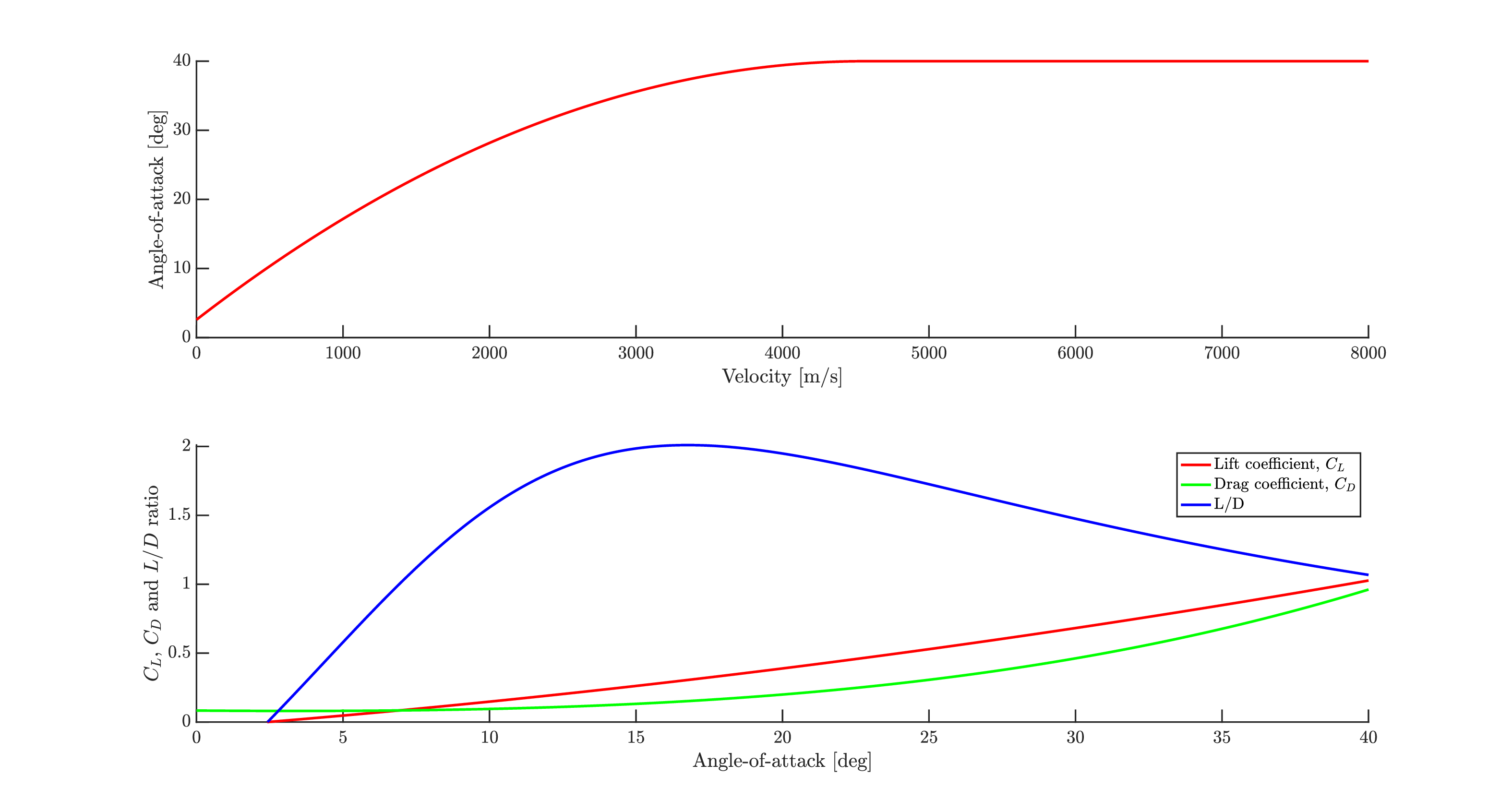}
    \caption{\footnotesize Look-up tables are provided for computing lift and drag coefficients of the RLV model \cite{wang2018autonomous}, using velocity-dependent angle-of-attack profiles in Equation \eqref{eq:wang-alpha-lut}, and associated angle-of-attack-dependent lift-to-drag ratios in Equations \eqref{eq:wang-cl-lut} and \eqref{eq:wang-cd-lut}.}
    \label{fig:LD_vs_alpha}
\end{figure}

The control constraint for the bank angle control is: 
\begin{align}
    C(x,u)= 
        |\sigma|-\sigma_{\max } 
        \leq 0,
\label{eq:bank-cnst}
\end{align}
where $\sigma_{\max} = 80^\circ$ denotes the upper magnitude bound on the of bank angle $u=\sigma$. The control rate constraint is:
\begin{align}
    \mathcal{\dot{U}} \triangleq  \{ \dot{\sigma} \in \real \ | \ \ | \dot{ \sigma } | - \dot{\sigma}_{\max}  \leq 0 \},
\end{align}
where $\dot{\sigma}_{\max} = 5^\circ/\text{s}$ denotes the upper magnitude bound on the bank angle rate. Analytical look-up table functions for computing lift and drag coefficients of the RLV model are provided in \cite{wang2018autonomous}, as shown in Figure \ref{fig:LD_vs_alpha}. First, velocity-dependent angle-of-attack profiles are determined in units of degrees as $\alpha^{\circ}$ (note that this is the only formula in this paper that does not assume an angle given in radians):
\begin{align}
    \alpha^{\circ} &= \left\{\begin{array}{ll}
    K_{\alpha1}, & \text { if } v > V_{\text{lim}} , \\
    K_{\alpha1}-K_{\alpha2}(v - V_{\text{lim}})^{2}, & \text { else }, \\
    \end{array}\right.
    \label{eq:wang-alpha-lut}
\end{align}
assuming $V_{\text{lim}}=4570 \text{m}/\text{s}$, $K_{\alpha1}=40 ^{\circ}$, and $K_{\alpha2}=1.7910 \cdot 10^{-6} \ [^{\circ} \ \mathrm{s}^2 / \mathrm{m}^2]$. For a given velocity, the resulting angle-of-attack (remaining in degrees) can be used in turn to compute the aerodynamic coefficients: 
\begin{align}
    C_L &= K_\text{L1} + K_\text{L2} \alpha^{\circ} + K_\text{L3}\left( \alpha^{\circ} \right)^{2},
    \label{eq:wang-cl-lut}
    \\
    C_D &= K_\text{D1} + K_\text{D2} C_L + K_\text{D3} C_L^{2},
    \label{eq:wang-cd-lut}
\end{align}
where coefficients $K_\text{L1} = -0.041065$, $K_\text{L2} = 0.016292$, $K_\text{L3} = 0.0002602$, $K_\text{D1} = 0.080505$, $K_\text{D2} = -0.03026$, and $K_\text{D3} = 0.86495$. 

\begin{figure}
    \centering
     \includegraphics[width=1\linewidth,trim={5cm 0.5cm 5cm 0.5cm},clip]{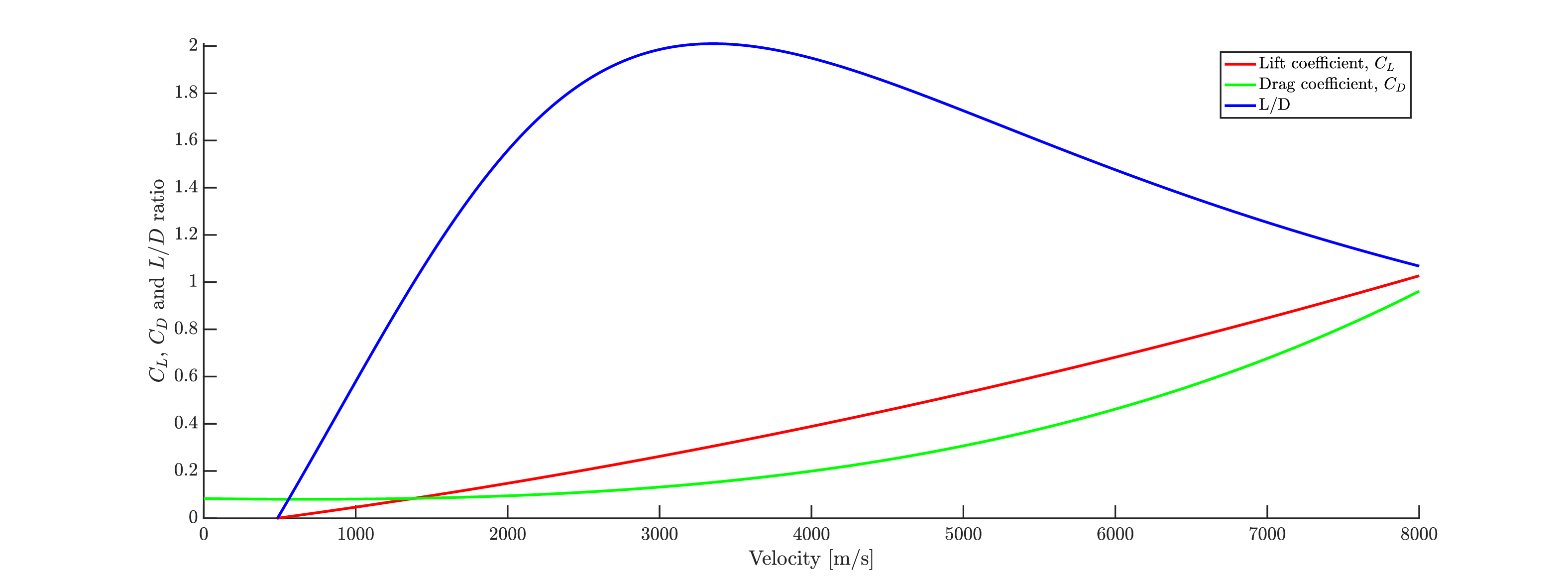}
    \caption{\footnotesize Equivalent lift and drag coefficient models as a direct function of velocity are formulated in Equations \eqref{eq:my-cl-lut} and \eqref{eq:my-cd-lut} for ease of taking analytical convex approximations of the dynamics.}
    \label{fig:LD_vs_v}
\end{figure}

To ease the process of producing convex approximations of the aerodynamics later on in this work--which requires taking analytical Jacobians of the dynamics--we reformulate the lift and drag coefficient look-up tables given in Equations \eqref{eq:wang-alpha-lut}--\eqref{eq:wang-cd-lut} into an explicit function of velocity, as shown in Figure \ref{fig:LD_vs_v}:
\begin{align}
                C_L(v) &= \left\{\begin{array}{ll}
                \bar{K}_{\text{L1}}, 
                & \text { if } v > V_{\text{lim}}\\
                \bar{K}_{\text{L1}}+\bar{K}_{\text{L2}}\left(v-V_{\text {lim }}\right)^{2}+\bar{K}_{\text{L3}}\left(v-V_{\text {lim }}\right)^{4},
                & \text { else }
                \end{array}\right.,
                \label{eq:my-cl-lut}
                \\
                C_D(v) &= \left\{\begin{array}{ll}
                K_{D 1} + \bar{K}_{\text{L1}}(K_{D 2}+K_{D 3} \bar{K}_{\text{L1}}), 
                & \text { if } v>V_{\text {lim }} \\
                K_{D 1}+K_{D 2} C_L+K_{D 3} C_L^{2}, 
                & \text { else } \\
                \end{array}\right.,
                \label{eq:my-cd-lut}
\end{align}
where the updated coefficients $\bar{K}_{\text{L1}} = K_{\text{L1}} + K_{\text{L2}} K_{\alpha1} + K_{\text{L3}} K_{\alpha1}^2$,
$\bar{K}_{\text{L2}}  = -K_{\text{L2}} K_{\alpha1} - 2 K_{\text{L3}} K_{\alpha1} K_{\alpha2}$, and
$\bar{K}_{\text{L3}}  = K_{\text{L3}} K_{\alpha2}^2$. Conveniently, the partial derivatives for the lift and drag coefficients can now be computed explicitly as a piece-wise function of velocity:
\begin{align}
                \frac{\partial C_L}{\partial v} &= \left\{\begin{array}{ll}
                0, 
                & \text { if } v>V_{\text {lim}}\\
                2 \bar{K}_{L 2}\left(v-V_{\text {lim }}\right)+4 \bar{K}_{L 3}\left(v-V_{\text {lim }}\right)^{3},
                & \text { else }
                \end{array}\right.,
                \\
                \frac{\partial C_D}{\partial v} &= \left\{\begin{array}{ll}
                0, 
                & \text { if } v>V_{\text {lim }} \\
               K_{D 2} \frac{C_L}{\partial v}+2 K_{D 3} C_L \frac{\partial C_L}{\partial v}, 
                & \text { else } 
                \end{array}\right.,
\end{align}
which can be used in turn to compute the partial derivatives for lift and drag with respect to velocity:
\begin{align}
                \frac{\partial L}{\partial v} &= \left\{\begin{array}{ll}
                \frac{R_\oplus S_{\text{ref}}}{m} C_L(v) \rho(r) v, 
                & \text { if } v>V_{\text {lim}}\\
                \frac{R_\oplus S_{\text{ref}}}{m} C_L(v) \rho(r) v+\frac{\partial C_L}{\partial v} \frac{R_\oplus S_{\text{ref}}}{2 m} p(r) v^{2},
                & \text { else } 
                \end{array}\right., \\
                \frac{\partial D}{\partial v} &= \left\{\begin{array}{ll}
                \frac{R_\oplus S_{\text{ref}}}{m} C_D(v) \rho(r) v, 
                & \text { if } v>V_{\text {lim }} \\
               \frac{R_\oplus S_{\text{ref}}}{m} C_D(v) \rho(r) v+\frac{\partial C_D}{\partial v} \frac{R_\oplus    S_{\text{ref}}}{2 m} p(r) v^{2},
                & \text { else } 
                \end{array}\right.,
\end{align}
which are needed to compute Jacobians of the dynamics, used for convex approximations, as described in Section \ref{sec:convex-dyn}.

\subsubsection{Bank Angle and Angle-of-attack Control} \label{subsec:bankaoa-ctrl}

For the extended model assuming both bank angle and angle-of-attack as control inputs, the control constraints are given as:
\begin{align}
    C(x,u)= \left[ 
        \begin{array}{l}
        |\sigma|-\sigma_{\max } \\
        \alpha - \alpha_{\max}(v) \\
        \alpha_{\min}(v) - \alpha\\
        \end{array}
        \right]
        \leq 0,
        \label{eq:bankaoa-cnst}
\end{align}  
for $u=[\sigma,\alpha]^\top$. Bank angle limit $\sigma_{\text{max}}$ is defined as in Section \ref{subsec:bank-ctrl}. The upper and lower bounds on $\alpha$ are given as nonconvex functions of velocity:
\begin{align}
    \alpha_{\min}(v) &= \max\left(
    0, \quad K_{\alpha1} - K_{\alpha2}(\min(v,V_{\text{lim}}) - V_{\text{lim}})^{2} - \delta \alpha_{\text{slack}}
    \right), \\
    \alpha_{\max}(v) &= \min\left(
    K_{\alpha1}, K_{\alpha1} - K_{\alpha2}(\min(v,V_{\text{lim}}) - V_{\text{lim}})^{2} + \delta \alpha_{\text{slack}}
    \right),
\end{align}
which relaxes the original design profile for angle-of-attack given in Equation \eqref{eq:wang-alpha-lut} into an inequality constraint within $\pm \delta \alpha_{\text{slack}} = 5^\circ$ that saturates above velocity limit $V_{\text{lim}}$. The control rate constraints are modeled as:
\begin{align}
\mathcal{\dot{U}} \triangleq
    \left\{ (\dot{\sigma},\dot{\alpha}) \in \real 
        \left|
        \left[ 
        \begin{array}{l}
        | \dot{ \sigma } | -\dot{ \sigma }_{\max } \\
        | \dot{ \alpha } | -\dot{ \alpha }_{\max }
        \end{array}
        \right]
        \leq 0
        \right.
    \right\}
\end{align}
where $\dot{\sigma}_{\text{max}}$ is defined as in Section \ref{subsec:bank-ctrl}, and $\dot{\alpha}=1^\circ/\text{s}$.
After this reformulation, the partial derivatives for the lift and drag coefficients can now be computed explicitly as a piece-wise function of angle-of-attack:
\begin{align}
                \frac{\partial C_L}{\partial \alpha} &= \left\{\begin{array}{ll}
                0, 
                & \text { if } v>V_{\text {lim}}\\
                \hat{K}_{L2}+2 \hat{K}_{L3} \alpha^o,
                & \text { else }
                \end{array}\right.,
                \\
                \frac{\partial C_D}{\partial \alpha} &= \left\{\begin{array}{ll}
                0, 
                & \text { if } v>V_{\text {lim }} \\
               \hat{K}_{D2}+2 \hat{K}_{D3} \alpha^o+3 \hat{K}_{D4}\left(\alpha^o\right)^2+4 \hat{K}_{D5}\left(\alpha^o\right)^3, 
                & \text { else } 
                \end{array}\right.,
\end{align}
which can be used in turn to compute the piece-wise partial derivatives for lift and drag with respect to angle-of-attack:
\begin{align}
                \frac{\partial L}{\partial \alpha} &= \left\{\begin{array}{ll}
                0, 
                & \text { if } v>V_{\text {lim}}\\
                \frac{\partial C_L}{\partial \alpha} \frac{R_\oplus}{2m} p(r) v^{2},
                & \text { else } 
                \end{array}\right., \\
                \frac{\partial D}{\partial \alpha} &= \left\{\begin{array}{ll}
                0, 
                & \text { if } v>V_{\text {lim }} \\
                \frac{\partial C_D}{\partial \alpha} \frac{R_\oplus}{2m} p(r) v^{2},
                & \text { else } 
                \end{array}\right. .
\end{align}

The partial derivatives of lift and drag with respect to velocity are now given as:
\begin{align}
\frac{\partial L}{\partial v} &= \frac{R_\oplus S_{\text{ref}}}{m} C_L(v) \rho(r) v,  \\
\frac{\partial D}{\partial v} &= \frac{R_\oplus S_{\text{ref}}}{m} C_D(v) \rho(r) v.
\end{align}

\subsection{Path Constraints}

Path constraints, which are highly nonlinear functions of both $x$ and $u$, are imposed to ensure that the the peak heat rate, dynamic pressure, and normal load experienced by the vehicle remain within allowable limits across the trajectory. In addition, no-fly zones (NFZs) may be defined over Earth's surface, which specify regions that the vehicle must avoid in the latitude-longitude-plane. Together, these path constraints are represented after nondimensionalization as:
\begin{align}
    P(x,u)  = \left[
                \begin{array}{l}
                \bar{k}_{Q} \exp (-\frac{1}{2} \beta R_\oplus [r-1]) v^{3}-1 \\
                \bar{k}_{q} \exp (-\beta R_\oplus [r-1]) v^{2}-1 \\
                \bar{k}_{n} \exp (-\beta R_\oplus [r-1]) v^{2} \sqrt{C_L^2 + C_D^2}-1 \\
                R_{NFZ,j}^2 - \left(\theta-\theta_{NFZ,j}\right)^2 -  \left(\phi-\phi_{NFZ,j}\right)^2 
                \end{array}
            \right] \leq 0,
            \label{eq:path-constr-spherical}
\end{align}
where $\bar{k}_{Q}= \left( k_Q \rho_{\oplus} \left(\sqrt{g_{\oplus} R_{\oplus}}\right)^3 \right)/ \left( \Dot{Q}_{\text{max}} \right)$ denotes the peak heat rate coefficient, $\bar{k}_{q} = \left( \rho_{\oplus} g_{\oplus} R_{\oplus} \right)/ \left( 2 \ q_{\text{dyn,max}} \right)$ denotes the peak dynamic pressure coefficient, and $\bar{k}_{n} = \left( \rho_{\oplus} R_{\oplus} S_{\text{ref}} C_L^* / ^2 \right)/ \left( 2 \ m \ n_{g,\text{max}} \right)$ denotes the peak normal load coefficient. For the RLV model, the dimensional allowable limits used to normalize these path constraints are given as $\Dot{Q}_{\text{max}}=3.\Bar{3} \cdot 10^4 \text{W}/\text{m}^2$,  $ q_{\text{dyn,max}} = 18 \cdot 10^3 \text{N} / \text{m}^2 $, and $n_{g,\text{max}} = 2.5*g_{\oplus} \text{m}/\text{s}^2$. Constant $k_Q \approx 1.2036 \cdot 10^{-5}$. If NFZs are present, they are expressed as circular exclusion zones on Earth's surface with center at longitude and latitude coordinates $(\theta_j , \phi_j)$, radius $R_{j}^{NFZ}$, and infinite altitude  $\forall j \in \{1,...,n_{\text{NFZ}}\}$ \cite{quals_paper}.

\subsection{Continuous-time Nonconvex Hypersonic Reentry Problem} \label{sec:ncvx-ocp}

A generalized nonconvex optimal control problem is expressed in Problem \ref{prob:ncvx-entry}. This template can be used to formulate any hypersonic reentry trajectory optimization problem with compatible constraints. The goal is to minimize cost $J(\cdot)$ by finding an optimal solution trajectory $x^*$ corresponding to optimal control input $u^*$ over time horizon $t_{F}^*$. In this work, we consider the cost function we seek to minimize as the terminal vehicle velocity:
\begin{align}
    J(x(t),u(t),t_F) = v(t_F),
    \label{cost:min-vf}
\end{align}
presented by \cite{wang2018autonomous}. Other literature has minimized the downrange distance traveled \cite{Rizvi2015waverider}, the integrated heat load experience by the vehicle \cite{2021scvxentry}, or even more generally a running cost and terminal cost \cite{jorris2012}. All three of the optimization variables $(x,u,t_{F})$ are free parameters, but they are constrained to obey the convex state, convex control, path and time horizon constraints, as well as the vehicle dynamics. Boundary conditions on the initial and terminal state of the trajectory are given in Equations \eqref{eq:prob-ncvx-tc}. Additionally, rate limits on the control input may be imposed, as shown in Equations \eqref{eq:bank-cnst} and \ref{eq:bankaoa-cnst}.

\boxing{th!}{problem}{prob:ncvx-entry}{16cm}{Nonconvex Hypersonic Reentry Problem}{
\begin{subequations}
\begin{align}
\text{\underline{Objective}}: & \
\underset{x(t), u(t), t_F}{\text{minimize}}\;\; 
    v(t_{F}), \label{eq:nl-cost}
\qquad \qquad \qquad 
\qquad \\
& \qquad \quad \text { s.t. } \
t \in [t_I,t_F]\\
\text{\underline{Dynamics}}: &
\qquad \qquad \quad
\dot{x}(t)=f(x(t), u(t)), \label{eq:prob-ncvx-dyn} \\
\text{\underline{State constraints}:}
& \qquad \qquad \quad x(t) \in \mathcal{X}, \\
\text{\underline{Control constraints}:}
& \qquad \qquad \quad 
C(x(t),u(t)) \leq 0, \ \label{eq:prob-ncvx-ctrl} \\
\text{\underline{Control rate constraints}:}
& \qquad \qquad \quad 
\dot{u}(t) \in \mathcal{\dot{U}}, \ \label{eq:prob-ncvx-ctrl-rt} \\
\text{\underline{Path constraints}:}
& \qquad \qquad \quad 
P(x(t),u(t)) \leq 0, \label{eq:prob-ncvx-path}  \\
\text{\underline{Time horizon constraints}:}
& \qquad \qquad \quad 
t_{F,\text{min}} < t_{F} \leq t_{F,\text{max}}, \label{eq:prob-ncvx-time} \\
\text{\underline{Boundary conditions}}:
& \qquad \qquad \quad 
x(t_I)-x_{0}=0,  \quad 
x(t_F) \in \mathcal{X}_F. \label{eq:prob-ncvx-tc}  
\end{align}
\end{subequations}
}

\section{Finite-Dimensional Problem Construction}

The infinite-dimensional, free-final-time nonconvex optimal control problem, presented in Problem \ref{prob:ncvx-entry}, is reformulated as an approximate finite-dimensional problem before a solution method is applied. First, a set of discrete temporal nodes are selected in order to sample the continuous-time signals from the original problem. The approximation is due to a parameterization on the control input, which assumes that the continuous-time control signal can be modeled in closed form as an analytical function of the control at the sample points. This control parameterization is then used in a multiple-shooting approach. In order to model final time $t_{F}$ as an optimization variable, time-interval dilation is also applied. This overall approach is similar to the work presented in \cite{elango2024successive}. The resulting finite-dimensional, fixed-final-time nonconvex optimal control problem can then be solved with a direct method, such as sequential convex programming (SCP) \cite{mceowen2023scitech3DoF}. The specific algorithm considered in this work is successive convexification, a variant of SCP \cite{szmuk2018successive}.

\subsection{Discretization}

We consider continuous-time nonlinear system dynamics, such as those presented in Equation \eqref{eq:gen-ncvx-dyn} where physical time (in seconds) is the independent variable. We select a set of discrete temporal nodes and corresponding sample points for the state and control: 
\begin{align}
    \left( t_k, ~x_k \triangleq x(t_k), ~ u_k \triangleq u(t_k) \right),
    ~\forallnodes.
\end{align}
Note that this produces time-intervals with timestep horizons:
\begin{align}
    T_{k} = (t_{k+1}-t_{k}), ~\forallintervals
\end{align}
along the trajectory. The finite-dimensional decision variables become $T_k$, $x_k$, and $u_k$. Across each time-interval $[t_{k}, t_{k+1}]$, we can determine the corresponding trajectory segment as:
\begin{align}
    x(t) &= x(t_{k})+\int_{t_{k}}^{t} f(x(s), u(s)) d s, \quad ~t \in [t_{k}, t_{k+1}], ~\forallintervals.
    \label{eq:trajseg}
\end{align}

\begin{figure}[!ht]
    \centering
    \includegraphics[width=0.75\linewidth]{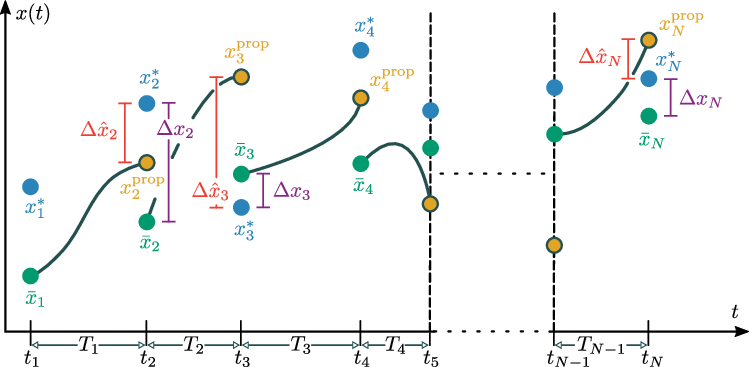}
    \caption{Variables introduced into the stitching constraints are depicted. Each continuous trajectory segment represents a piecewise integration across each timestep used in the exact discretization, where $\xbar_k$ represents the initial condition and $\xprop_{k+1}$ represents the final condition.}
    \label{fig:deviation_discretization}
\end{figure}

\subsection{Time-Interval Dilation and Control Parameterization}

We introduce normalized time $\tau$ across time-interval $k$, which can be mapped from a physical time $t$:
\begin{align}
    \tau(k, t) = \frac{t-t_{k}}{T_{k}}, 
\end{align}
which maps $t:\{1, \cdots, N-1\} \times[0,1] \rightarrow [0, t_{F}]$ across interval $k$. This mapping is bijective. Physical time within interval $k$ can be recovered from the normalized time $\tau$ as follows:
\begin{align}
    t(k, \tau) = t_{k} + T_{k} \tau, 
\end{align}
which inverts back to $\tau:\{1, \cdots, N-1\} \times [0, t_{F}] \rightarrow [0,1]$. Note that as physical time sweeps between $[t_{k},t_{k+1}]$ across the interval, normalized time sweeps $[0,1]$. At the start of each new interval, normalized time resets to $\tau=0$. This mapping can be used to convert any function of $t$ over time-interval $k$ to a new basis where $\tau$ is the independent variable (and $t$ becomes dependent):
\begin{align}
    \square(k, \tau) \triangleq \square (t_{k}+T_{k} \tau),
    \label{eq:taufunc}
\end{align}
where $\square$ is a placeholder for any time-varying variable (or signal). The derivative with respect to $\tau$ for any function in the form of Equation \eqref{eq:taufunc} is denoted:
\begin{align}
    \dotc{\square}(k, \tau) \triangleq \frac{d \square(k, \tau)}{d \tau}.
\end{align}
For example, the nonlinear system dynamics can now be expressed with normalized time as the independent variable across each interval: 
\begin{align}
\begin{split}
    \dotc{x} &\triangleq \frac{d x(k, \tau)}{d \tau}
    = \frac{d x(k, \tau)}{d t(k, \tau)} \frac{d t(k, \tau)}{d \tau}
    = f(x(k, \tau), u(k, \tau)) \cdot T_{k}
    \triangleq F(x(k, \tau), u(k, \tau), T_{k}),
\end{split}
\label{eq:nonlin-sys-tau}
\end{align}
where $F(\cdot)$ denotes the nonlinear dynamics with respect to normalized time. With a slight abuse of notation, we redefine the solution of the nonlinear system across the interval $k$:
\begin{align}
\begin{split}
    x(k, \tau) &\triangleq x\left(t_{k}+T_{k} \tau\right) 
            = x(k, 0) 
            + \int_{0}^{\tau} T_{k} f(x(k, s), u(k, s)) ds, 
            \equiv x(t),
\end{split}
\label{eq:xfin}
\end{align}
to be a function of normalized time $\forallintervals, ~t \in\left[t_{k}, t_{k+1}\right],$ and $\tau \in[0,1]$. Note that this is equivalent to Equation \eqref{eq:trajseg}, and can be seen as a change of basis. This reformulation is equivalent to adding physical time $t$ to the state, and applying control input $\hat{u}=\frac{d t}{d \tau}=T_{k}$ with a zero-order hold across each interval, similar to what is presented in \cite{elango2024successive}. Note that for the nonlinear system:
\begin{align}
    x(k, 0) &\equiv x(t_{k})=x_{k}, \\ 
        x(k, 1) &= x(k+1,0) \equiv x\left(t_{k+1}\right)=x_{k+1}
        \triangleq \xprop_{k+1},
    \label{eq:xprop}
\end{align}
 where $\xprop$ is the multiple-shooting solution across interval $k$ with parameterized control, as shown in Figure \ref{fig:deviation_discretization}. In order to propagate the nonlinear system, we assume a first-order-hold (FOH) parameterization on the control across the interval:
\begin{align}
    u(k, \tau) \triangleq u\left(t_{k}+T_{k} \tau\right)=(1-\tau) u_{k}+\tau u_{k+1}, \label{eq:foh} 
\end{align}
also a function of normalized time $\forallintervals, ~t \in\left[t_{k}, t_{k+1}\right],$ and $\tau \in[0,1]$.

\subsection{Finite-Dimensional Nonconvex Hypersonic Reentry Problem}

The finite-dimensional nonconvex optimal control problem presented in Problem \ref{prob:ncvx-entry-finite} is a discrete approximation of Problem \ref{prob:ncvx-entry}. The solution of this approximate, finite-dimensional problem can be solved with direct methods.


\boxing{th!}{problem}{prob:ncvx-entry-finite}{16cm}{Finite-Dimensional Nonconvex Hypersonic Reentry Problem}{
\begin{subequations}
\begin{align}
\text{\underline{Objective}}: & \
\underset{\substack{x_{k}, u_{k} ~\forallnodes \\ \qquad T_{k} ~\forallintervals }}{\text{minimize}}\;\; 
    v(t_{N}), \label{eq:prob-fin-nl-cost}
\qquad \qquad \qquad 
\qquad \\
& \qquad \quad \text { s.t. } \ \\
\text{\underline{Dynamics}}: &
\qquad \qquad
x_{k+1} = x_{k}+\int_{0}^{1} T_{k} f(x(k, s), u(k, s)) ds, \quad ~\forallintervals, \label{eq:prob-fin-ncvx-dyn} \\
\text{\underline{State constraints}:}
& \qquad \qquad x_{k} \in \mathcal{X}, \qquad \qquad ~\forallnodes,\\
\text{\underline{Control constraints}:}
& \qquad \qquad 
C(x_k,u_k) \leq 0, \ \quad ~\forallnodes, \ \label{eq:prob-fin-ncvx-ctrl} \\
\text{\underline{Control rate constraints}:}
& \qquad \qquad 
\frac{u_{k+1} - u_{k}}{T_{k}} \in \mathcal{\dot{U}} \subseteq \real^{\nnu}, \ \ \ ~\forallnodes, \ \label{eq:prob-fin-ncvx-ctrl-rt} \\
\text{\underline{Path constraints}:}
& \qquad \qquad 
P(x_k,u_k) \leq 0, \quad  ~\forallnodes, \label{eq:prob-fin-ncvx-path}  \\
\text{\underline{Time horizon constraints}:}
& \qquad \qquad 
t_{F,\min} \leq \sum_{k=1}^{(N-1)} T_k \leq t_{F,\max}, \ \
~T_{\text{min}} < T_{k} \leq T_{\text{max}}, ~\forallintervals,\label{eq:prob-fin-ncvx-time} \\
\text{\underline{Boundary conditions}}:
& \qquad \qquad 
x_1-x_{I}=0, \qquad 
x_N \in \mathcal{X}_F. \label{eq:prob-fin-ncvx-tc}  
\end{align}
\end{subequations}
}
\section{Auto-tuned Primal-dual Successive Convexification (\autoscvx)}

The {\autoscvx} algorithm is outlined in the sections below. First, a generic template for nonconvex optimal control problems is posed. Then, the methodology extensions for optimizing both primal and dual variables are described. Convex approximate subproblems are formed and solved iteratively until convergence to a solution. Conveniently, the dual variable solution has a closed form for each iteration \cite{bertsekas2014constrained}. For the primal variable optimization, this enhanced sequential convex programming (SCP) algorithm applies a deviation variable approach to construct a sequence of finite-dimensional, convex subproblems that form local approximations of the optimal control problem for hypersonic reentry. 
The full proposed method is laid out in Algorithm \ref{alg:autoscvx}.

\begin{figure}[ht!]
    \centering
    \includegraphics[width=\linewidth,trim={0cm 3cm 0cm 1.5cm},clip]{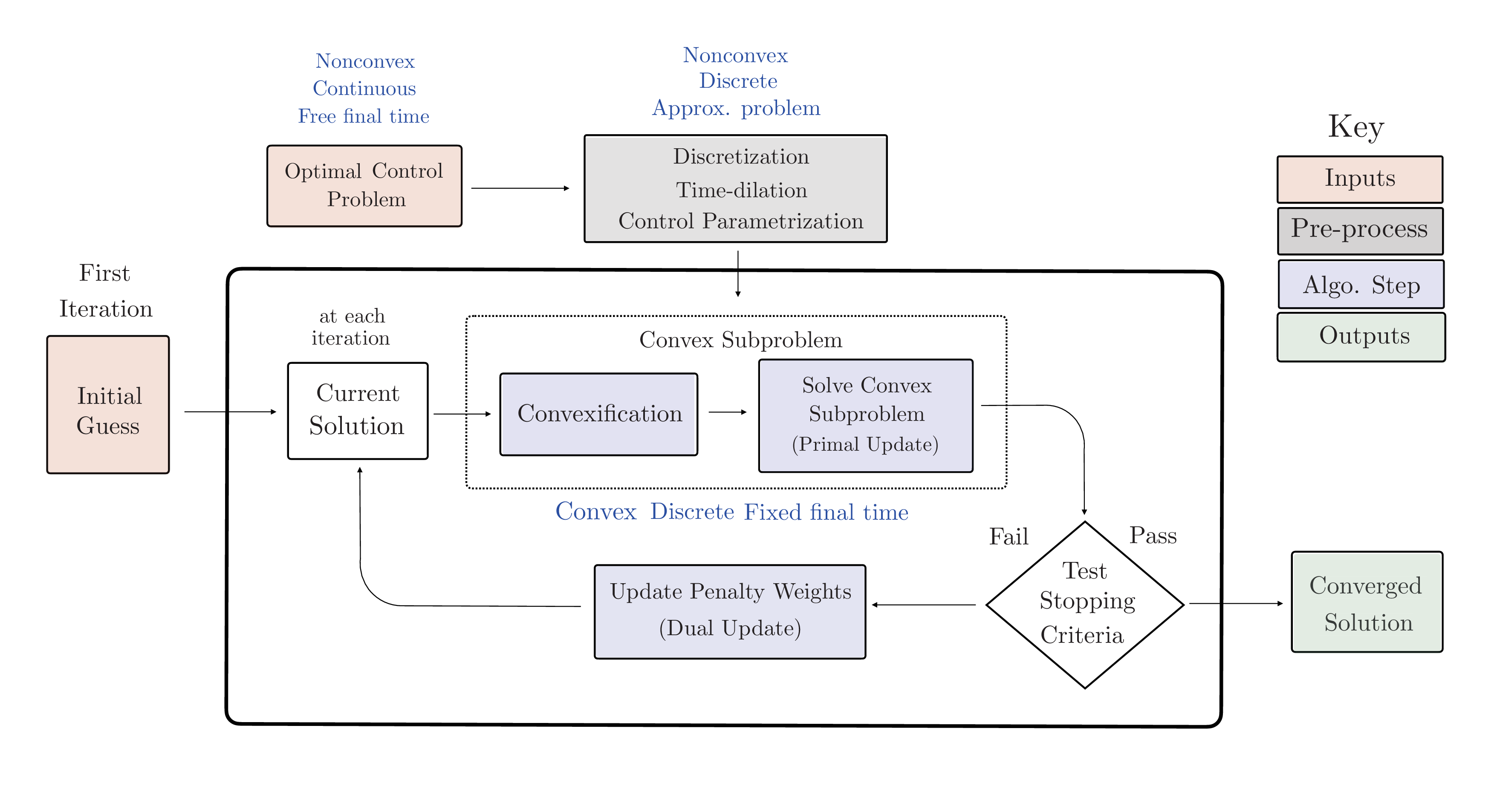}
    \caption{A block diagram for the {\autoscvx} framework for trajectory optimization. 
    }
    \label{fig:scp_framework}
\end{figure}


\subsection{Generic Nonconvex Optimal Control Problem} \label{sec:generic-nocp}

In Problem \ref{prob:ncvx-generic}, we pose a generic template for nonconvex optimal control problems that naturally subsumes the specific formulation for hypersonic reentry modeled in Problem \ref{prob:ncvx-entry}. We group all optimization variables into a single primal variable $z \triangleq [x_1^\top, \dots, x_N^\top, u_1^\top, \dots, u_N^\top, T_{1}, \dots, T_{N-1}]^\top \in \real^{\nz} $ where $\nz = (n + m + 1) N - 1 $, all convex equality and inequality constraints into $\hc(z)$ and $\gc(z)$ (respectively), and all nonconvex equality and inequality constraints into $\hn(z)$ and $\gn(z)$ (respectively). These parsed constraints are represented with vector-valued constraint functions:
 \begin{align}
     H(z) & \triangleq 
        \begin{bmatrix}
            \hc(z) \\
            \hn(z)
        \end{bmatrix} \in \real^{\nH}, \quad
    G(z) \triangleq 
        \begin{bmatrix}
            \gc(z) \\
            \gn(z)
        \end{bmatrix} \in \real^{\nG}.
        \label{eq:generic-H-G}
\end{align}   

A subset of these constraints are again grouped into:
\begin{align}
    h(z) \triangleq& \Mh H(z) \in \real^{\nequ} , \quad
    g(z) \triangleq \Mg G(z) \in \real^{\nineq} , 
    \label{eq:generic-h-g}
\end{align}
where $\Mh \in \real^{\nequ \times \nH}$ and $\Mg \in \real^{\nineq \times \nG}$ are binary (zero-one) matrices with only one nonzero element per row, chosen to select which rows of the original constraints from Equations \eqref{eq:generic-H-G} appear in the subset $h(z)$ and $g(z)$. The remaining constraints are captured by the closed set:
\begin{align}
    \mathcal{Z} \triangleq& \{ z \in \real^{n_z} \ | \  \MH H(z) = 0, \ \MG G(z) \leq 0 \},
    \label{eq:generic-Z}
\end{align}
where once again $\MH \in \real^{(\nH-\nequ) \times \nH}$ and $\MG \in  \real^{(\nG-\nineq) \times \nG}$ are binary matrices with only one nonzero element per row that select which rows of the original constraints appear in set $\mathcal{Z}$. These reformatted constraints are used to define a generic, discrete nonconvex optimal control problem, as shown in Problem \ref{prob:ncvx-generic}. We separate the constraint set $\mathcal{Z}$ because these constraints will be enforced directly, while the equality and inequality constraints in Equations \eqref{eq:prob-ncvx-eq} and \eqref{eq:prob-ncvx-ineq} (respectively) will be the focus of our analysis.

\boxing{th!}{problem}{prob:ncvx-generic}{16cm}{Generic Discrete Nonconvex Optimal Control Problem}{
\begin{subequations}
\begin{align}
\text{\underline{Objective}}: & \ 
\underset{z \in \mathcal{Z}}{\text{minimize}}\;\; 
    J(z), \label{eq:nl-cost2}
\qquad \qquad \qquad 
\qquad \\
& \qquad \ \ \ \text { s.t. } \ \\
\text{\underline{Equality constraints}}:
& \qquad \qquad 
h(z)=0, \label{eq:prob-ncvx-eq} \\
\text{\underline{Inequality constraints}:}
& \qquad \qquad 
g(z) \leq 0. \label{eq:prob-ncvx-ineq}
\end{align}
\end{subequations}
}

In the standard case, the matrices ($\Mh, \Mg, \MH, \MG$) would be chosen such that $h(z) = \hn$, $g(z) = \gn$, and $\mathcal{Z} = \{ z \in \real^{n_z} \ | \  \hc = 0, \ \gc \leq 0 \}$. In other words, the convex constraints would typically be grouped within the set $\mathcal{Z}$, while the nonconvex constraints are isolated into $h(z)$ and $g(z)$ such that they can be handled directly via successive convexification \cite{szmuk2020successive}. However, in specialized cases, the structure of functions $h(z)$ and $g(z)$ and set $\mathcal{Z}$ may take on a different form where the convex and nonconvex constraints are divided differently \cite{mceowen2023scitech3DoF}. More details of the specialized case implemented to solve the hypersonic reentry problem is discussed in Section \ref{sec:cvx-buff}.

\subsection{\autoscvx: Methodology Extensions}

Problem \ref{prob:ncvx-generic} is equivalent to:
\begin{subequations}
\begin{align}
    \underset{{\substack{z \in \mathcal{Z} \\
    p,q \geq 0}}}{\mathrm{min.}} & \quad J(z) 
    + \frac{1}{2}  p^{\top} \Wh p + \frac{1}{2}  q^{\top} \Wg  q 
\\
\text { s.t. } & 
 \quad h(z) = p, \quad g(z) \leq q,
\quad p=0 , \quad q=0,
\end{align}
\label{eq:m3-ocp}
\end{subequations}
which can be re-expressed in the equivalent form:
\begin{subequations}
\begin{align}
    \underset{\substack{z \in \mathcal{Z} \\
    p,q \geq 0}}{\mathrm{min.}}  
    \left[ \underset{\lambda, \mu \geq 0}{\mathrm{max.}} \right.
    & \left. \quad  J(z)+ \frac{1}{2}  p^{\top} \Wh p + \frac{1}{2}  q^{\top} \Wg  q + \lambda^{\top} p+\mu^{\top} q \right] \label{eq:m3-ncvx-cost}
    \\
    \text { s.t. }  &
    \quad h(z)=p, 
    \quad g(z) \leq q , 
\end{align}
\label{eq:m3-ocp2}
\end{subequations}
where constraint violation is quantified by virtual buffer variables $p$ and $q$. This problem remains equivalent to Problem \ref{prob:ncvx-generic}. Note that for Equations \eqref{eq:m3-ocp2}, the problem will always remain feasible with infinite cost, unlike Equations \eqref{eq:m3-ocp} which may become infeasible. Note that the objective being minimized with respect to primal variables $(z,p,q)$ contains a maximization with respect to dual variables $(\lambda, \mu)$. We introduce $\Lncvx(p,q,\lambda,\mu) \triangleq \lambda^{\top} p+\mu^{\top} q$ to isolate the bilinear terms involving the dual variables, which comprise the nonconvex portion of the cost in Equation \eqref{eq:m3-ncvx-cost}. The full cost function:
\begin{subequations}
\begin{align}
\Lagr_1(z, \lambda, \mu, p, q) &\triangleq 
J(z) 
+ \frac{1}{2}  p^{\top} \Wh p + \frac{1}{2}  q^{\top} \Wg  q
+ \lambda^{\top} p+\mu^{\top} q, \\
&= J(z) 
+ \frac{1}{2}  p^{\top} \Wh p + \frac{1}{2}  q^{\top} \Wg  q
+ \Lncvx(p,q,\lambda,\mu),
\end{align}
\label{eq:m3-lagr}
\end{subequations}
is the partial Lagrangian of Problem \ref{eq:m3-ocp}, with dual variables $\lambda$ and $\mu$ introduced for the equality constraints on primal variables $p$ and $q$, respectively. For convenience, we define $\nu^\top \triangleq [p^\top, q^\top, \lambda^\top, \mu^\top]$. Given a reference solution $(\zbar, \nubar)$ we next define the deviation quantities:
\begin{subequations}
\begin{align}
    \dz &\triangleq z - \zbar, \\
    \dnu &\triangleq
        [ 
          \dpp^{\top} , 
          \dq^{\top} , 
          \dlam^{\top} , 
          \dmu^{\top}
        ]^{\top} 
    \triangleq
        [
          (p - \pbar)^{\top} , 
          (q - \qbar)^{\top} , 
          (\lambda - \lambar)^{\top} , 
          (\mu - \mubar)^{\top}
        ]^{\top} , 
\end{align}
\end{subequations}
as well as the first-order Taylor series expansion operator:
\begin{align}
    \lin{c}{y} \triangleq c(\bar{y}) +  \left. \frac{\partial c}{\partial y} \right|_{\bar{y}} \Delta y,
    \label{eq:taylor-series-approx}
\end{align}
that produces a linear approximation in proximity to any point $\bar{y}$ for any nonconvex function $c(y)$ and small perturbation $\Delta y$. In this way, we produce first-order approximations of $\Lncvx(p,q,\lambda,\mu)$, $h(z)$ and $g(z)$ and notice:
\begin{subequations}
\begin{align}
    \begin{split}
    \arg \underset{\substack{z \in \mathcal{Z^{'}} \\
    p,q \geq 0}}{\mathrm{min.}}  
    \left[ \underset{\substack{\lambda, \mu \geq 0}}{\mathrm{max.}} \right.
    & \quad J(z)+ \frac{1}{2}  p^{\top} \Wh p + \frac{1}{2}  q^{\top} \Wg  q + \linfull{\Lncvx}{\nu} 
    \\ & \left. \qquad \quad + \prox{s_{z}}{z} + \prox{s_{\lambda}}{\lambda} + \prox{s_{\mu}}{\mu} \right]
    \end{split}
    \\
    \text{ s.t. } & \ \  \linfull{h}{z}=p, \ \linfull{g}{z} \leq q,
    \\
    & \qquad \qquad \equiv \notag \\
    \begin{split}
    \arg \underset{\substack{\dz + \zbar  \in \mathcal{Z^{'}} \\
    p,q \geq 0}}{\mathrm{min.}}  
    \left[ \underset{\substack{\dlam, \dmu + \mubar \geq 0}}{\mathrm{max.}} 
    \right.
    & \quad  J(\dz + \zbar) + \frac{1}{2}  p^{\top} \Wh p + \frac{1}{2}  q^{\top} \Wg  q + \lambar^{\top} p + \mubar^{\top} q +  \pbar^{\top} \dlam + \bar{q}^{\top} \dmu \\
    & \left. \qquad \quad + \proxdev{s_{z}}{z} + \proxdev{s_{\lambda}}{\lambda} + \proxdev{s_{\mu}}{\mu} \right]
    \end{split}
    \\
    \text{ s.t. } & \ \  \lin{h}{z}=p, \ \lin{g}{z} \leq q,
    \\
    & \qquad \qquad \equiv \notag \\
\begin{split}
    \arg \underset{\substack{\dz + \zbar \in \mathcal{Z^{'}} \\
    p,q \geq 0 \\ \dlam, \dmu + \mubar \geq 0}}{\mathrm{min.}}
    & \quad  J(\dz) + \frac{1}{2}  p^{\top} \Wh p + \frac{1}{2}  q^{\top} \Wg  q + \lambar^{\top} p + \mubar^{\top} q -  \pbar^{\top} \dlam - \bar{q}^{\top} \dmu \\
    & \qquad \quad + \proxdev{s_{z}}{z} + \proxdev{s_{\lambda}}{\lambda} + \proxdev{s_{\mu}}{\mu} 
    \end{split}
    \\
    \text{ s.t. } & \ \  \lin{h}{z}=p, \ \lin{g}{z} \leq q,
\end{align}
 \label{eq:m3-argmin}
\end{subequations}
where proximal costs have been added for $(\dz,\dlam,\dmu)$ to keep the solutions within a trust region of the linearization point and avoid artificial unboundedness from below. In addition, the set $\mathcal{Z}^{'}$ contains linear approximations of any nonconvex approximations included in set $\mathcal{Z}$. In this process, we discover an equivalent, local convex approximation of the problem in Equations \eqref{eq:m3-ocp} converted into a single minimization. We arrive at the approximate convex subproblem:
\begin{subequations}
    \begin{align}
    \begin{split}
    \underset{\substack{\dz + \zbar \in \mathcal{Z^{'}}, p, q \geq 0, \\
    \dlam,\dmu + \mubar \geq0}}{\mathrm{min.}}   
    & \ \ J(\dz + \zbar) + \frac{1}{2}  p^{\top} \Wh p + \frac{1}{2}  q^{\top} \Wg  q 
    + \lambar^{\top} p + \mubar^{\top} q -  \pbar^{\top} \dlam - \bar{q}^{\top} \dmu \\
    & \quad \quad 
    + \proxdev{s_{z}}{z}  + \proxdev{s_{\lambda}}{\lambda} + \proxdev{s_{\mu}}{\mu}  
    \end{split}
    \\
    \text { s.t. } & \ \  \lin{h}{z}=p, \ \  
    \lin{g}{z} \leq q.
    \end{align}
\label{eq:m3-cocp1}
\end{subequations}

The optimal primal variable solution in Equation \eqref{eq:m3-cocp1} can be solved independently:
\begin{subequations}
\begin{align}
\begin{split}
(\dzstar, \pstar, \qstar) = \arg \underset{\substack{\dz + \zbar \in \mathcal{Z^{'}}, \\ p, q \geq 0}}{\mathrm{min.}}
& \quad  J(\dz + \zbar) + \frac{1}{2}  p^{\top} \Wh p + \frac{1}{2}  q^{\top} \Wg  q + \lambar^{\top} p + \mubar^{\top} q 
+ \proxdev{s_{z}}{z}, 
\end{split}
\\
\text { s.t. } & \quad \lin{h}{z}=p, \quad 
\lin{g}{z} \leq q,
\end{align}
\label{eq:m3-primal}
\end{subequations}
from the dual variable solution. Once separated, the optimal dual variables can be determined analytically:
\begin{subequations}
\begin{align}
\dlamstar & =\argmin{\dlam}         
\proxdev{s_{\lambda}}{\lambda} - \dlam^{\top} \pstar 
= s_{\lambda} \pstar, 
\\
\dmustar  &=\argmin{\dmu + \mubar \geq 0} 
\proxdev{s_{\mu}}{\mu} - \dmu^{\top} \qstar 
=\max \left( -\mubar, s_{\mu} \qstar \right) ,
\end{align}
\label{eq:m3-dual}
\end{subequations}
which gives an update rule for the linear coefficients $(\lambar,\mubar)$ in Equation \eqref{eq:m3-primal} (similar to augmented Lagrangian \cite{bertsekas2014constrained,wan2022entry}). However, it would still be desirable to determine the quadratic penalty weights $(\Wh, \Wg)$ for the virtual buffer terms. 

With this motivation in mind, we take the full Lagrangian of Equations \eqref{eq:m3-primal}:
\begin{align}
\begin{split}
    \Lagr_2(z, p, q)  \triangleq J(\dz + \zbar) & 
    +\frac{1}{2} p^{\top} \Wh p + \frac{1}{2} q^{\top} \Wg q + \lambar^{\top} p+\mubar^{\top} q + \proxdev{s_{z}}{z} 
    \\
    & +\lamhat^{\top}(\lin{h}{z}-p)+\muhat^{\top}(\lin{g}{z}-q),
\end{split}
\label{eq:m3-lagr2}
\end{align}
where we see that dual variables $(\lamhat, \muhat)$ drive the value of the constraints towards the value of the virtual buffer variables. Due to this, in order to achieve feasibility of Problem \eqref{eq:m3-primal}, we seek positive-definite weights $(\Wh,\Wg)$ that drive virtual buffers $(p,q)$ to zero, respectively. Taking the stationarity conditions of \eqref{eq:m3-lagr2} yields \cite{nocedal1999}:
\begin{subequations}
    \begin{align}
        \frac{\partial \Lagr_{2}}{\partial p} &= 
        {\pstar}^{\top} \Wh + (\lambar-\lamhat)^{\top} = 0, \\
        \frac{\partial \Lagr_{2}}{\partial q} &= 
        {\qstar}^{\top} \Wg + (\mubar-\muhat)^{\top} = 0, 
    \end{align}
\end{subequations}
which reveals a relationship between the buffer variables, weights and dual variables:
\begin{align}
    \Wh \pstar &= \lamhat - \lambar, \\
    \Wg \qstar &= \muhat - \mubar,
\end{align}
where each positive-definite matrix is constructed with diagonal elements 
\begin{align}
        \Wh = 
        \begin{bmatrix}
            \omega_{h_{1}} & & \\
            & \ddots & \\
            & & \omega_{h_{\nequ}} 
        \end{bmatrix}, \quad
        \Wg = 
        \begin{bmatrix}
            \omega_{g_{1}} & & \\
            & \ddots & \\
            & & \omega_{g_{\nineq}} 
        \end{bmatrix}.
\end{align}
If we aim to drive our virtual buffers towards a prescribed feasibility tolerance:
\begin{align}
\pstar_{i} &= \frac{\lamhat_{i} - \lambar_{i}}{\wh{j}} = \epsh{i}, \\
\qstar_{i} &= \frac{\muhat_{i} - \mubar_{i}}{\wg{j}} = \epsg{j}, 
\end{align}
%
$\forall i \in [1,\dots,\nequ]$ and $\forall j \in [1,\dots,\nineq]$, then we may exploit this relationship to choose a corresponding weight update for each diagonal element of $(\Wh,\Wg)$: 
\begin{align}
    \wh{i} &\gets \frac{\lamhat_{i}-\lambar_{i}}{\epsh{i}} 
    = \frac{\wh{i} \pstar_{i}}{\epsh{i}} , \\
    \wg{j} &\gets \frac{\muhat_{j}-\mubar_{j}}{\epsg{j}} 
    = \frac{\wg{j} \qstar_{j}}{\epsg{j}}.
\end{align}
Note that the update scheme can be represented in closed form as a function of the primal solution variables and hyperparameters of the convex subproblem.


\subsection{Summary: {\autoscvx} Algorithm}

The full {\autoscvx} algorithm is displayed in Algorithm \ref{alg:autoscvx}. After initialization, the primal variables, dual variables and quadratic penalty weights are solved iteratively until convergence. The ``Solve Convex Subproblem (Primal Update)'' block in Figure \ref{fig:scp_framework} is represented by Step \ref{alg-primal-update}. The ``Update Penalty Weights (Dual Update)'' block in the same diagram subsumes Steps \ref{alg-quadratic-update} and \ref{alg-dual-update}. Note that both the dual variables and the quadratic penalty weights are updated in closed-form. 

\begin{figure}[th!]
  \makebox[\linewidth]{%
  \begin{minipage}{\dimexpr\linewidth-5em}
    \begin{method}[H]
    \caption{Auto-tuned Primal-dual Successive Convexification (\autoscvx)} 
    \label{alg:autoscvx}
    \begin{algorithmic}[1]
        \Require Initial guess $\zbar_0$, step sizes ($s_{z}, s_{\lambda}, s_{\mu}$), convergence tolerances $(\epsoptz, \epsoptJ, \epsfeas)$, desired constraint residual $(\epsh{},\epsg{})$, minimum quadratic weight $\varepsilon_{\text{w,}\min}$.  
        \State Initialize $\Wh = \diag{\one{\nequ}}$, $\Wg = \diag{\one{\nineq}}$, $\lambar = \zero{\nequ}$ and $\mubar = \zero{\nineq}$.
        \While{
        $\left( | \dz | \geq \epsoptz \right.$ 
        \Or \ $\left. \infnorm{[\lin{h}{z},\lin{g}{z}]} \geq \epsfeas \right)$ \label{alg-conv1} \\
        \qquad \qquad \And \quad
        $\left( | J(\dz) | \geq \epsoptJ \right.$ 
        \Or \ $\left. \infnorm{[h(\zstar),g(\zstar)]} \geq \epsfeas \right)$ \label{alg-conv2} 
        }
            \State Update primal variables: \label{alg-primal-update}
            \begin{align*}
                \begin{split}
                \begin{pmatrix}
                \dzstar, \pstar, \qstar
                \end{pmatrix}
                = \argmin{\substack{\dz + \zbar \in \mathcal{Z^{'}}, \\ p, q \geq 0}} & \ J(\dz + \zbar) 
                + \frac{1}{2} p^{\top} \Wh p
                +\frac{1}{2} q^{\top} \Wg q + \lambar^{\top} p + \mubar^{\top} q 
                + \proxdev{s_{z}}{z} 
                \end{split}
                \\
                \text{s.t. } & \ \lin{h}{z}=p , \quad 
                \lin{g}{z} \leq q.
            \end{align*}
            \State Update penalty weights (enforce $(\Wh, \Wg) \geq \varepsilon_{\text{w,}\min}$):
            \begin{align*}
                {\Wh} \gets & \diag{ \frac{\pstar_1 \wh{1} }{\epsh{1}} ,\dots, \frac{\pstar_{\nequ} \wh{\nequ}}{\epsh{\nequ}}} , \quad 
                {\Wg} \gets \diag{\frac{\qstar_{1} \wg{1}}{\epsg{1}}, \dots,  \frac{\qstar_{\nineq} \wg{\nineq}}{\epsg{\nineq}} }.        
            \end{align*}
            \label{alg-quadratic-update}
            \State Update dual variables:
            \begin{align*}
                \quad \dlamstar \gets s_{\lambda} \pstar, \qquad 
                \dmustar \gets  \max(-\mubar,  s_{\mu} \qstar).
            \end{align*}
            \label{alg-dual-update}
            \State Update reference values for next iteration:
            \begin{align*}
                \quad \zbar \gets \dzstar + \zbar, \qquad
                \lambar \gets \dlamstar + \lambar, \qquad
                \mubar \gets \dmustar + \mubar.
            \end{align*}
        \EndWhile
        \Ensure Converged solution $\zstar = \dzstar + \zbar$ (final value at convergence).
    \end{algorithmic}
\end{method}
\end{minipage}
}%
\end{figure}


\subsection{Primal Variable Update: Implementation Details}

In this section, the formulation of the discrete, convex subproblem in Equations \eqref{eq:m3-primal} is discussed for the hypersonic reentry optimal control problem. The resulting Problem \ref{prob:cvx-entry} is solved iteratively to update the primal variables until convergence. Here we will revert back to the variables $(x_{k},u_{k}) ~\forallnodes$ and $T_{k} ~\forallintervals$ from Problem \ref{prob:ncvx-entry-finite} to clearly present the convexification of the dynamics constraint in Equation \eqref{eq:prob-fin-ncvx-dyn} with respect to deviation variables. 

\subsubsection{Convex Approximation of the Dynamics with Deviation Variables} \label{sec:convex-dyn}

The discretized nonlinear system dynamics are approximated as a linear time varying (LTV) system by solving an initial value problem. This variational approach is inverse-free and exact (due to multiple shooting) \cite{bock1984multiple,lin2014surveyopt}, meaning the nonlinear system dynamics constraint is satisfied to arbitrary precision at convergence (when the first approximation terms drop to zero). The details of the procedure for computing the linearized matrices $(\Ak,\Bmk,\Bpk,\Sk)$ are given in the Appendix. These are used in the discretized, convex dynamics constraint of the primal subproblem:
\begin{align}
    \dxhat_{k+1} = \dx_{k+1} + \xbar_{k+1} - \xprop_{k+1} &= \Ak \dx_{k} + \Bmk \du_{k} + \Bpk \du_{k+1} + \Sk \dT_{k},
    \label{eq:ltv-dynamics}
\end{align}
 where $(\xbar_{k}, \ubar_{k}, \Tbar_{k}) ~\forallintervals$ are a nominal reference trajectory. The value of state $\xprop_{k+1}$, defined in Equation \eqref{eq:xprop}, is achieved from multiple-shooting of the nonlinear system with a FOH on $\ubar_{k}$ and $\ubar_{k+1}$ across each interval $k$. The two distinct deviation quantities:
\begin{subequations}
\begin{align}
    \dx_k &\triangleq x_k - \xbar_k, \\
    \dxhat_k &\triangleq x_k - \xprop_k,
\end{align}
\end{subequations}
capture the difference between $x$ and the reference $\xbar$, and $x$ and propagated state $\xprop$, respectively. These deviation quantities are depicted in Figure \ref{fig:deviation_discretization}. This multiple-shooting approach \cite{bock1984multiple} produces a set of piecewise-continuous trajectory segments, where $\xbar_k$ represents the initial condition and $\xprop_{k+1}$ represents the terminal condition, and is an example of a variational method \cite{lin2014surveyopt}. Due to this relationship, the optimal solution at each iteration is equivalent to:
\begin{align}
    x^* &= (\dx^*_{k} + \xbar_{k}) = (\dxhat^*_{k} + x_{k}^{\text{prop}}).
    \label{eq:stitch1}
\end{align}

By embedding the state obtained from the nonlinear propagation from the discretization step into the linear dynamics constraint of the subproblem, the piecewise-smooth trajectory segments are driven to become continuous (to numerical precision of the integrator) everywhere at convergence:
\begin{subequations}
\begin{align}
    \dx_{k} \rightarrow 0 &\iff  x_{k}^{*} \rightarrow \xbar_{k}, \\
    \dxhat_{k} \rightarrow 0 &\iff  x_{k}^{*} \rightarrow x_{k}^{\text{prop}}, \quad \forallnodes .
\end{align}
\end{subequations}
The optimal solution becomes dynamically-consistent with a one-shot nonlinear integration of the vehicle dynamics, i.e. continuous-time dynamic feasibility of the nonlinear system is achieved, even over a sparse time grid. 

\boxing{ht!}{problem}{prob:cvx-entry}{16cm}{Discrete, Convex Primal Update Subproblem for Hypersonic Reentry}{
\begin{subequations}
\begin{align}
\begin{split}
\text{\underline{Objective}}: \quad & \underset{\dx, \du, \dT, p, q}{\text{minimize}}\;\; 
    v_N
    +  \underset{J_{\text{buff}}}
    {\underbrace{ \begin{bmatrix} p \\ q \end{bmatrix}^\top
    \begin{bmatrix} W_h & 0 \\ 0 & W_g \end{bmatrix}
    \begin{bmatrix} p \\ q \end{bmatrix}
    + \begin{bmatrix} \lambar \\ \mubar \end{bmatrix}^\top
     \begin{bmatrix} p \\ q \end{bmatrix} }}
    + \underset{J_{\text{tr}}}
    {\underbrace{ \proxdev{s_{x}}{x} + \proxdev{s_{u}}{u}}} 
\end{split}
\\ 
& \text { s.t. }
\\
\text{\underline{Dynamics}:}
& \qquad \dx_{k+1} + \xbar_{k+1} - \xprop_{k+1} = \Ak \dx_{k}+\Bmk \du_{k}+\Bpk \du_{k+1}+\Sk \dT_k, 
\\ & \qquad 
\ ~\forallintervals, \notag \\
\text{\underline{State constraints}:}
& \qquad \xbar_{k} \ + \dx_{k} + {q_{\mathcal{X},k}} \in \mathcal{X}, ~\forallnodes,
\\ \text{\underline{Control constraints}:} 
& \qquad \lintwo{C}{x_k}{u_k} + {q_{\text{ctrl},k}} \leq 0 , ~\forallnodes, \\
\\ \text{\underline{Control rate constraints}:}  
& \qquad \frac{((\ubar_{k+1}+\du_{k+1})-(\ubar_k+\du_k))}{\Tbar_k} \in \mathcal{\dot{U}}, ~\forallintervals, \label{eq:prob-cvx-ctrl-rt} \\
\text{\underline{Path constraints}:} \label{eq:prob-cvx-path}
& \qquad \lintwo{P}{x_k}{u_k} + {q_{\text{ncvx}}} \leq 0, ~\forallnodes,
\\
\text{\underline{Time horizon constraints}:}
& \qquad t_{F,\min} \leq \sum_{k=1}^{(N-1)} (\Tbar_k + \dT_k) \leq t_{F,\max}, \\
& \qquad ~\dT_{\text{min}} < \dT_{k} \leq \dT_{\text{max}}, ~\forallintervals, \notag \\
\text{\underline{Boundary conditions}:} 
& \qquad  \xbar_1 + \dx_{1} =  x_0, \qquad
\begin{bmatrix}
            \xbar_{N} + \dx_{N} + p_{\mathcal{X_F}} \\
            \xbar_{N} + \dx_{N} + q_{\mathcal{X_F}}
        \end{bmatrix} \in \mathcal{X_F}\\ 
\notag \\ \notag 
& \qquad \qquad \text{where:} \ 
    p \triangleq \left[ p_{\mathcal{X_F}} \right], \
    q \triangleq \left[ q_{\mathcal{X_F}}^\top, q_{\mathcal{X}}^\top , q_{\text{ctrl}}^\top, q_{\text{ncvx}}^\top \right]^\top. \notag
\end{align}
\end{subequations}
}

\subsubsection{Specialized Constraint Buffering and Initial Guess} \label{sec:cvx-buff}
As noted in Section \ref{sec:generic-nocp}, selecting which constraints appear in functions $h(z)$ and $g(z)$ (Equations \eqref{eq:generic-h-g}) and set $\mathcal{Z}$ (Equation \eqref{eq:generic-Z}) is an implementation detail left to the user. 
The convex subproblem in Equation \eqref{eq:m3-primal} enforces linear approximations of both $\lin{h}{z}=p$ and $\lin{g}{z}\leq q$, where virtual buffers $p$ and $q$ are penalized in the cost. Violation of these constraints is therefore is permitted in early iterations, but driven beneath a specified feasibility tolerance as algorithm iterations progress. In contrast, the constraints within the set $\mathcal{Z}^{'}$ (the linear approximation of set $\mathcal{Z}$) are enforced directly. In practice, it is standard to select $h(z)$ and $g(z)$ such that they contain all nonconvex constraints and $\mathcal{Z}$ such that it contains all convex constraints \cite{szmuk2020successive}. When all linearized nonconvex constraints are assigned virtual buffers, then the convex constraints do not require buffering in order to achieve a feasible solution to the subproblem.

However, in specific cases, we may wish to forego buffering a subset of the nonconvex constraints; in these instances, we buffer a subset of the convex constraints instead. In addition, we supply an initial guess to the {\autoscvx} algorithm that is feasible with respect to the unbuffered nonconvex constraints. In this work, we avoid buffering the nonconvex dynamics constraints by including the discrete dynamics in set $\mathcal{Z}$, such that the linearized dynamics are satisfied at every {\autoscvx} iteration. Instead, we add the terminal boundary conditions to $h(z)$ and $g(z)$. We supply the algorithm with a dynamically feasible initial guess. This is found by integrating the nonlinear system with a naive control input $\sigma=0^\circ$ from the initial boundary condition.

In this way, we preserve the fidelity of the LTV dynamics model and achieve high-accuracy solutions at convergence that match closely when compared against single-shooting trajectories propagated with the optimal control. This is demonstrated in numerical results in Section \ref{sec:results}. The hypersonic reentry dynamics are particularly sensitive. With small changes in the state, the state derivative may change drastically (and even hit a singularity, or become unstable). 
In addition, quadratic-linear penalty for constraint violation does not encourage sparsity of the virtual buffer variables. 
The recursive discrete dynamics constraints may accumulate constraint violation as virtual buffers may compound. These considerations motivate the design choice of leaving the nonconvex dynamics unbuffered.

\subsubsection{Convex Subproblem (Primal Update)}

The primal variables can be updated iteratively for the hypersonic reentry problem by solving the discrete, convex subproblem as shown in Problem \ref{prob:cvx-entry} via \autoscvx. This is completed in Step \ref{alg-primal-update} of Algorithm \ref{alg:autoscvx}. Both linear and quadratic penalties are present in the augmented cost $J_{\text{buff}}$ to drive the buffered constraints towards feasibility as the SCP algorithm iterates. A quadratic proximal trust region is enforced on the state and control via $J_{\text{tr}}$, while a hard trust region on the change in time of flight is introduced due to additional linearization sensitivity to this term. 
This subproblem is a quadratic program, amenable to real-time first-order solvers such as OSQP \cite{stellato2020osqp}. Due to the FOH parameterization on the control, the control rate constraint may be enforced directly as in Equation \eqref{eq:prob-cvx-ctrl-rt}, and will be exactly satisfied at convergence. 

\subsubsection{Convergence Criterion}
The convergence criteria evaluates satisfaction of two distinct quantities: an optimality tolerance, and a feasibility tolerance. The algorithm is considered to converge in the event that one of two separate criteria are satisfied (Steps \ref{alg-conv1} and \ref{alg-conv2} in Algorithm \ref{alg:autoscvx}). The first criteria is that the deviation in state variables falls within a given optimality tolerance and buffered (convex approximate) constraints are satisfied within a feasibility tolerance:
\begin{align}
    | \Delta z_{k} | \leq& \ \epsilon_{\Delta z} \ \triangleq \epsoptz, \\
    \begin{bmatrix}
        \infnorm{\lin{h}{z_k}} \\ 
        \infnorm{\lin{g}{z_k}}
    \end{bmatrix}
    \leq& \begin{bmatrix} \epsilon_{h} \\ \epsilon_{g} \end{bmatrix} \triangleq \epsilon_{\text{feas}}, 
    \quad \forallnodes.
    \label{eq:planar-conv}
\end{align}

The second alternative criteria is that the deviation in the cost falls within a given optimality tolerance, and buffered (nonconvex) constraints are satisfied within the same feasibility tolerance:
\begin{align}
    | J(\dz_k) | \leq& \ \epsilon_{\text{cost}} \ \triangleq \epsoptJ, \\
    \begin{bmatrix}
        \infnorm{h(\dz_k + \zbar_k)} \\ 
        \infnorm{g(\dz_k + \zbar_k)}
    \end{bmatrix}
    \leq& \begin{bmatrix} \epsilon_{h} \\ \epsilon_{g} \end{bmatrix} \triangleq \epsilon_{\text{feas}}, 
    \quad \forallnodes.
    \label{eq:conv}
\end{align}

\section{Numerical Results} \label{sec:results}
The proposed {\autoscvx} algorithm is demonstrated alongside existing methods in various numerical studies. These studies are performed on two hypersonic reentry optimal control problem examples. This first example models bank angle as the only control input, as described in Section \ref{subsec:bank-ctrl}. The second example models both bank angle and angle-of-attack as control variables, as described in Section \ref{subsec:bankaoa-ctrl}. The mission parameters for the constraints are described in Table \ref{tab:ex1-mission-params}. The convergence criteria is described in Table \ref{tab:ex1-mission-conv-hyperparams}. The trust region step-size hyperparameters (which determine how quickly the primal dual variables are updated in each {\autoscvx} iteration) are described in Table \ref{tab:ex1-mission-ptr-weights}. For an initial guess, $\sigma=0^\circ$ was propagated from the initial boundary condition over an initial time horizon of $t_{F}=1700$s. A grid size of $N=40$ temporal nodes were used for discretization.

\vspace{0.25in}

\begin{minipage}[b]{.49\textwidth}
\centering
\begin{tabular}{lrl}
\hline 
    Parameter       & Value                 & Description      \\ 
    \hline                                          
    $h_I$  & $100 \mathrm{~km} $            & Init. alt. \\
    $\theta_I$      & $0^\circ$             & Init. long. \\
    $\phi_I$        & $0^\circ$             & Init. lat.  \\
    $\mathrm{v}_I$  & $7450 \mathrm{~m/s}$  & Init. speed     \\
    $\gamma_I$      & $-0.5^\circ$          & Init. f.p.a. \\
    $\psi_I$        & $0^\circ$             & Init. hdg.       \\
    \hline
    $h_{F}$  &    \begin{tabular}{@{}c@{}} $15 \mathrm{~km}$ \\ $ [15,35] \mathrm{~km}$  \end{tabular}   &  \begin{tabular}{@{}l@{}} Term. alt., bank \\ `` '' bank/a.o.a.  \end{tabular}    \\
   
    $\theta_F$      & $12^\circ $           & Term. long.    \\
    $ \phi_F $      & $70^\circ$            & Term. lat.     \\
    $\gamma_F$      & $-10^\circ$           & Term. f.p.a. \\
    $\psi_F$        & $90^\circ$            & Term. hdg.   \\
    \hline
    $\sigma_I$      & $0^\circ$             & Init. bank angle \\ 
    $\sigma_{\max}$ & $80^\circ$            & Bank angle limit \\
    $\dot{\sigma}_{\max }$ & $10^\circ / \mathrm{s}$  & Bank rate limit \\ 
     $\alpha_{\min}$ & $0^\circ$        & Min. a.o.a. limit \\
     $\alpha_{\max}$ & $40^\circ$       & Max. a.o.a. limit \\
    $\dot{\alpha}_{\max }$      &  $5^\circ / \mathrm{s}$                   & a.o.a. rate limit   \\ \hline
    $\Dot{Q}_{\text{max}}$ & $33.\bar{3} \cdot 10^{3} \text{W}/\text{m}^2$ & Heat rate \\
    $q_{\text{dyn,max}}$ & $18 \cdot 10^{3} \text{N} / \text{m}^2 $ & Dyn. press. \\
    $n_{g,\text{max}}$ & $2.5*g_{\oplus}$ & Normal load \\
    $(\theta , \phi, R)_{\text{NFZ},1}$ & $(5,30,5)^\circ$ & NFZ \#1 (pos., rad.) \\
    $(\theta , \phi, R)_{\text{NFZ},2}$ & $(-6.5,50,5)^\circ$ & NFZ \#2 (pos., rad.) \\
    \hline 
\end{tabular}
\captionof{table}{Mission parameters.}
\label{tab:ex1-mission-params}
\vspace{0.18in}
\end{minipage}
\begin{minipage}[b]{.49\textwidth}
\begin{tabular}{lrl}
\hline 
    Convergence. tol.       & Value                 & Description      \\ 
    \hline   
    \multicolumn{3}{c}{\textbf{\underline{Opt. tolerances:}}} \\ 
    $\epsoptin{J,v_{f}}$ & 5 m/s & Cost (term. vel.) \\ \hline
   $\epsoptin{h}$  &    \begin{tabular}{@{}c@{}} $5 \mathrm{~km}$ \\ $10 \mathrm{~km}$  \end{tabular}   &  \begin{tabular}{@{}l@{}} Alt., bank \\ `` '' bank/a.o.a.  \end{tabular}    \\
    $\epsoptin{\theta}, \epsoptin{\phi}$      & $1^\circ$            & Lon., lat. \\
    $\epsoptin{v}$  &    \begin{tabular}{@{}c@{}} $30 \mathrm{~m/s}$  \\ $50 \mathrm{~m/s}$   \end{tabular}   &  \begin{tabular}{@{}l@{}} Vel., bank \\ `` '' bank/a.o.a.  \end{tabular}    \\
    $\epsoptin{\gamma}, \epsoptin{\psi}$      & $5^\circ$            & F.p.a., hdg. \\
    \hline
    \multicolumn{3}{c}{\textbf{\underline{Feas. tolerances:}}} \\
    $\epsfeasin{h_F}$  & $2 \mathrm{~km}$     & Term. altitude     \\
    $\epsfeasin{\theta}, \epsfeasin{\phi}$      & $2^\circ$            & Lon., lat. \\
     $\epsfeasin{\gamma}, \epsfeasin{\psi}$      & $6^\circ$            & F.p.a., hdg. \\
    \hline
    $\epsfeasin{\Dot{Q}}$ & $10^{-2}\cdot\Dot{Q}_{\text{max}}$ & Heat rate \\
    $\epsfeasin{q}$ & $10^{-2}\cdot q_{\text{max}}$ & Dyn. press. \\
    $\epsfeasin{n}$ & $10^{-2}\cdot n_{g,\text{max}}$ & Normal load \\
    $\epsfeasin{\text{NFZ}}$ & $0.1^\circ$ & NFZ constraint \\
    \hline 
\end{tabular}
\captionof{table}{Convergence hyperparameters.}
\label{tab:ex1-mission-conv-hyperparams}

\begin{tabular}{lrl}
\hline 
    Parameter       & Value                 & Description      \\ 
    \hline                                          
    $s_x$  & $0.5$      & Primal state variable \\
    $s_u$  & $10$       & Primal control variable \\
    $s_{\lambda}$       & $0.1$     & Equality-constraint dual variable \\
    $s_{\mu}$           & $1$       & Inequality-constraint dual variable \\
    $\varepsilon_{\text{w,}\min}$  & $10^{-3}$       & Min. lower threshold for $(\Wh,\Wg)$ \\
    \hline 
\end{tabular}
\captionof{table}{Trust region step-sizes.}
\label{tab:ex1-mission-ptr-weights}
\end{minipage}

\subsection{Bank Angle Control Example} 


\subsubsection{Standalone {\autoscvx} }

We first present standalone numerical results for {\autoscvx} for an example where bank angle is modeled as the only control input. Angle-of-attack is assumed to adhere to the pre-designed velocity-dependent profile given in \cite{wang2018autonomous}. We solve the optimal control problem shown in Problem \ref{prob:ncvx-entry} with the mission parameters given in Table \ref{tab:ex1-mission-params}. The constraint penalty weight hyperparameters are tuned automatically within the {\autoscvx} framework. However, the trust region step sizes (on both the primal and dual variables), the optimality and feasibility convergence tolerances, and the minimum threshold for the quadratic weight update must still be selected. 
Larger trust region step sizes permit subproblem solutions farther from the given reference trajectory, resulting in fewer iterations to convergence. Tighter optimality tolerances enforce smaller optimization variable or cost variations for algorithm termination, while tighter feasibility tolerances enforce smaller constraint violation residuals for termination. As the quadratic penalty weights decay, it is desirable that these matrices remain positive definite; the weight update is lower bounded by an arbitrarily small threshold.
The values for these additional hyperparameters are displayed in Tables \ref{tab:ex1-mission-conv-hyperparams} and \ref{tab:ex1-mission-ptr-weights}. 

\begin{figure}[htb!]
    \centering
    \includegraphics[width=1\columnwidth,trim={0cm 0cm 0cm 0cm},clip]{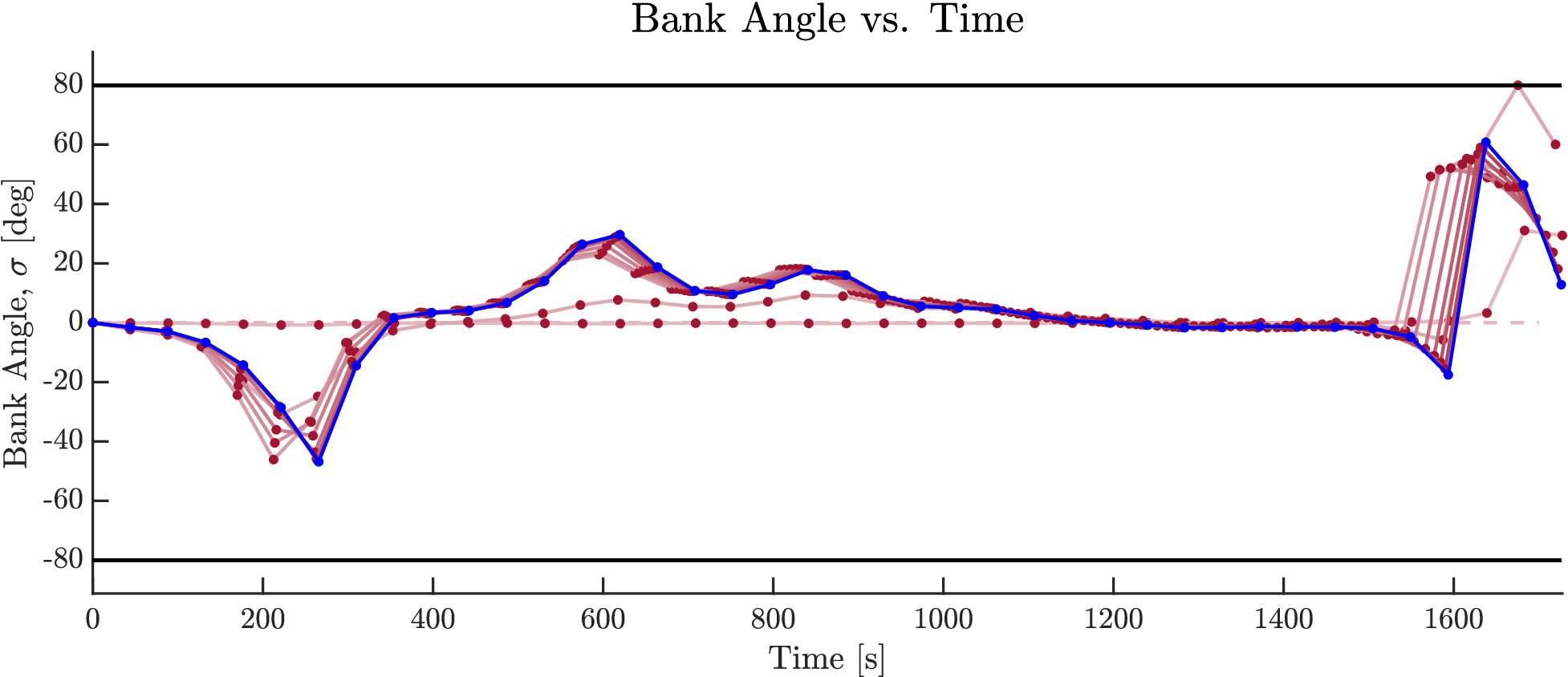}
    \caption{The bank angle solution for {\autoscvx} in Example A is shown converging to the optimal solution across iterations, the converged solution shown in blue.}
    \label{fig:ex1-ctrl}
\end{figure}

\begin{figure}[htb!]
    \centering
    \includegraphics[width=1\columnwidth,trim={1cm 0cm 1cm 0cm},clip]{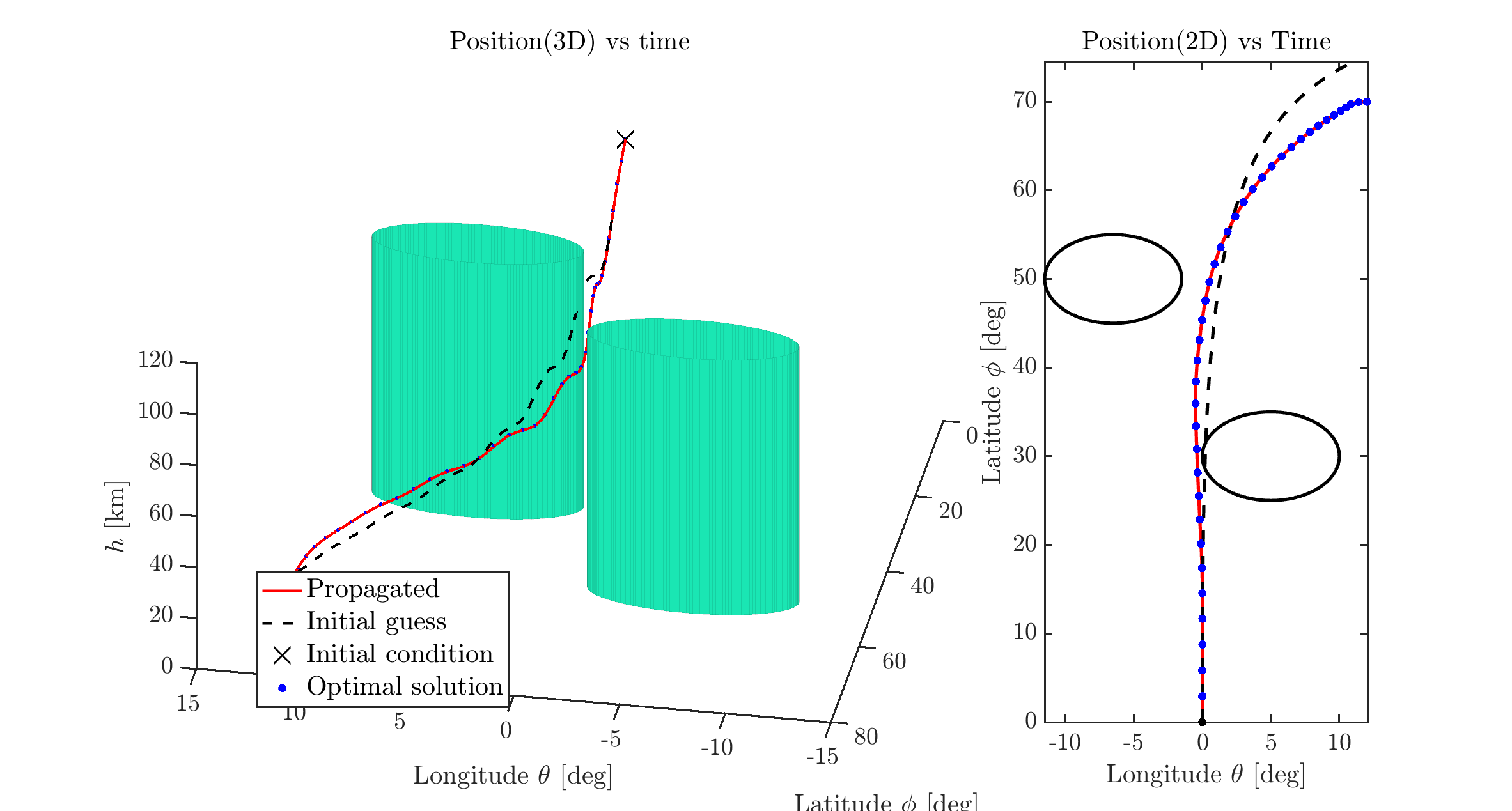}
    \caption{The optimal solution from {\autoscvx} for vehicle position is displayed for Example A. 
    }
    \label{fig:ex1-position}
\end{figure}

\begin{figure}[htb!]
    \centering
    \includegraphics[width=1\columnwidth,trim={0cm 0cm 0cm 0cm},clip]{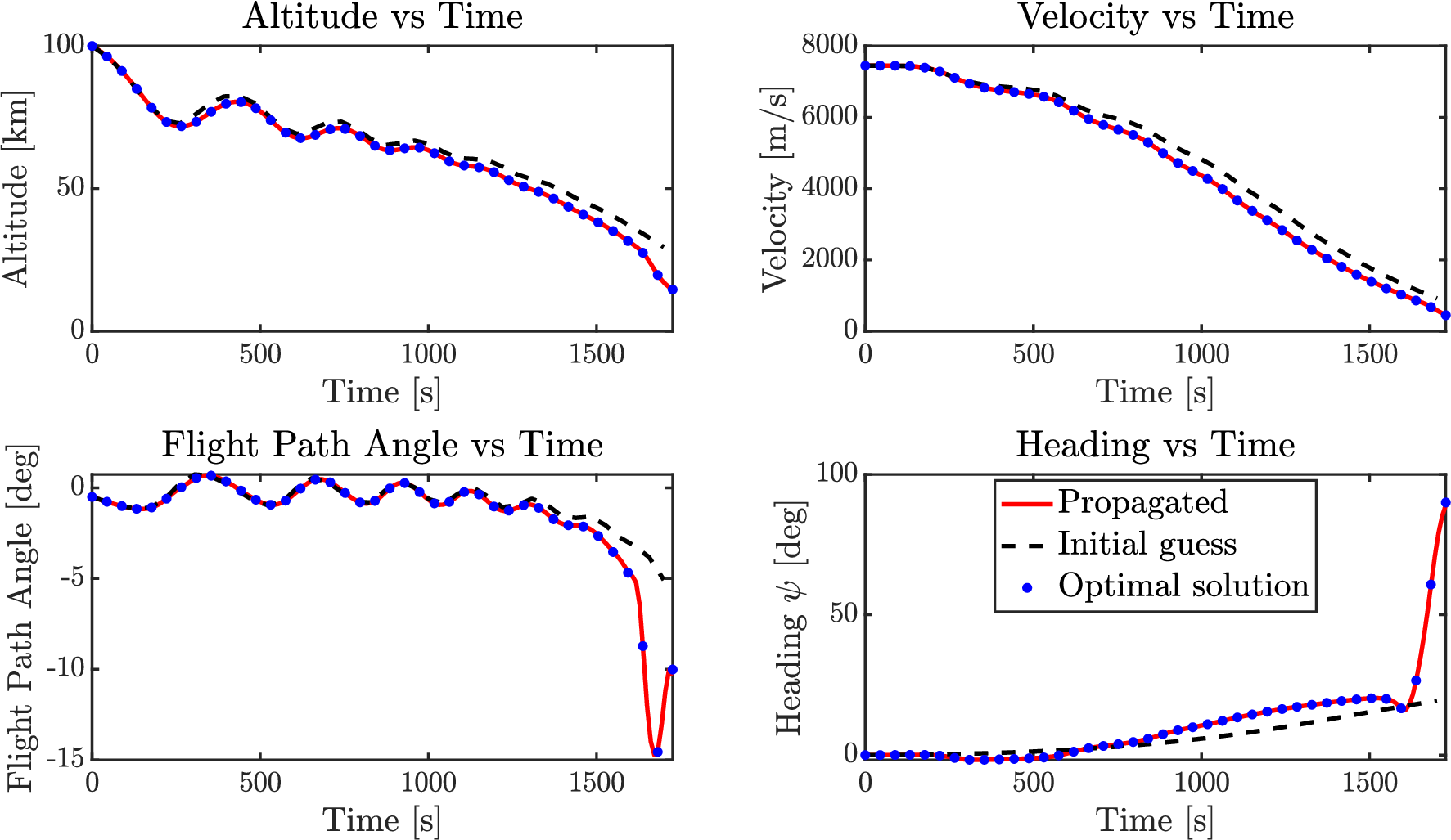}
    \caption{The discrete optimal solution for {\autoscvx} is displayed for vehicle altitude, velocity, flight path angle and heading are displayed for Example A. 
    }
    \label{fig:ex1-states}
\end{figure}

The bank angle control is shown converging to the optimal solution across {\autoscvx} iterations in Figure \ref{fig:ex1-ctrl}. This optimal control corresponds to an optimal state solution (the vehicle's position and velocity in polar coordinates). To validate that the optimal state and control are dynamically feasible, a single-shooting propagation of the system dynamics is performed via a first-order-hold of the optimal control, starting from the initial state boundary condition (given in Table \ref{tab:ex1-mission-params}) over the optimal time horizon solution $t_{F}\approx1714.93 \ [\text{s}]$. The converged optimal state solutions are shown to overlay this single-shooting propagation in Figures \ref{fig:ex1-position} and \ref{fig:ex1-states}, and compared against the initial guess. 

\begin{figure}[htb!]
    \centering
    \includegraphics[width=1\columnwidth,trim={2.5cm 0cm 2.5cm 0cm},clip]{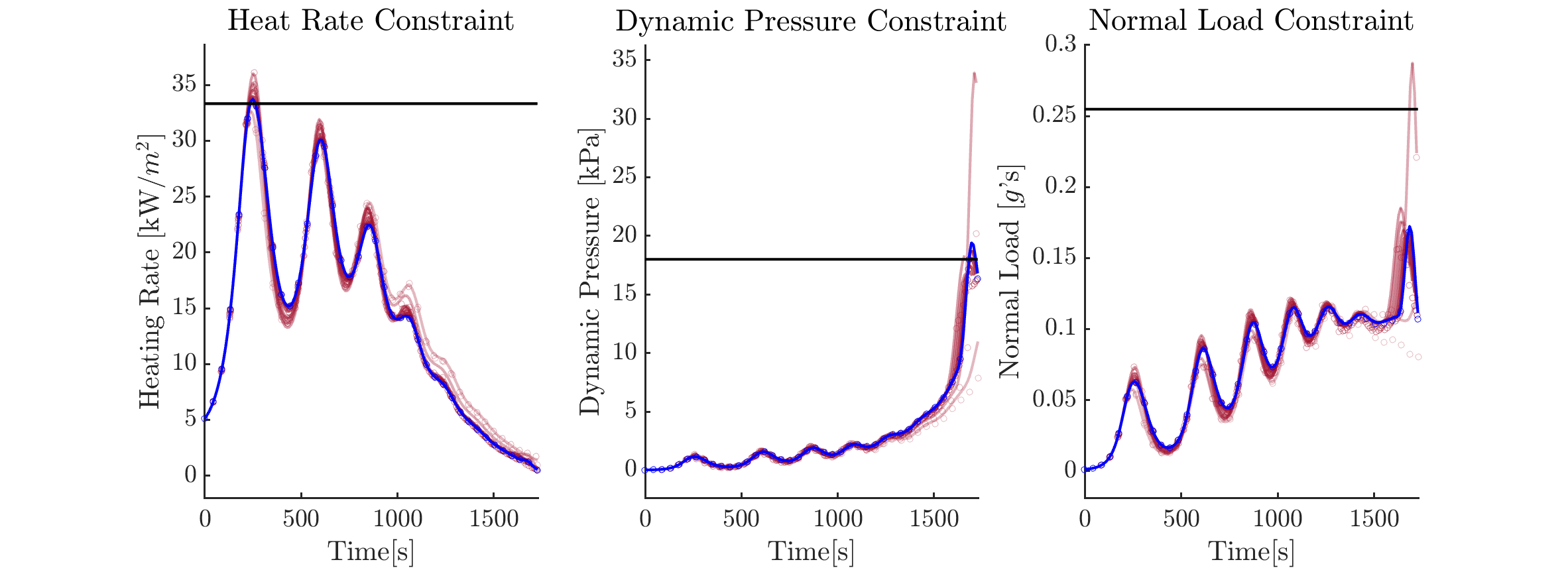}
    \caption{Path constraints are shown converging across {\autoscvx} iterations for Example A, the converged solution shown in blue.}
    \label{fig:ex1-cnst-conv}
\end{figure}

\begin{figure}[htb!]
    \centering
    \includegraphics[width=1\columnwidth,trim={0cm 0cm 0cm 0cm},clip]{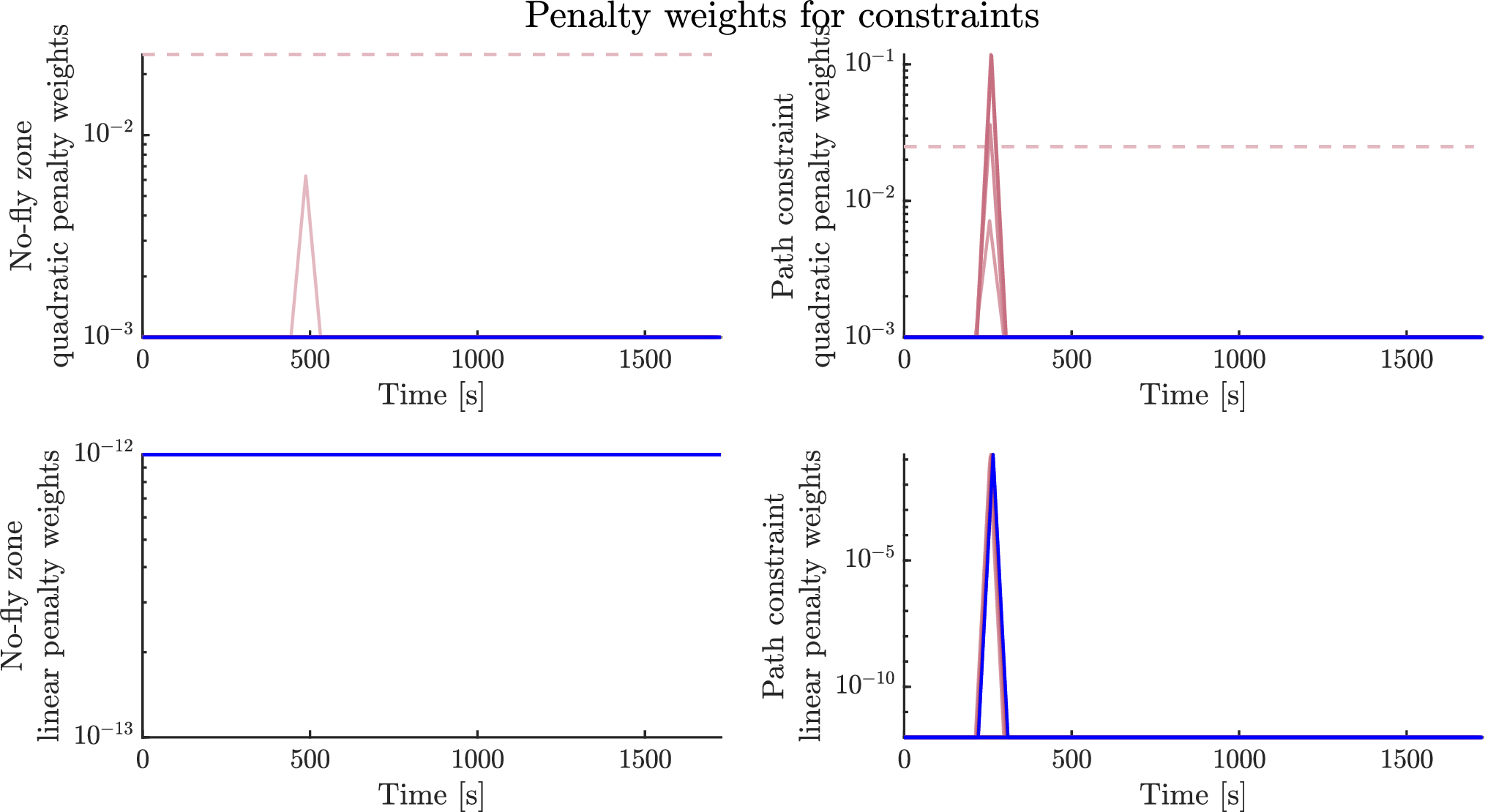}
    \caption{Quadratic and linear penalty weights ($\Wh,\Wg$ and $\lambda,\mu$, resp.) are shown converging across {\autoscvx} iterations as constraint violation evolves, the converged solution shown in blue.}
    \label{fig:ex1-hyperparams-conv}
\end{figure}

\begin{figure}[htb!]
    \centering
    \includegraphics[width=1\columnwidth,trim={0cm 0cm 0cm 0cm},clip]{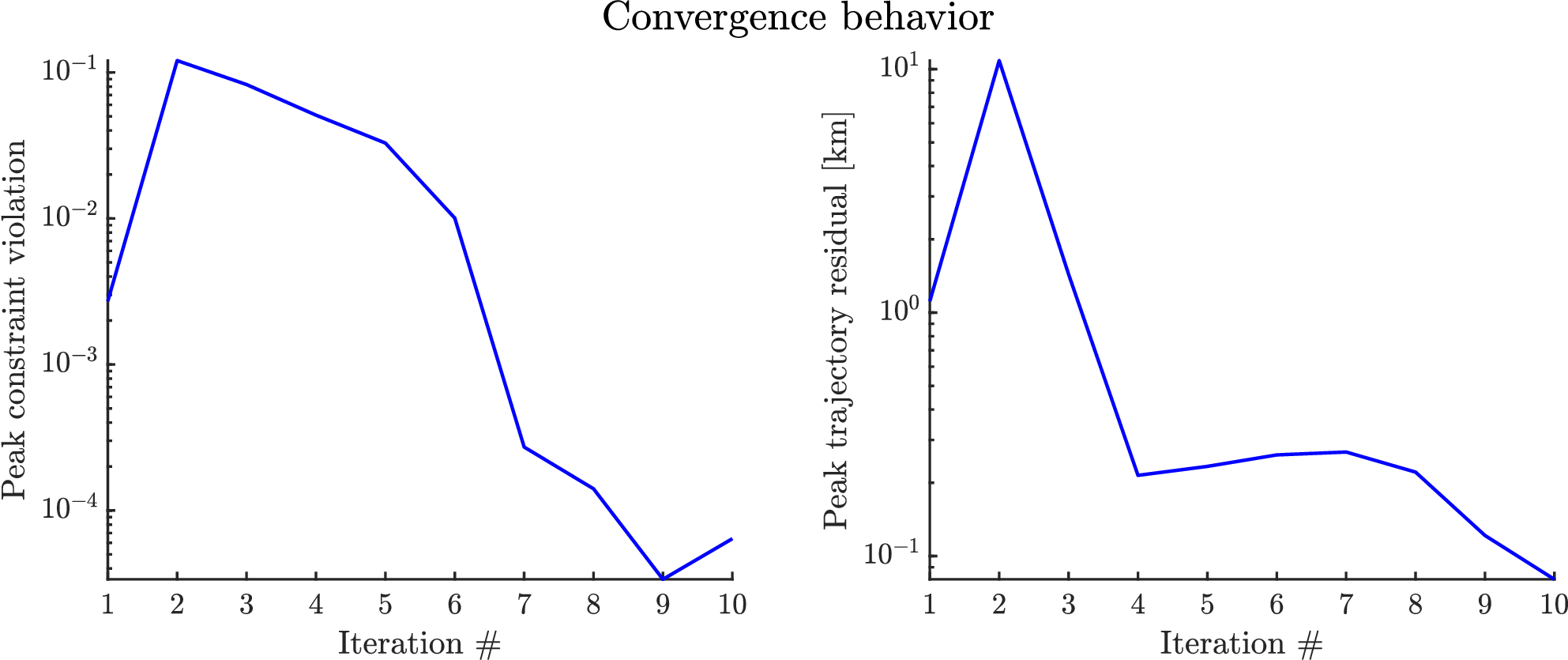}
    \caption{Nondimensionalized constraint violation is shown converging to the feasibility tolerance within {\autoscvx} for Example A.}
    \label{fig:ex1-violation}
\end{figure}

The constraint violation across {\autoscvx} iterations is shown in Figure \ref{fig:ex1-cnst-conv}. In the early iterations, when the changes in the state and control are large, the convex approximations of both the dynamics and path constraints have a large discrepancy when compared against the true nonconvex functions. In addition, the constraint violation (virtual buffers) have significant magnitude, indicating violation of the linearized constraints, as shown in Figure \ref{fig:ex1-violation}. 
These factors result in spikes in the constraint violation when evaluating the feasibility of both the subproblem and the single-shooting propagation trajectory in these intermediate iterations for all buffered constraints. When these buffered subproblem constraints are violated, the corresponding penalty weights for these problematic constraints and time indices begin to grow, as shown in Figure \ref{fig:ex1-hyperparams-conv}. Here, the penalty weights for all path constraints and no fly zone constraints over time are plotted together, as well as their convergence across iterations.

When the constraints become satisfied, the quadratic penalty weights begin to shrink due to the multiplicative update; they rapidly decay towards a minimum threshold bound (given in Table \ref{tab:ex1-mission-ptr-weights}) at all time indices where corresponding constraints become satisfied. If constraints remain tight, the corresponding quadratic penalty weights will remain constant. The linear penalty weights, on the other hand, integrate constraint violation across {\autoscvx} iterations; these linear penalties remain constant when constraints become either tight or feasible. 

\begin{figure}[htb!]
    \centering
    \includegraphics[width=1\columnwidth,trim={0cm 0cm 0cm 0cm},clip]{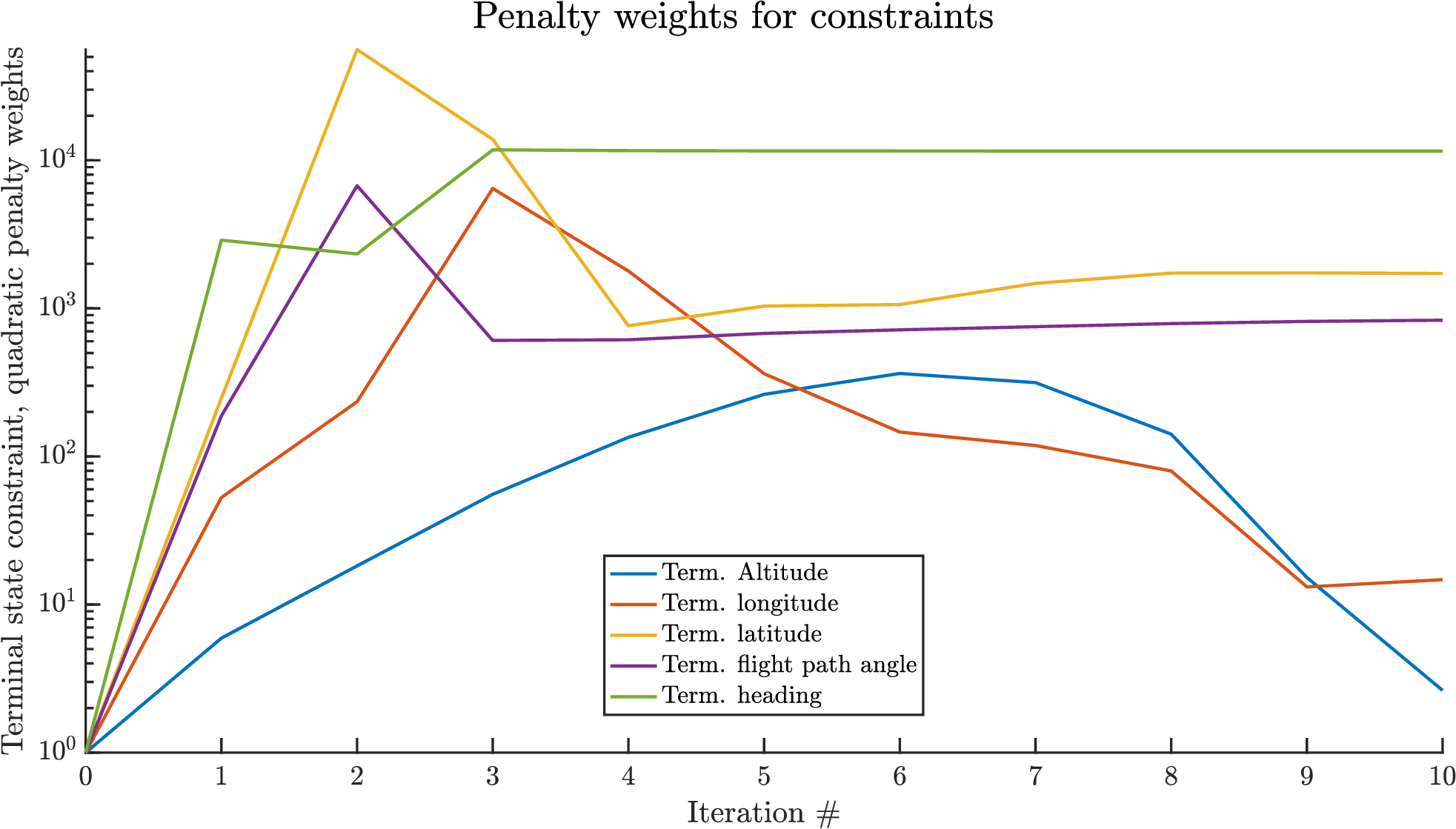}
    \caption{Each element of the terminal boundary condition quadratic penalty weight is shown converging across {\autoscvx} iterations.}
    \label{fig:ex1-hyperparams-term}
\end{figure}

As mentioned in Section \ref{sec:cvx-buff}, in this work the nonconvex dynamics are left unbuffered; 
for such a sensitive system, the presence of small but nonzero virtual control terms in the recursive LTV dynamics constraint imposed in the convex subproblem can result in large defects between the optimal state solution and the single-shooting propagation. 
Instead, the algorithm is supplied with a dynamically feasible initial guess found by propagating zero bank angle from the initial boundary condition across an initial guess of the time horizon. The convex terminal state constraints are supplied with virtual buffers, which are assigned to corresponding penalty weights, and the algorithm iterates until convergence to a solution satisfying the terminal boundary condition. The quadratic penalty weights for the terminal state boundary conditions are shown across across {\autoscvx} iterations in Figure \ref{fig:ex1-hyperparams-term}. 
These weights grow and decay alongside the terminal boundary condition violation, and converge to a large value as the equality constraint becomes tight, i.e. when the virtual buffers drop below the feasibility tolerance (given in Table \ref{tab:ex1-mission-conv-hyperparams}). 

\subsubsection{Comparison: {\autoscvx} vs. {\ptr} Solution} \label{subsec:ex1-compare}

Next, the converged solution of the {\autoscvx} algorithm is compared against solutions for the penalized trust region (\ptr) algorithm with various constant quadratic penalty weight selections. Typically, the quadratic penalty weights in the {\ptr} algorithm are tuned with a single fixed scalar value that remains constant across all constraints and time indices:
\begin{align}
    \Wh &= \omega ~\diag{\one{\nequ}}\succ 0 , \quad 
    \Wg = \omega ~\diag{\one{\nineq}}
    \succ 0, \text{ where } 
    1 \ll \omega \in \mathbb{R}^{++} , \\
    \lambar &= \one{\nequ}, \quad 
    \mubar = \one{\nineq}.
\end{align}
In prior work, these penalty weights were completely hand-tuned by a trial and error process, and typically guessed to be positive numbers of large magnitude with the hopes of coaxing the algorithm to converge to a feasible solution \cite{szmuk2018jgcdarxiv,szmuk2020successive,Reynolds2019,malyuta2021convex}. 
In order to improve performance, experienced researchers may also hand-tune a selection of weights (associated with different constraints) to have different relative magnitudes, essentially guessing at which constraints will become tight across the trajectory and push back most on a feasible solution. Although not algorithmic, when compared against {\autoscvx} it can be seen that selecting a hand-tuned set of weights for the {\ptr} algorithm is akin to guessing the infinity norm of the dual variables associated with each constraint. While the {\ptr} algorithm is able to produce feasible solutions for well-tuned penalty weights, the tuning process may be time-consuming, and the large quadratic weights across all constraints and time result in worsened problem conditioning. 
As a result, increases in the penalty weight order of magnitude may result in increased cost of the optimal solution returned by {\ptr}, discussed below.

\begin{figure}[htb!]
    \centering
    \includegraphics[width=1\columnwidth,trim={0cm 0cm 0cm 0cm},clip]{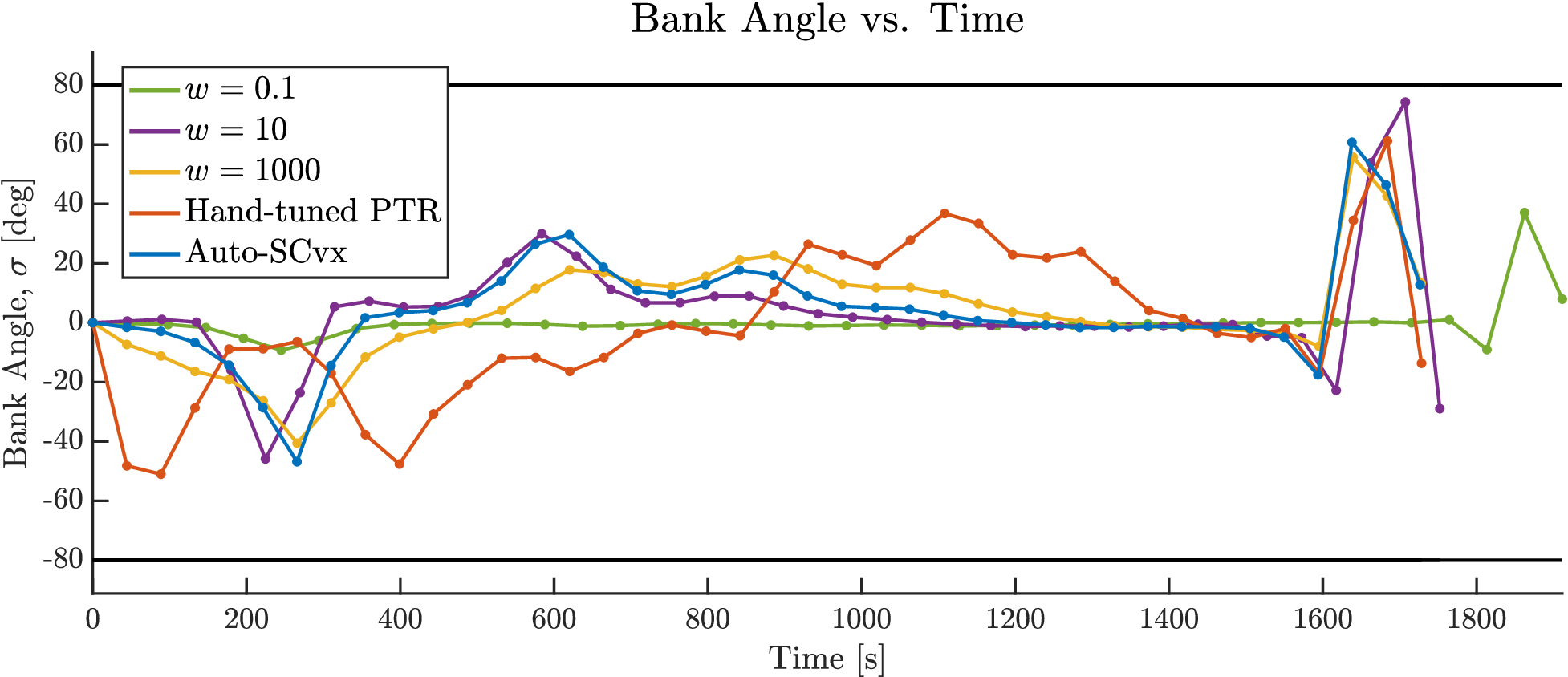}
    \caption{Optimal bank angle solution profiles for Example A are compared between {\autoscvx} and several {\ptr} runs, each with a different fixed-weight tuning.}
    \label{fig:ex1-ctrl-compare}
\end{figure}

\begin{figure}[htb!]
    \centering
    \includegraphics[width=1\columnwidth,trim={1cm 0cm 1cm 0cm},clip]{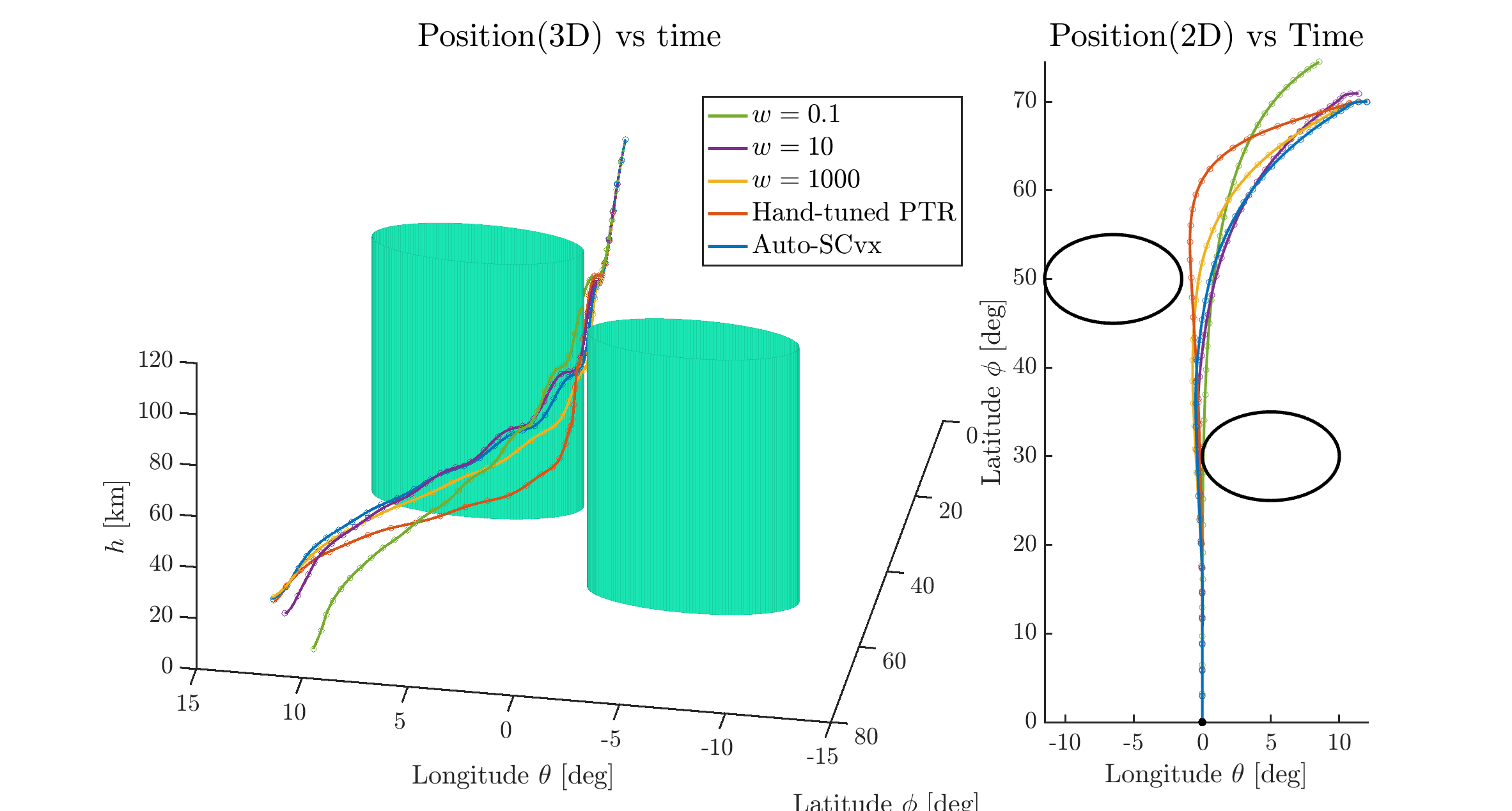}
    \caption{Position trajectory solutions for Example A are compared between {\autoscvx} and various fixed-weight {\ptr} runs.}
    \label{fig:ex1-position-compare}
\end{figure}

\begin{figure}[htb!]
    \centering
    \includegraphics[width=1\columnwidth,trim={0cm 0cm 0cm 0cm},clip]{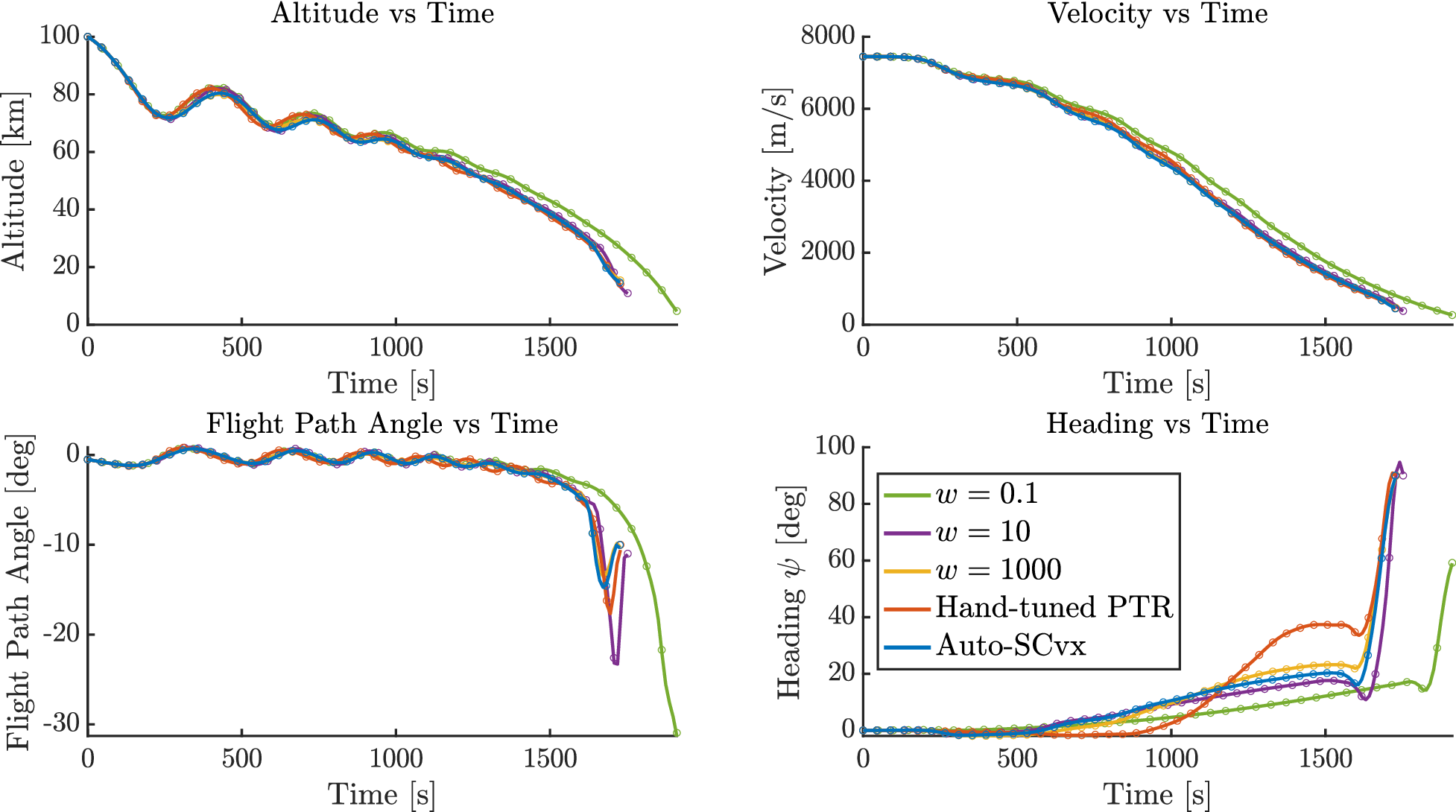}
    \caption{State solutions are compared between {\autoscvx}  and various fixed-weight {\ptr} runs for Example A.}
    \label{fig:ex1-state-compare}
\end{figure}

\begin{figure}[htb!]
    \centering
    \includegraphics[width=1\columnwidth,trim={0cm 0cm 0cm 0cm},clip]{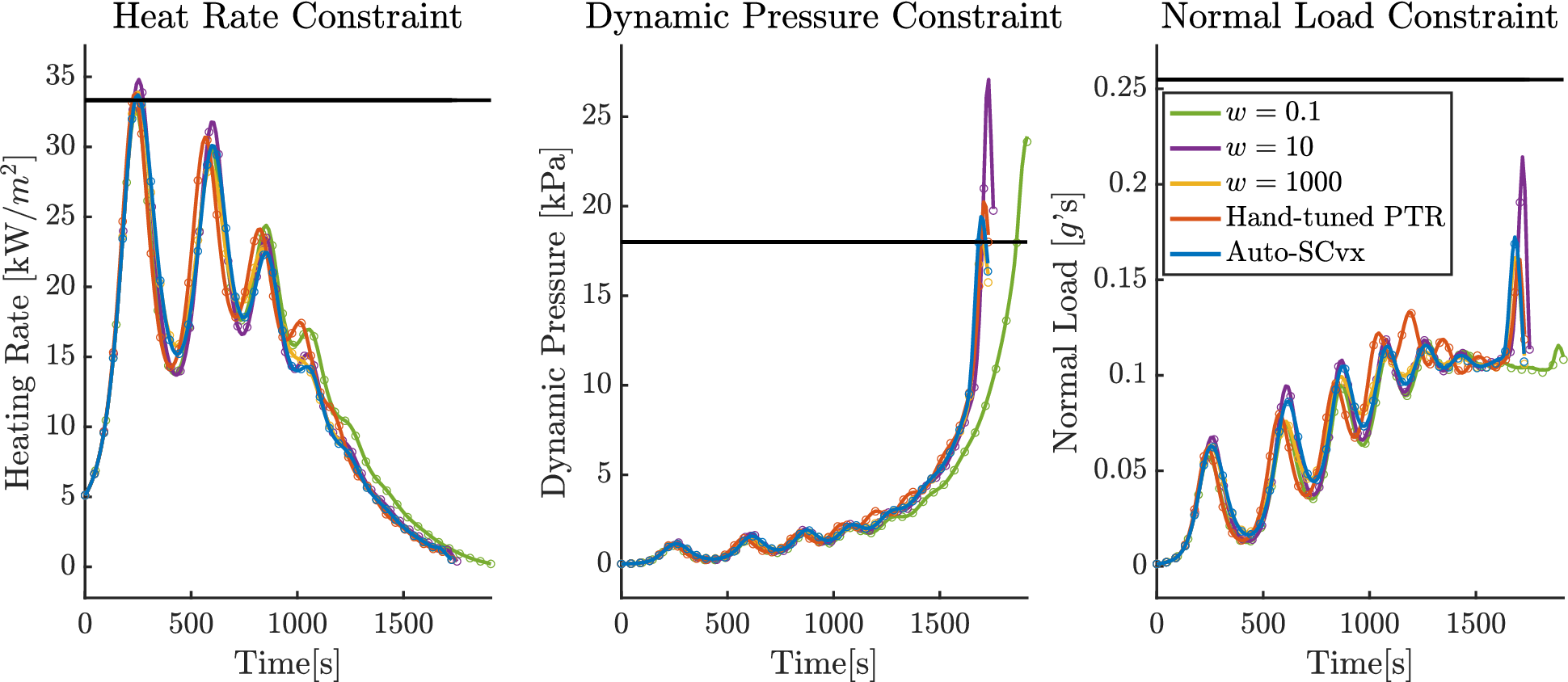}
    \caption{Path constraint solutions are compared between {\autoscvx} and {\ptr} methods with various fixed-weight tuning for Example A.}
    \label{fig:ex1-cnst-compare}
\end{figure}

The optimal bank angle control of the {\autoscvx} solution is compared against four different {\ptr} solutions in Figure \ref{fig:ex1-ctrl-compare}. Three {\ptr} solutions use a constant weight of increasing magnitude across a log scale, each normalized by the number of temporal nodes $N=40$ to reduce solution sensitivity to the number of variables. The fourth {\ptr} solution uses hand-tuned weights with different relative magnitudes for different constraint penalties, again normalized by the number of temporal nodes. The corresponding optimal state solution trajectories are shown in Figures \ref{fig:ex1-position-compare} and \ref{fig:ex1-state-compare}. The {\autoscvx} algorithm converges to a feasible solution in 10 iterations, i.e. terminates with a trajectory satisfying all constraints to the specified feasibility tolerance. In contrast, the {\ptr} algorithms with penalty weights of insufficient magnitude never converge to a feasible solution. Both the terminal state boundary condition (given in Table \ref{tab:ex1-mission-params}) and the dynamic pressure constraint (displayed in Figure \ref{fig:ex1-cnst-compare}) remain violated. Due to this, these algorithm instances do not converge beneath the feasibility tolerance but instead terminate without a solution at the specified maximum iteration limit (20 iterations). 

\begin{figure}[htb!]
    \centering
    \includegraphics[width=1\columnwidth,trim={0cm 0cm 0cm 0cm},clip]{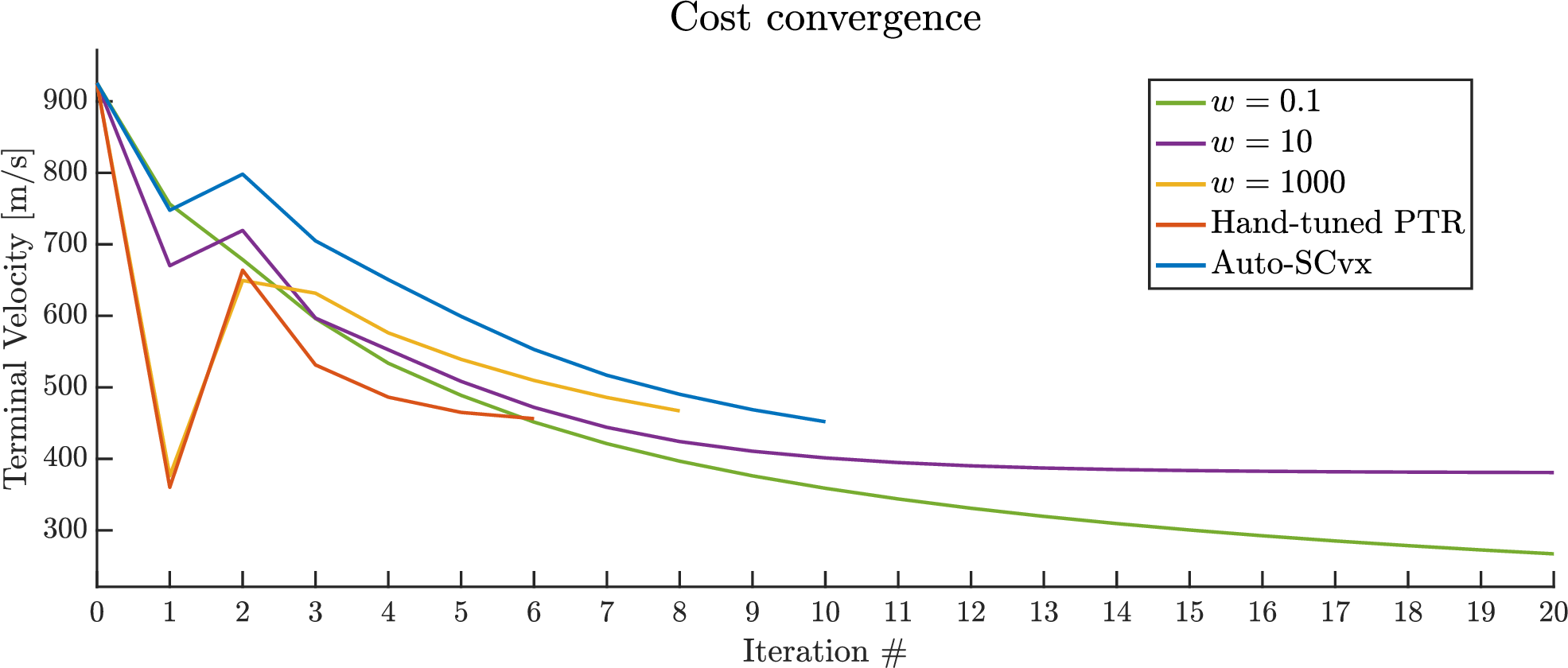}
    \caption{A comparison of cost convergence across iterations between {\autoscvx} and {\ptr} with various weight-tuning is shown for Example A.
    Note that the {\ptr} instances without sufficiently large penalty weights never achieve feasibility, and terminate at a maximum iteration limit.
    }
    \label{fig:ex1-cost-compare}
\end{figure}

\begin{figure}[htb!]
    \centering
    \includegraphics[width=1\columnwidth,trim={0cm 0cm 0cm 0cm},clip]{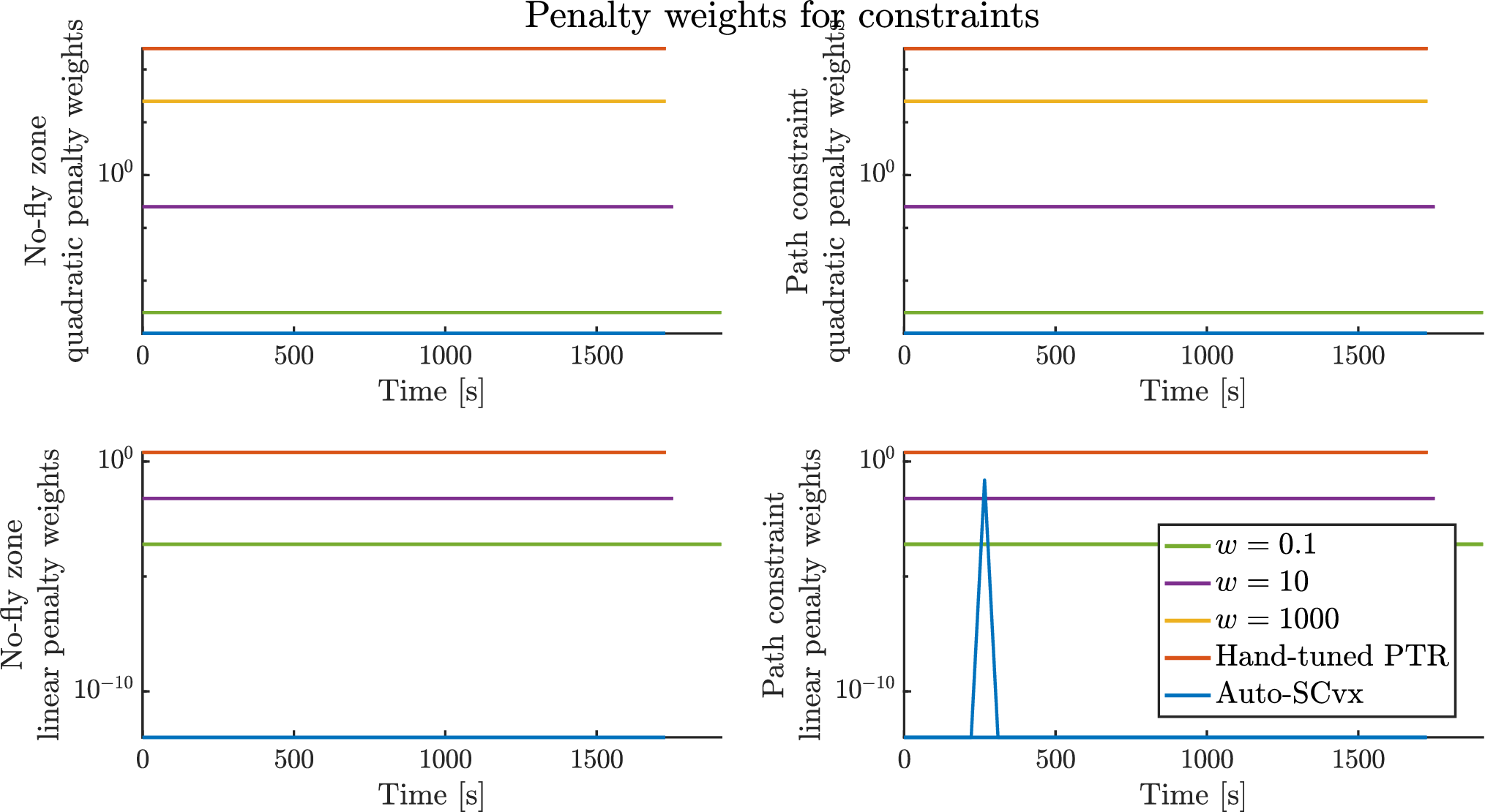}
    \caption{Hyperparameters at convergence are compared between {\autoscvx} and {\ptr} implementations for Example A.}
    \label{fig:ex1-hyperparams-compare}
\end{figure}

\begin{figure}[htb!]
    \centering
    \includegraphics[width=1\columnwidth,trim={0cm 0cm 0cm 0cm},clip]{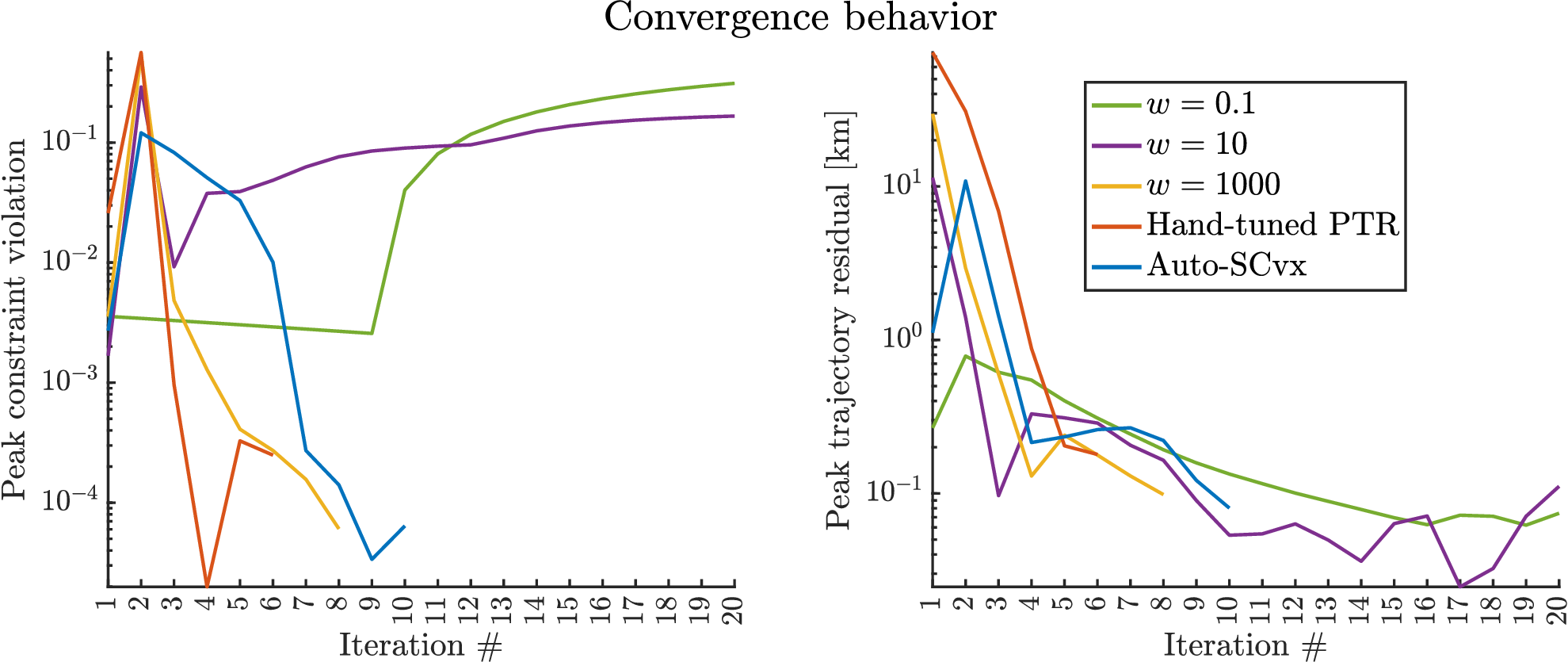}
    \caption{A comparison of nondimensionalized constraint violation convergence behavior between {\autoscvx} and {\ptr} for Example A. 
    Note that dynamic defect is decreasing as $\dx$ becomes small.
    }
    \label{fig:ex1-violation-compare}
\end{figure}

The {\ptr} algorithms for both constant weight $\omega=1000$ and the hand-tuned weights (shown in Table \ref{tab:ex1-mission-ptr-weights}) converge to a feasible solution. Inspection of the cost convergence in Figure \ref{fig:ex1-cost-compare} shows that these solutions are similar to the solution returned by {\autoscvx}. Numerically, it is shown that {\autoscvx} achieves the best performance, with a terminal velocity cost  of approximately $v_f\approx451.88$m/s, while the convergent {\ptr} algorithms return approximately and $v_f\approx469.04$m/s for the $\omega=1000$ case and $v_f\approx456.22$m/s for the hand-tuned case. 
Intuitively, this makes sense: if the quadratic penalty weights for all constraints at every time index have large magnitude relative to the gradient of the true cost, then the problem conditioning will suffer. 
Because {\autoscvx} harnesses dual variable information to update these penalty weights, the penalty terms with large magnitudes can be targeted at specific constraints and time indices along the trajectory required to achieve feasibility. 
As the algorithm converges to the feasible region of the nonconvex constraints, these quadratic penalty weights decay towards zero and the solver can begin to prioritize optimality. Converged hyperparameters for {\autoscvx} are compared against the static weights of {\ptr} in Figure \ref{fig:ex1-hyperparams-compare}, and the convergence of normalized constraint residuals are displayed in Figure \ref{fig:ex1-violation-compare}. Statistical results for ten different problems, each with a dispersed initial boundary condition, are displayed in Table \ref{tab:ex1-autoscvx-compare-10-mc}.

\begin{table}[htb!]
\centering
\small
\begin{tabular}{lrrrrr}
\hline
Mean Statistics: & \multicolumn{5}{c}{Method:} \\\cmidrule{2-6}
 & \multicolumn{1}{l}{\begin{tabular}[c]{@{}l@{}} \ptr, $w = 0.1$ \end{tabular}}
 & \multicolumn{1}{l}{\begin{tabular}[c]{@{}l@{}} \ptr, $w = 10$ \end{tabular}}
 & \multicolumn{1}{l}{\begin{tabular}[c]{@{}l@{}} \ptr, $w = 1000$ \end{tabular}}
 & \multicolumn{1}{l}{\begin{tabular}[c]{@{}l@{}} \ptr, Hand-tuned  \end{tabular}}
 & \multicolumn{1}{l}{\autoscvx} \\ \hline
\textbf{\# Iterations}                     & 
    20 & 
    20 & 
    9.2 & 
    6.7 & 
    10.9 \\ \hline
\textbf{Solve time per iteration {[}ms{]}} & 
    20.63 & 
    43.29 & 
    104.11 & 
    223.30 & 
    141.88 \\ \hline
\textbf{\begin{tabular}[c]{@{}l@{}} Cost [m/s] \\ (Terminal velocity) \end{tabular}} &
  267.00 &
  380.76 &
  461.40 &
  458.45 &
  449.32 \\ \hline
\textbf{\begin{tabular}[c]{@{}l@{}}Non-dimensionalized \\ constraint violation residual\end{tabular}} &
  0.31 &
  0.17 &
  4.7e-3 &
  4.2e-4 &
  3.9e-3 \\ \hline
\textbf{Convergence \% (across all runs)} & 
0 \% 
& 0 \% 
& 90 \% 
& 100 \% 
& 100 \% \\ \hline
\end{tabular}
  \caption{Mean performance metrics of {\ptr} with static weights vs. \texttt{\autoscvx}.}
  \label{tab:ex1-autoscvx-compare-10-mc}
\end{table}

\subsection{Bank Angle and Angle-of-attack Control Example}


In this section, we explore an example where both bank angle and angle-of-attack are modeled as control inputs for the vehicle. In this model, the control constraints bounding angle-of-attack magnitude are nonconvex functions of vehicle velocity as shown in Equation \eqref{eq:bankaoa-cnst}. Here, the nominal design profile from Example A has been relaxed into an inequality constraint by a margin of $\pm 5^\circ$, with hard limits on the angle-of-attack magnitude between $[0^\circ,40^\circ]$. As a result, in this example the lift-to-drag ratio (i.e. lift and drag coefficients) are directly modulated by the control input. In order to allow {\ptr} to converge, two changes were made to the problem parameters. First, the terminal altitude constraint was relaxed into an inequality, shown in Table \ref{tab:ex1-mission-params}. Second, the convergence criteria for altitude and velocity perturbations were relaxed, shown in Table \ref{tab:ex1-mission-conv-hyperparams}.

\begin{figure}[htb!]
    \centering
    \includegraphics[width=1\columnwidth,trim={0cm 0cm 0cm 0cm},clip]{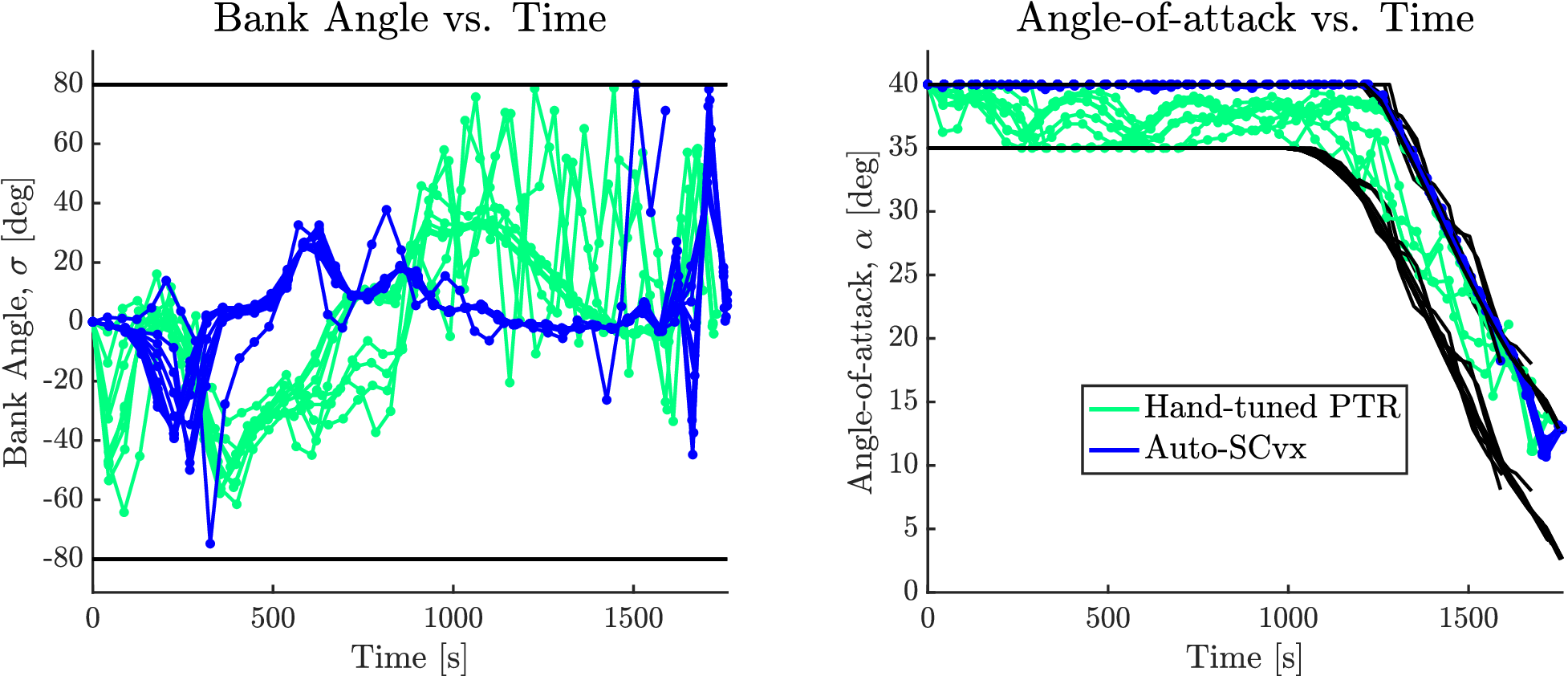}
    \caption{Optimal control solutions for bank angle and angle-of-attack are compared between {\autoscvx} and {\ptr} for Example B across dispersed initial position.
    Because the angle-of-attack bounds are velocity-dependent, a different set of control limits are present for each dispersed run. The {\autoscvx} angle-of-attack solution is seen to become tight against the upper bound, while the hand-tuned {\ptr} algorithm chatters between the feasible bounds.
    }
    \label{fig:ex2-ctrl}
\end{figure}

Converged solutions from the {\autoscvx} algorithm are compared against the highest-performing {\ptr} algorithm from Example A (Section \ref{subsec:ex1-compare}) for 10 different problem statements, where each initial boundary condition has been dispersed. Optimal control profiles for both bank angle and angle-of-attack are shown in Figure \ref{fig:ex2-ctrl}. Across dispersed problem statements, the optimal control profiles for {\autoscvx} are smoother and more consistent in comparison to the {\ptr} algorithm. Notably, the {\ptr} algorithm angle-of-attack profile oscillates wildly within the feasible region, while rarely touching the constraint boundary. In contrast, the {\autoscvx} solution remains tight for most of the trajectory at the highest feasible angle-of-attack, which corresponds to the highest lift-to-drag ratio. 
Intuitively, increasing the lift-to-drag ratio would allow increased control authority over decelerating and reducing vehicle velocity. It is interesting to note that there is overall less numerical chatter in the {\autoscvx} optimal control solutions, even when observing the bank angle profiles.

\begin{figure}[htb!]
    \centering
    \includegraphics[width=1\columnwidth,trim={1cm 0cm 1cm 0cm},clip]{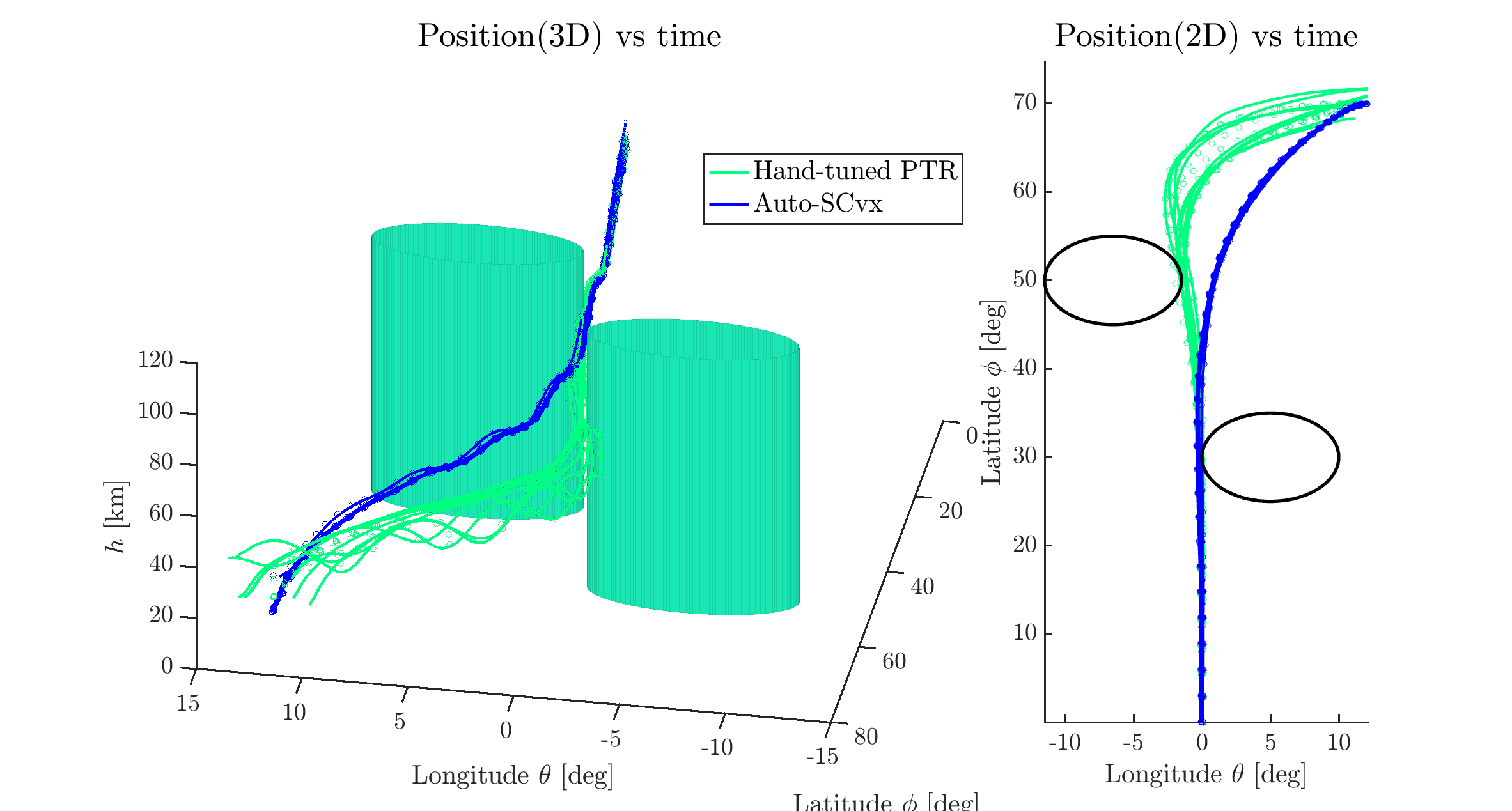}
    \caption{Dispersed optimal position trajectories are displayed comparing {\autoscvx} and {\ptr}  for Example B.}
    \label{fig:ex2-position}
\end{figure}

\begin{figure}[htb!]
    \centering
    \includegraphics[width=1\columnwidth,trim={0cm 0cm 0cm 0cm},clip]{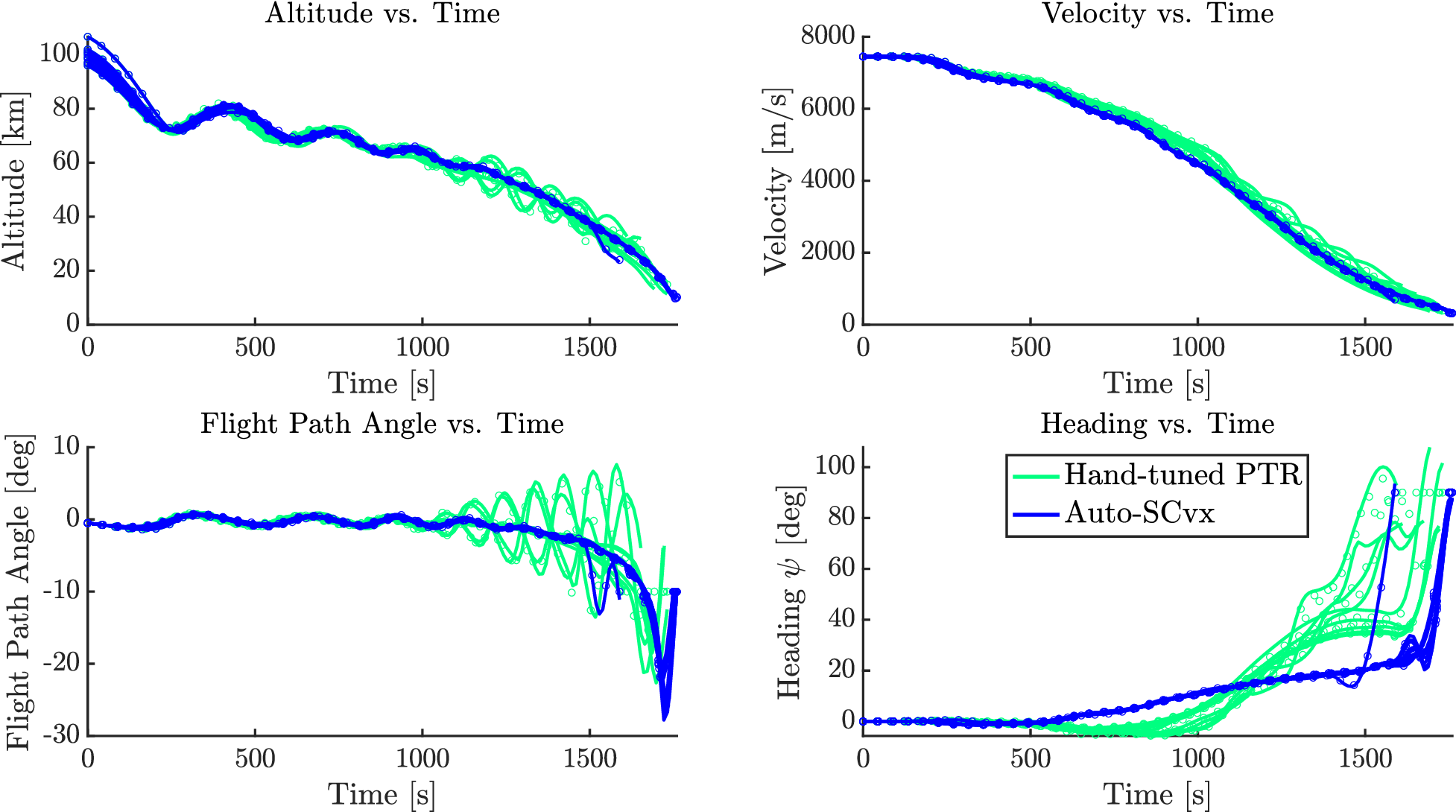}
    \caption{Optimal state trajectories are shown for Example B across dispersed initial condition between {\autoscvx} and \ptr.}
    \label{fig:ex2-states}
\end{figure}

\begin{figure}[htb!]
    \centering
    \includegraphics[width=1\columnwidth,trim={0cm 0cm 0cm 0cm},clip]{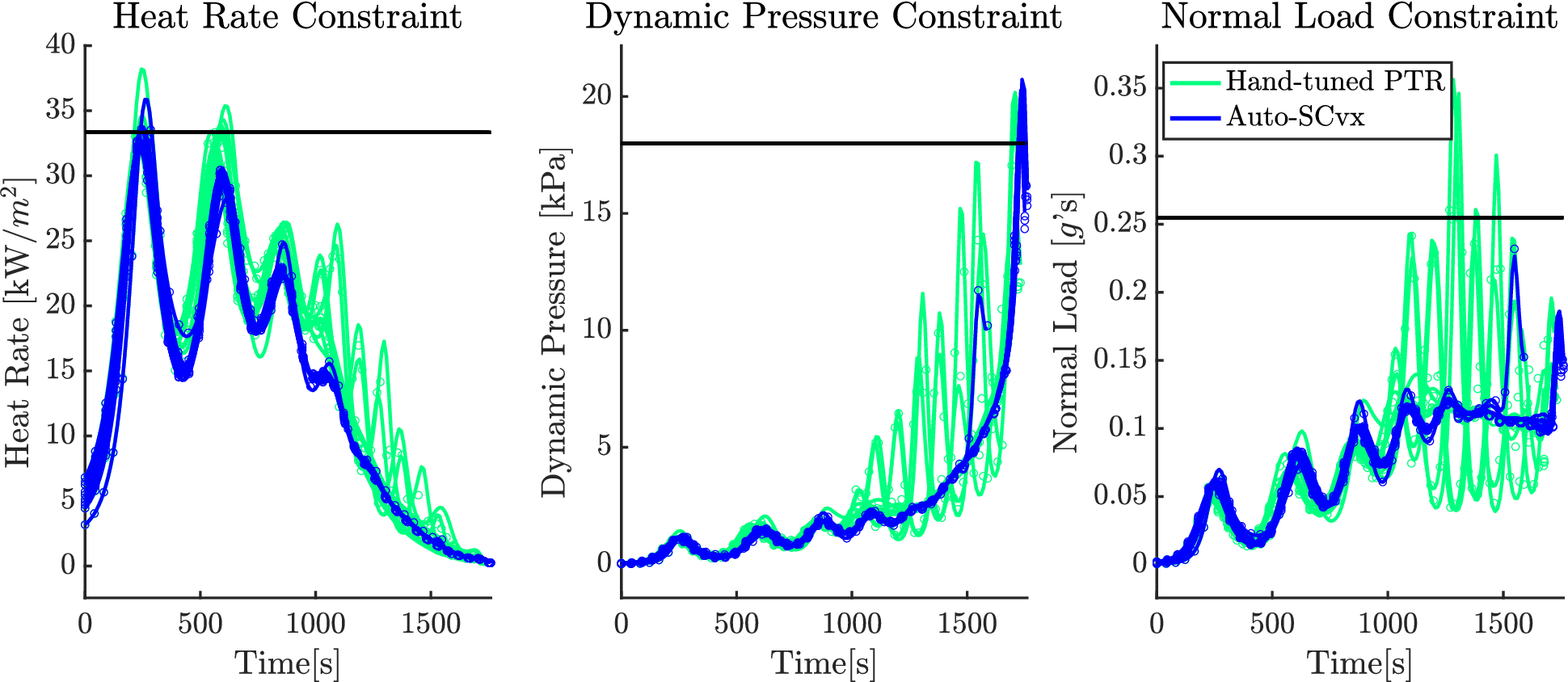}
    \caption{Path constraints are displayed for Example B for dispersed initial guess comparing {\autoscvx} and {\ptr}.}
    \label{fig:ex2-cnst}
\end{figure}

\begin{figure}[htb!]
    \centering
    \includegraphics[width=1\columnwidth,trim={0cm 0cm 0cm 0cm},clip]{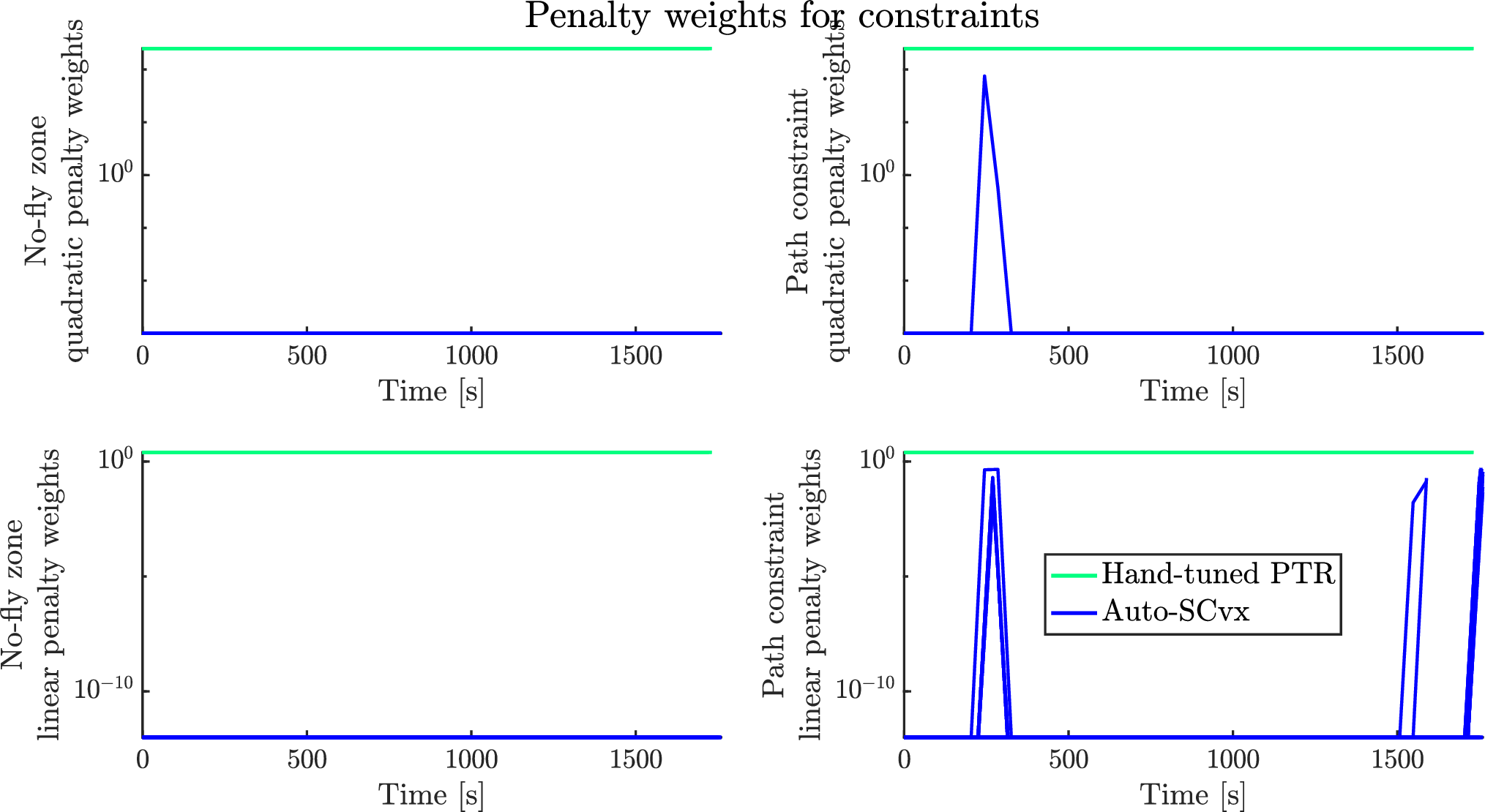}
    \caption{Converged quadratic and linear penalty weights are compared between {\autoscvx} and {\ptr} for dispersed solutions. 
    Note that {\ptr} has constant, hand-tuned quadratic penalty weights.
    }
    \label{fig:ex2-hyperparams}
\end{figure}

\begin{figure}[htb!]
    \centering
    \includegraphics[width=1\columnwidth,trim={0cm 0cm 0cm 0cm},clip]{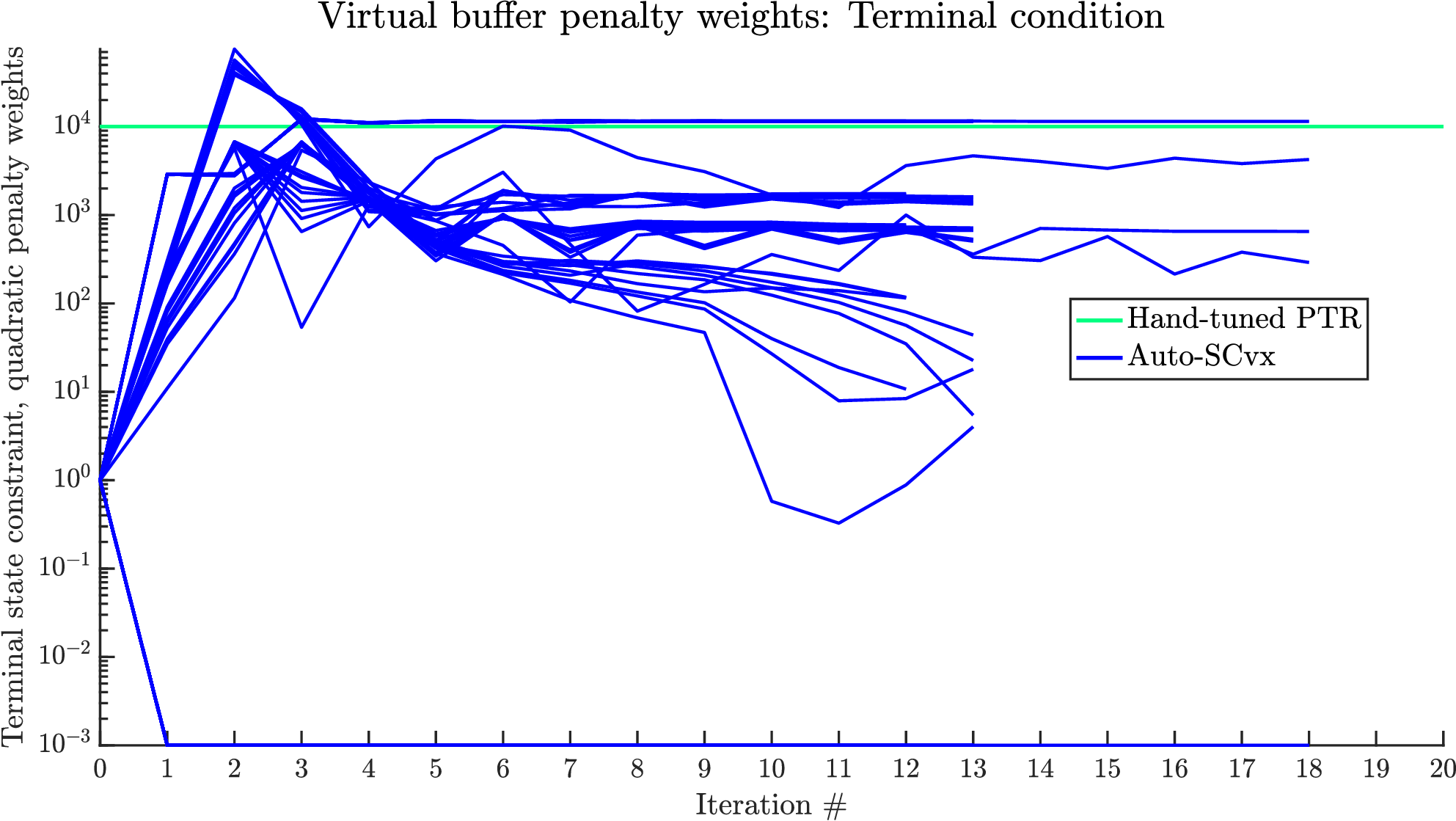}
    \caption{Quadratic penalty weights for the terminal boundary condition are displayed between {\autoscvx} and {\ptr} across dispersed solutions for Example B.}
    \label{fig:ex2-hyperparams-term}
\end{figure}

The corresponding optimal position trajectories are compared in Figure \ref{fig:ex2-position}. Once again, it can be observed that the family of solutions from {\autoscvx} are smoother and more consistent than the counterpart solutions from \ptr; in addition, the PTR solutions are shown to veer aggressively to the west to perform a prolonged side-slip maneuver seemingly to adjust heading to the terminal condition bound. In the process, a large intersample constraint violation occurs for one of the no-fly zone constraints. Solutions for all remaining vehicle states are displayed in Figure \ref{fig:ex2-states}.
Nonconvex constraints are shown in Figure \ref{fig:ex2-cnst}. Intersample violation is a common weakness of direct methods when optimizing trajectories in the presence of nonconvex constraints. Such intersample violations are visible for both algorithms here. However, the constraints themselves have a less dramatic oscillatory magnitude for the {\autoscvx} solution, likely due to the smoother control profile; the frequency and severity of intersample violation for the heating rate and normal load constraints seems reduced. The dynamic pressure constraint exhibits intersample violation for both algorithms. 
Corresponding penalty weights between the two algorithms are displayed in Figures \ref{fig:ex2-hyperparams} at the final iteration. Convergence behavior of the penalty weights associated with terminal state boundary condition violation are displayed in Figure \ref{fig:ex2-hyperparams-term}.

\begin{figure}[htb!]
    \centering
    \includegraphics[width=1\columnwidth,trim={0cm 0cm 0cm 0cm},clip]{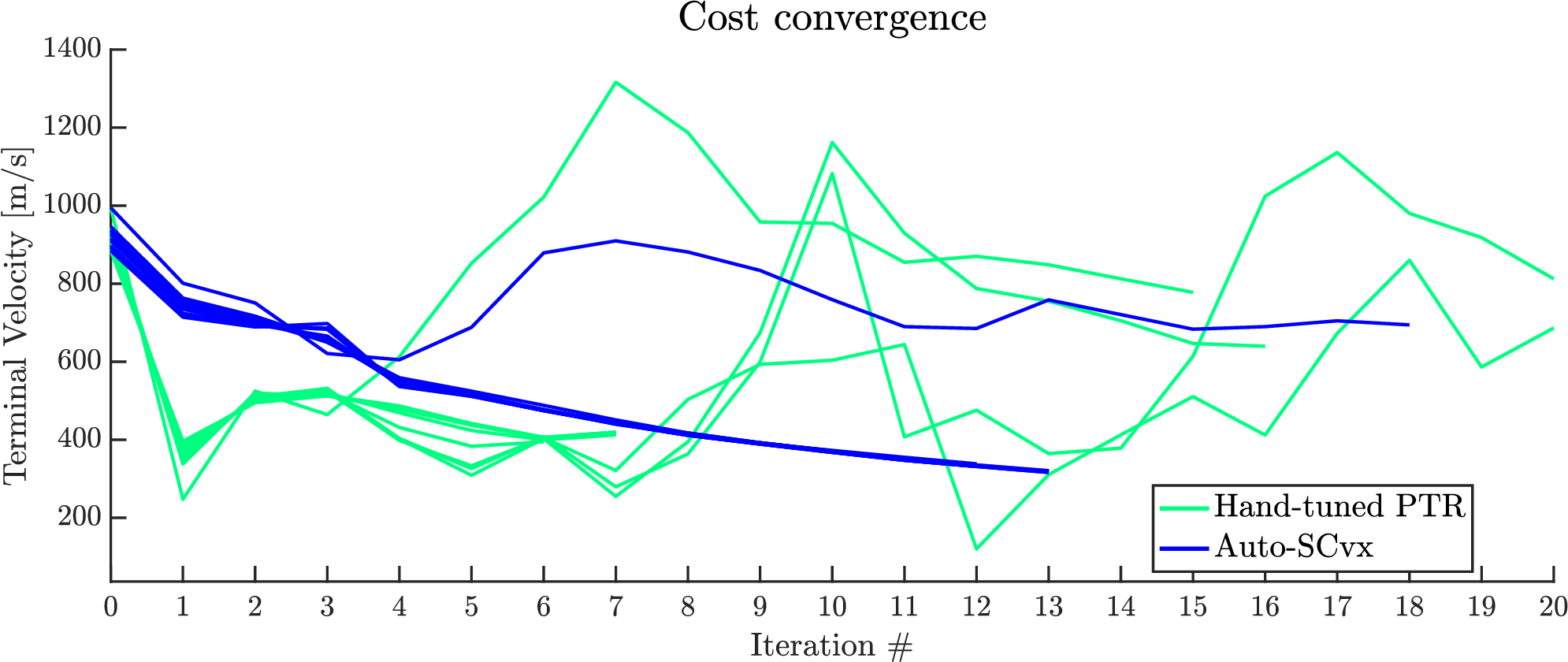}
    \caption{Cost convergence behavior is compared across dispersed solution between {\autoscvx} and {\ptr} for Example B. The important value is the cost at the final iteration.
    Note that the {\autoscvx} most often achieves a more optimal cost and converges in all cases, while {\ptr} converges in 80\% of the cases.
    }
    \label{fig:ex2-cost}
\end{figure}

\begin{figure}[htb!]
    \centering
    \includegraphics[width=1\columnwidth,trim={0cm 0cm 0cm 0cm},clip]{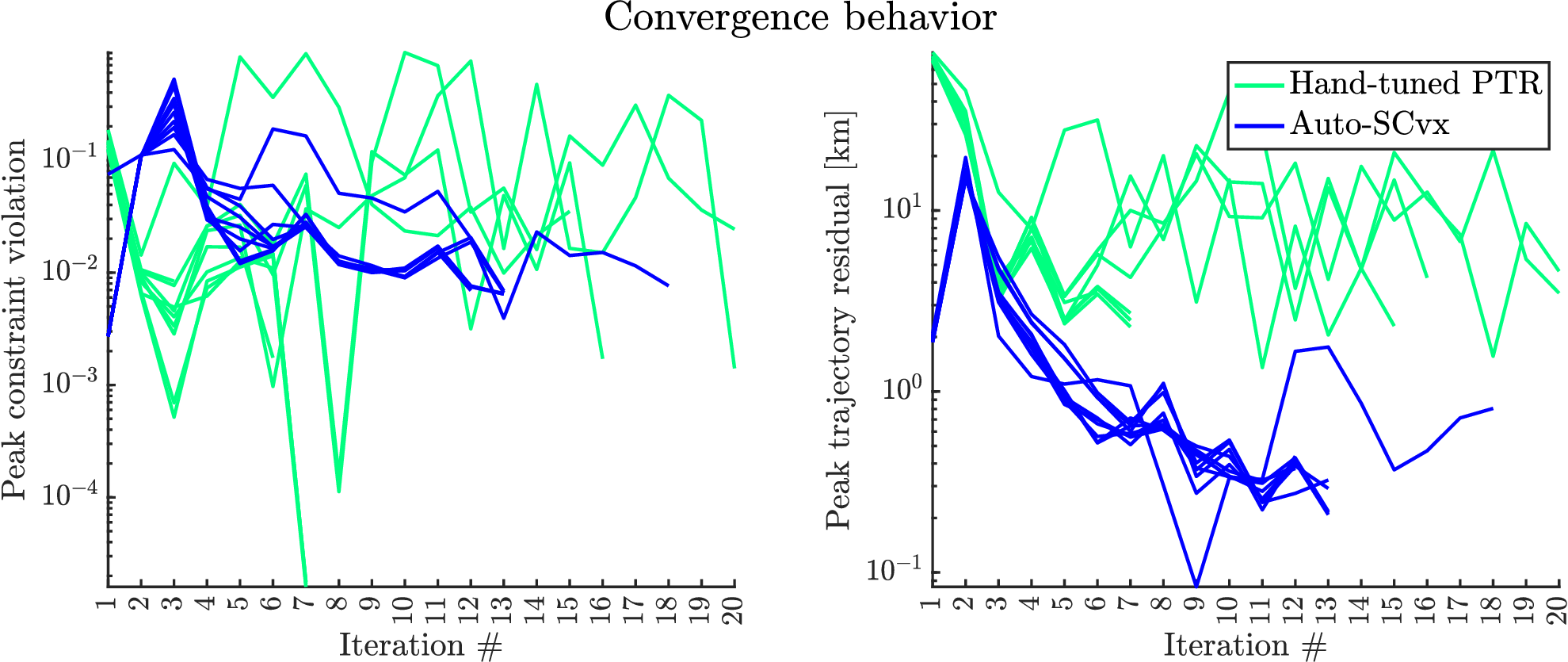}
    \caption{Nondimensionalized constraint residual convergence behavior is displayed across dispersed solutions between {\autoscvx} and {\ptr} for Example B. 
    }
    \label{fig:ex2-violation}
\end{figure}

Finally, comparisons of the optimality between solutions are displayed in Figure \ref{fig:ex2-cost}, and feasibility convergence behavior comparisons are displayed in Figure \ref{fig:ex2-violation}. In almost all cases, {\autoscvx} converges to a feasible solution with a lower terminal velocity cost than \ptr. In addition, {\autoscvx} appears to average lower constraint residuals at convergence. {\ptr} converges to a feasible solution for 8 of the 10 problems, but with a higher terminal velocity cost. For two of the problems, {\ptr} fails to converge to a feasible solution at all, and instead terminates at the max iteration bound (20 iterations). Mean statistics between the two algorithms are displayed in Table \ref{tab:ex2-autoscvx-compare-10-mc}.

\begin{table}[htb!]
    \centering
    \begin{tabular}{lrr}
        \hline
        Mean Statistics: & \multicolumn{2}{c}{Method:} \\\cmidrule{2-3}
         &  \multicolumn{1}{l}{\ptr, Hand-tuned}
         &   \multicolumn{1}{l}{ \textbf{\autoscvx} } \\ \hline
        \textbf{\# Iterations} & 
        10.8 & 
        13.1         \\ \hline
        \textbf{Solve time per iteration {[}ms{]}} & 
        419.64 & 
        152.96 \\ \hline
        \textbf{\begin{tabular}[c]{@{}l@{}} Cost [m/s] \\ (Terminal velocity) \end{tabular}} &
        549.05 &
        362.8 \\ \hline
        \textbf{\begin{tabular}[c]{@{}l@{}}Non-dimensionalized \\ constraint violation residual\end{tabular}} &
        7.3e-3 &
        6.9e-3 \\ \hline
        \textbf{Convergence \% (across all runs)} & 
        80 \%  & 
        100 \%         \\ \hline
    \end{tabular}
    \caption{Mean performance metrics of hand tuned weights and \autoscvx.}
    \label{tab:ex2-autoscvx-compare-10-mc}
\end{table}

\subsection{Bank Angle Dispersion Study}

\begin{figure}[htb!]
    \begin{minipage}{.5\textwidth}
    \centering
    \includegraphics[width=1\columnwidth,trim={0cm 0cm 12cm 2cm},clip]{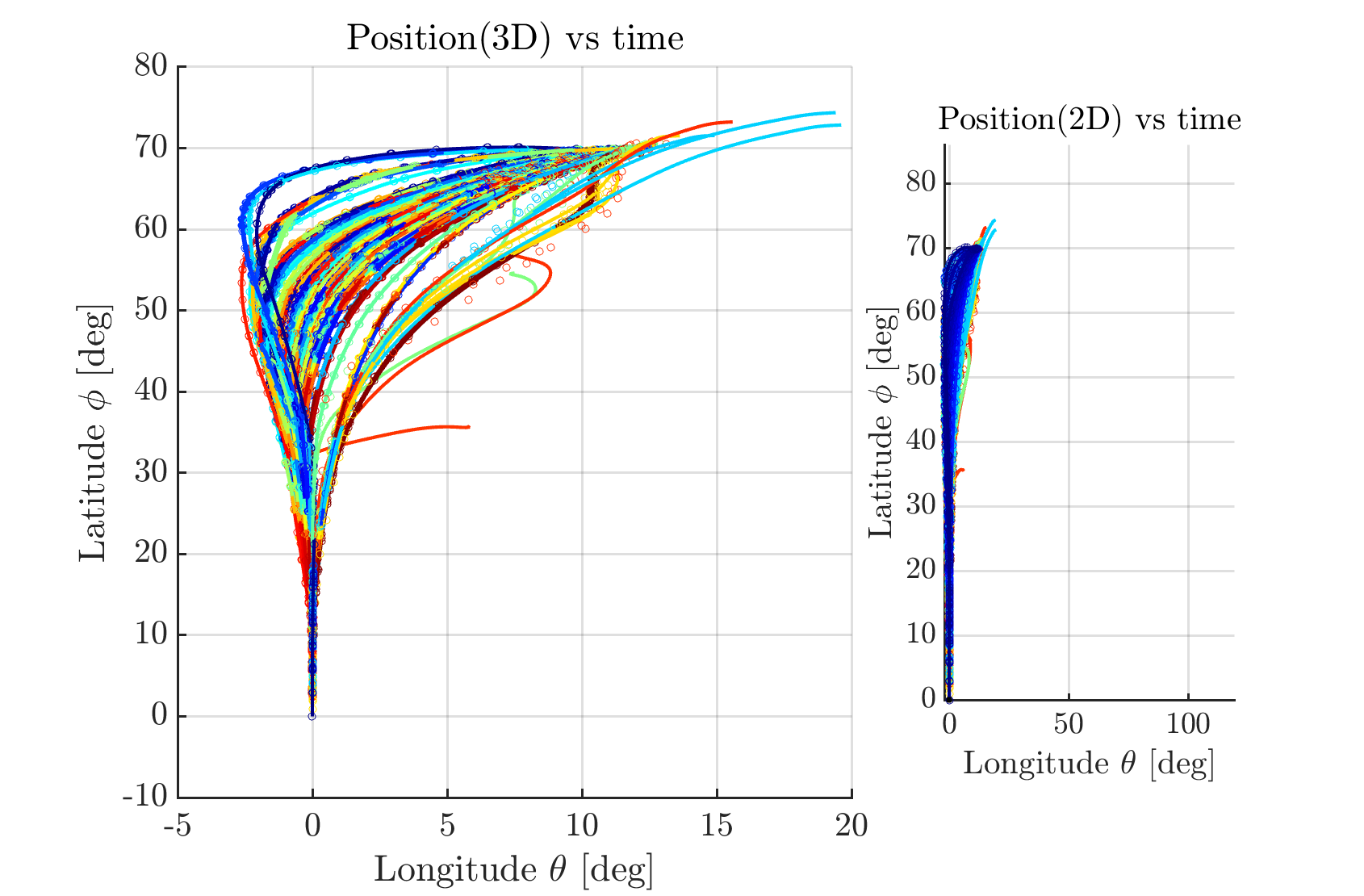}
    \end{minipage}%
    \begin{minipage}{.5\textwidth}
    \centering
    \includegraphics[width=1\columnwidth,trim={0cm 0cm 12cm 2cm},clip]{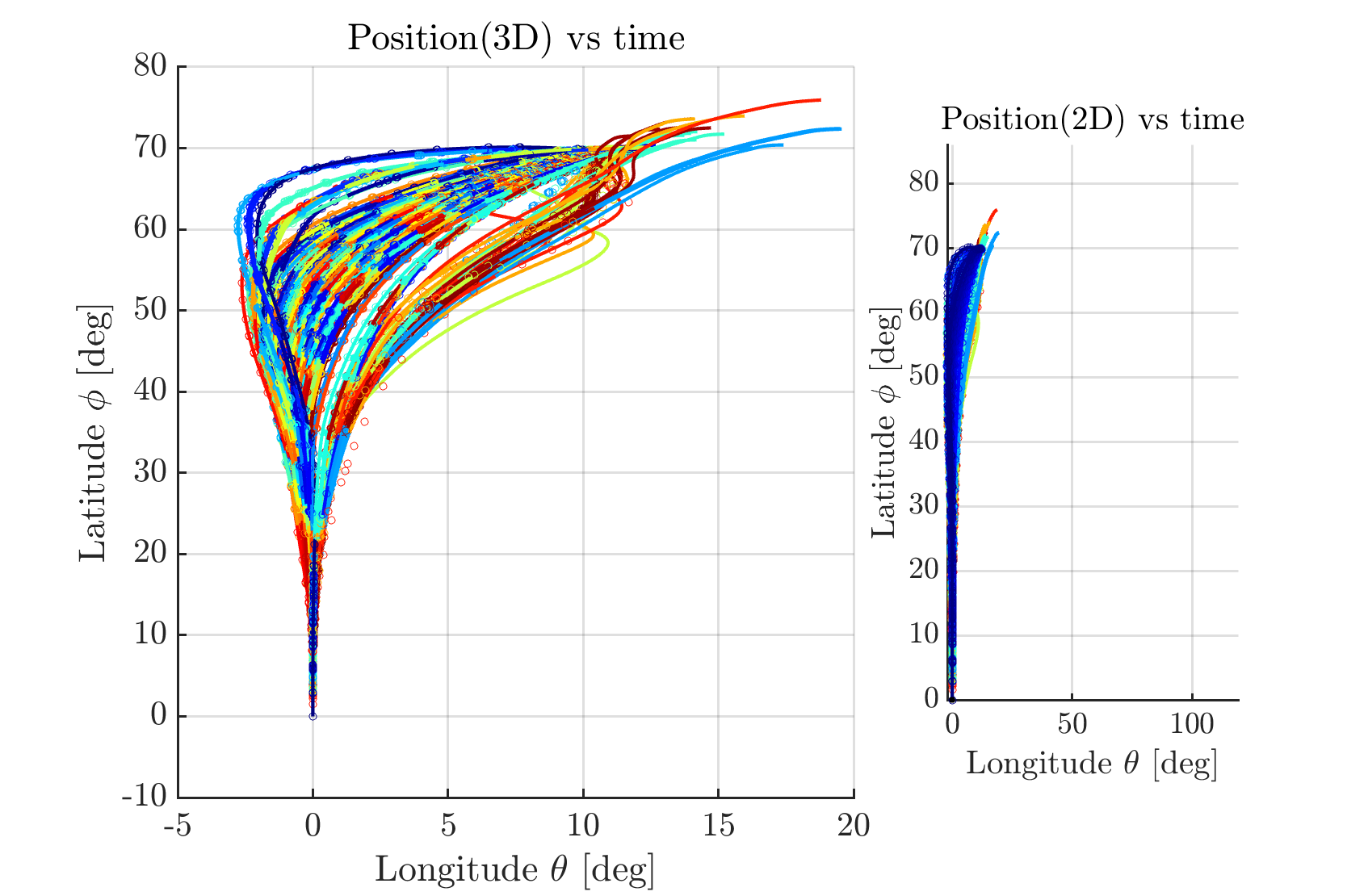}
    \end{minipage}%
    \caption{Left: Low-density uniform dispersion study of problem parameters over a coarse grid (216 cases). 
    Right: High-density uniform dispersion study over a fine grid (490 of 1230 cases displayed).
    }
    \label{fig:montecarlo}
\end{figure}

Assuming bank angle as the only control input, two more comprehensive studies were conducted where problem parameters were dispersed uniformly between parameter bounds. The parameters dispersed were initial altitude, initial planet-relative velocity, initial flight path angle, and vehicle mass. The bounds are given in Table \ref{tab:dispersions}. All other parameters and hyperparameters are given as in Tables \ref{tab:ex1-mission-params}, \ref{tab:ex1-mission-conv-hyperparams} and \ref{tab:ex1-mission-ptr-weights}. NFZ constraints were omitted. First, a coarse dispersion was done, with 216 cases. This was followed up with a separate dense dispersion, with 1230 cases. Altogether, 1446 unique cases were run. Mean statistics for performance of the {\autoscvx} algorithm are given in Table \ref{tab:montecarlo}. For the coarse dispersion, $93.5 \%$ of the 216 cases converged to a feasible, locally optimal solution. For the fine dispersion, $92.9 \%$ of the 1230 cases converged to a solution. Both the $6.5 \%$ of cases that did not meet criteria for dispersion in the coarse dispersion, and the $7.1 \%$ of cases in the fine dispersion, were ruled out based on hitting a maximum algorithm iteration limit without satisfying the convergence criteria; when hitting 20 SCP iterations, those runs were terminated. Note that this is not a particularly large number of iterations, and there remains the possibility that increasing the limit may result in a solution. For nonconvex problems, it is unclear whether infeasibility is due to pre-mature termination of the algorithm because of the chosen limit for maximum number of SCP iterations (i.e. 20) that doesn't allow convergence within specified tolerance, or because of having a dynamically infeasible problem instance. Currently, there is not a general way to certify feasibility for a general nonconvex problem such as hypersonic reentry. Future work will involve investigating certificates of infeasibility in cases where a converged solution is not achieved. These uniform cases, whether or not they met the convergence criteria prior to hitting the maximum iteration limit, are displayed in Figure \ref{fig:montecarlo}.

\begin{table}[htb!]
\centering
\small
\begin{tabular}{lrrrrr}
\hline
    \textbf{Dispersed Problem Parameter:} &
    \textbf{Nominal Value:} &
    \textbf{Uniform Perturbation Range:} \\ \hline
\textbf{Initial orbital altitude $h_{I}$} & $100,000 ~\text{m}$        & $[-10,000,10,000] ~\text{m}$       \\ \hline
\textbf{Initial orbital velocity $v_{I}$} & $7,450 ~\text{m/s}$        & $[-50,150] ~\text{m/s}$       \\ \hline
\textbf{Initial orbital flight path angle $\gamma_{I}$} & $-0.5^\circ$     & $[-0.3,0.4]^\circ$       \\ \hline
\textbf{Vehicle mass $m$} & $104,305 ~\text{kg}$ & $[-1000,1000] ~\text{kg}$       \\ \hline
\end{tabular}
  \caption{Uniformly dispersed problem parameters for Monte Carlo analysis of the \texttt{\autoscvx} algorithm.}
  \label{tab:dispersions}
\end{table}

\begin{table}[htb!]
\centering
\begin{tabular}{lrrrrr}
\hline
Mean Statistics: & \multicolumn{2}{c}{\# Dispersed Cases:} \\\cmidrule{2-3}
    & 216
    & 1230 \\ \hline
\textbf{\# Iterations} 
    & 9.7    
    & 9.37       \\ \hline
\textbf{Solve time per iteration {[}ms{]}} 
    & 131.0.5 
    & 129.09  \\ \hline
\textbf{\begin{tabular}[c]{@{}l@{}} Cost [m/s] \\ (Terminal velocity) \end{tabular}} 
    &325.13 
    &320.12  \\ \hline
\textbf{\begin{tabular}[c]{@{}l@{}}Non-dimensionalized \\ constraint violation residual\end{tabular}} 
    & 0.047 
    & 0.041 \\ \hline
\textbf{Convergence \% (across all runs)} 
    & 93.5 \%  
    & 92.9 \%     \\ \hline
\end{tabular}
  \caption{Mean performance metrics of \texttt{\autoscvx} in Monte Carlo analysis with coarse and fine dispersions.}
  \label{tab:montecarlo}
\end{table}

\section{Conclusion}

This work proposed {\autoscvx}, or \textit{auto-tuned primal-dual successive convexification}, as a framework for solving nonconvex optimal control problems such as the hypersonic reentry guidance. This method optimizes dual variables in closed-form within the SCP framework in order to update the penalty hyperparameters used in the primal variable update. A benefit of this method is that it is auto-tuning, and requires no hand-tuning by the user with respect to the constraint penalty weights. This method is motivated with duality theory, and after derivation of a closed-form solution for the penalty weight updates, the primal convex discrete subproblem for the hypersonic reentry problem is formed. These primal and dual updates are iterated until convergence. The full method is presented in Algorithm \ref{alg:autoscvx}.

Several example hypersonic reentry problems are posed and solved with the proposed approach, and comparative studies against the existing {\ptr} algorithm are conducted. In these studies, {\autoscvx} demonstrates an ability to reliably and flexibly solve a wide array of complex problems. When introducing the additional nonconvex control constraint where angle-of-attack is introduced as a control variable, the quality of the best-performing hand-tuned {\ptr} solution struggles to produce a high-quality solution. 
The angle-of-attack control oscillates wildly, and the dispersed trajectories veer aggressively across the earth's surface while adjusting vehicle heading to match the terminal condition.
{\autoscvx}, however, consistently produces less-aggressive maneuvers (even under dispersion) with smoother angle-of-attack profiles that remain mostly tight along the course of the trajectories.
In addition, even when {\ptr} converges to a solution qualitatively-similar to {\autoscvx}, the solutions from {\autoscvx} more consistently has an more optimal cost. Because penalty weights for {\autoscvx} only spike where constraints are violated or tight, this algorithm may adaptively become better conditioned as the SCP subproblems converge to the feasible region. Once feasibility is achieved, the {\autoscvx} can prioritize optimality, and the gradients driving down the true cost begin to dominate. 

When parameter dispersions are introduced for the {\autoscvx} algorithm, such as 216 or 1230 runs over a uniform parameter sweep, the algorithm still generates solutions between $92.9-93.5 \%$ of these cases. In the cases where convergence was not achieved, the algorithm was terminated at a maximum iteration limit before the convergence criteria was satisfied. In general, no clear techniques exist to determine whether or not a nonconvex problem has a feasible solution for a given set of parameters. Future work will involve investigating certificates of infeasibility in cases where a converged solution is not achieved. 
This work builds and improves on existing sequential convex programming methods for nonconvex problems such as hypersonic reentry guidance. Historically, SCP methods involve hand-tuning for each individual set of problem parameters. The ability to auto-tune the penalty hyperparameters eases the burden from the trajectory designer when applying {\autoscvx} to different problems. Extensions to penalty functions of different forms and integration with continuous-time constraint satisfaction techniques to reduce inter-sample violation are future directions.

\section*{Appendix}

The state-transition matrices (STMs) of the discrete LTV system are defined across each interval:
\begin{subequations}
\begin{align}
    \PhiA(k, \tau) &= \left. \frac{\partial x(k, \tau)}{\partial x_{k}}  \right|_{\Xk} , \\
    \PhiBm(k, \tau) &= \left. \frac{\partial x(k, \tau)}{\partial u_{k}} \right|_{\Xk} , \\
    \PhiBp(k, \tau) &= \left. \frac{\partial x(k, \tau)}{\partial u_{k+1}} \right|_{\Xk}, \\
    \PhiS(k, \tau) &= \left. \frac{\partial x(k, \tau)}{\partial T_{k}} \right|_{\Xk},
\end{align}
\label{eq:stm-defn}
\end{subequations}
$\forallintervals$ and $\tau \in [0,1]$, where $\Xk \triangleq (\xbar_{k}, \ubar_{k}, \ubar_{k+1}, \Tbar_{k})$ is a nominal reference trajectory. We define an approximation of Equations \ref{eq:xfin} and \ref{eq:foh}:
\begin{subequations}
    \begin{align}
    \xbar(k, \tau) \triangleq  
    \left. x(k,\tau) 
   \right|_{\Xk} 
    &= \xbar_k
            + \int_{0}^{\tau} \Tbar_{k} f(x(k, s), \ubar(k, s)) ds  
            \label{eq:xfin-bar} \\ 
    \ubar(k, \tau) \triangleq 
    \left. u(k,\tau) \right|_{\Xk} &= 
            ( 1-\tau) \ubar_{k}+\tau \ubar_{k+1}, \label{eq:foh-bar}
\end{align}
\end{subequations}
$\forall \tau \in[0,1]$ and $\forallnodes$. Differentiating Equations \ref{eq:stm-defn} yields the LTV system of equations \cite{lin2014surveyopt}:
\begin{subequations}
\begin{align}
\dotc{\Phi}_{A}(k, \tau) 
&= \eval{ \frac{\partial \dotc{x}(k, \tau)}{\partial x_{k}} }{\Xk}
= \eval {\frac{\partial F\left(x(k, \tau), u(k, \tau), T_{k}\right)}{\partial x_{k}} }{\Xk}
= \Tbar_{k} \evalfull{ \frac{\partial f(x(k, \tau), u(k, \tau))}{\partial x(k, \tau)} \frac{\partial x(k, \tau)}{\partial x_{k}} }{\Xk}, \\
\begin{split}
    \\
    \dotc{\Phi}_{B}^{-}(k, \tau) 
    &= \eval{ \frac{\partial \dotc{x}(k, \tau)}{\partial u_{k}}}{\Xk} 
    = \eval{ \frac{\partial F(x(k, \tau), u(k, \tau), T_{k})}{\partial u_{k}} }{\Xk} \\
    &= \Tbar_{k} \evalfull{ \frac{\partial f(x(k, \tau), u(k, \tau))}{\partial x(k, \tau)} \frac{\partial \xbar(k, \tau)}{\partial u_{k}} }{\Xk}
    + \Tbar_{k} \evalfull{ \frac{\partial f(x(k, \tau), u(k, \tau))}{\partial u(k, \tau)}  \frac{\partial \ubar(k, \tau)}{\partial u_{k}} }{\Xk} , 
\end{split}
    \\
\begin{split}
    \\
    \dotc{\Phi}_{B}^{+}(k, \tau) 
    &= \eval{ \frac{\partial \dotc{x}(k, \tau)}{\partial u_{k+1}} }{\Xk}
    = \eval{ \frac{\partial F\left(x(k, \tau), u(k, \tau),T_{k}\right)}{\partial u_{k+1}} }{\Xk} \\
    &= \Tbar_{k} \evalfull{ \frac{\partial f(x(k, \tau), u(k, \tau))}{\partial x(k, \tau)} \frac{\partial x(k, \tau)}{\partial u_{k+1}} }{\Xk}
    + \Tbar_{k} \evalfull{ \frac{\partial f(x(k, \tau), u(k, \tau))}{\partial u(k, \tau)} \frac{\partial \ubar(k, \tau)}{\partial u_{k+1}} }{\Xk}, 
\end{split}
    \\
\begin{split}
    \\ 
    \dotc{\Phi}_{S}(k, \tau) 
    &= \eval{ \frac{\partial \dotc{x}(k, \tau)}{\partial T_{k}}}{\Xk}
    = \eval{ \frac{\partial F\left(x(k, \tau), u(k, \tau), T_{k}\right)}{\partial T_{k}} }{\Xk}
    = \eval{ \frac{\partial(T_{k} f(x(k, \tau), u(k, \tau)))}{d T_{k}} }{\Xk}\\
    & = \Tbar_{k} \evalfull { \frac{\partial f(x(k, \tau), u(k, \tau))}{\partial x(k, \tau)} \frac{\partial x(k, \tau)}{\partial T_{k}} + f(x(k, \tau), u(k, \tau)) }{\Xk} ,
\end{split}
\end{align}
\end{subequations}
which simplifies into the differential LTV system:
\begin{subequations}
\begin{align}
    \dotc{\Phi}_{A}(k, \tau) &= \Tbar_{k} A(k, \tau) \PhiA(k, \tau), \\ 
    \dotc{\Phi}_{B}^{-}(k, \tau) &= \Tbar_{k} A(k, \tau) \bar{\Phi}_{B}^{-}(k, \tau)+\Tbar_{k} B(k, \tau)(1-\tau), \\
    \dotc{\Phi}_{B}^{+}(k, \tau) &= \Tbar_{k} A(k, \tau) \PhiBp(k, \tau)+ \Tbar_{k} B(k, \tau) \tau , \\
    \dotc{\Phi}_{S}(k, \tau) &= \Tbar_{k} A(k, \tau) \PhiS(k, \tau)+S(k, \tau).
\end{align}
\label{eq:stm-ltv}
\end{subequations}
where we define the Jacobians of the original continuous nonlinear system from Equation \ref{eq:gen-ncvx-dyn} with respect to the state and control:
\begin{subequations}
\begin{align}
    A(k,\tau) 
    &\triangleq \eval{ \frac{\partial f(x(k,\tau),u(k,\tau))}{\partial x(k,\tau)} }{\Xk}, \\ 
    B(k,\tau) 
    &\triangleq \eval{ \frac{\partial f(x(k,\tau),u(k,\tau))}{\partial u(k,\tau)} }{\Xk}, 
\end{align}
\label{eq:ltv-jacs}
\end{subequations}
as well as the Jacobian of the continuous nonlinear system reformulated as a function of normalized time from Equation \ref{eq:nonlin-sys-tau} with respect to timestep horizon:
\begin{align}
     S(k,\tau) 
     &\triangleq \eval{ \frac{\partial F(x(k,\tau),u(k,\tau),T_{k}) }{\partial T_{k}} }{\Xk}
     = \eval{ f(x(k,\tau),u(k,\tau)) }{\Xk}.
\end{align}

We solve the initial value problem for Equations \ref{eq:stm-ltv} assuming initial boundary conditions:
\begin{subequations}
\begin{align}
    \PhiA(k,0)  &= \diag{\one{\nx}}, \\
    \PhiBm(k,0) &= \diag{\zero{\nnu}}, \\
    \PhiBp(k,0) &= \diag{\zero{\nnu}}, \\
    \PhiS(k,0)  &= \diag{\zero{\nx}},
\end{align}
\label{eq:ltv-stm-bc}
\end{subequations}
for each interval by taking the integral:
\begin{subequations}
\begin{align}
    \PhiA(k, \tau) &= I+\int_{0}^{\tau} \Tbar_{k} A(k, s) \PhiA(k,s) ds, \\
    \PhiBm(k, \tau) &= \int_{0}^{\tau} \Tbar_{k} A(k, s) \PhiBm(k, s)+\Tbar_{k} B(k, s)(1-s) ds, \\
    \PhiBp(k, \tau) &= \int_{0}^{\tau} \Tbar_{k} A(k, s) \PhiBp(k, s)+\Tbar_{k} B(k, s) ~s ~ds, \\
    \PhiS(k, \tau) &= \int_{0}^{\tau}  \Tbar_{k} A(k, s) \PhiS(k, s)+S(k, s) ds.
\end{align}
\label{eq:ltv-stm-integral}
\end{subequations}
where $\tau \in [0,1]$. The solution to this initial value problem across the interval yields the discrete LTV matrices to approximate the system dynamics:
\begin{subequations}
\begin{align}
    \Ak       &= \PhiA(k, 1), \\
    \Bmk   &= \PhiBm(k, 1), \\
    \Bpk   &= \PhiBp(k, 1), \\
    \Sk       &= \PhiS(k, 1),
\end{align}
\label{eq:ltv-discrete-ltv}
\end{subequations}
$\forallintervals$.

\section*{Acknowledgments}
Support for studying the convergence properties of the successive convexification framework was provided by the Office of Naval Research grants N00014-20-1-2288. This research was supported by NASA grant NNX17AH02A and was partially carried out at the NASA Johnson Space Center; Government sponsorship is acknowledged. A special thanks to Abhinav Kamath and Samuel Buckner for their work in preparing the conference version of the manuscript upon which this work was based.

\clearpage

\bibliography{references}

\end{document}